\documentclass[11pt,a4paper]{article}
\usepackage{amsmath,empheq,graphicx,float,bbold,tikz,enumitem,amsmath,dsfont,bbm}
 \usepackage{caption,bm,amssymb,calc,commath,lipsum,etoolbox,hyperref}
\newtheorem{theorem}{Theorem}
\newtheorem{lemma}{Lemma}
 
\newtheorem{corollary}{Corollary}
 \newtheorem{proposition}{Proposition}
\newcommand{\eop}{\hfill{$\blacksquare$}}
\usepackage[square,sort,comma,numbers]{natbib}
\newcommand{\Friends}{{\mathbb F}}
\newcommand{\bpgf}{{\bf { \overline f}}} \newcommand{\bpgftwo}{{\bf {\bar g}}}
\newcommand{\CPj}{CP-$j$ }  
\newcommand{\idot}{{\bm \cdot}}
\usepackage[a4paper,bindingoffset=0.6in,%
            left=0.7 in,right=0.7in,top=0.7in,bottom=0.7in,%
            footskip=.95in]{geometry}
  \newcommand{\ignore}[1]{}
\newcommand{\sxi }{{\bar \xi}}  
  
\begin{document}
\date{}
\title {A Viral Timeline Branching
Process to Study a Social Network }

\author{Ranbir Dhounchak  \\ IEOR, IIT Bombay, India  \\ ranbirsingh@iitb.ac.in    \and
        Veeraruna Kavitha \\ IEOR, IIT Bombay, India \\ vkavitha@iitb.ac.in  \and
        Eitan Altman \\  INRIA, France \\ eitan.altman@sophia.inria.fr
} 
\maketitle 
\textbf{Table of contents}
\begin{itemize}
\item Part-I -  Viral Marketing Branching Processes in OSNs
\item Part-II - Competitive Viral Marketing Branching Processes in OSNs
\end{itemize}

\newpage
\vspace*{10cm}
{\Large {\textbf{Part-I: Viral Marketing Branching Processes in OSNs}}}
\newpage
\centerline{\Large \bf 
Viral Marketing Branching Processes in OSNs}
\begin{abstract}
We consider the inherent timeline structure of the appearance of contents in online social networks (OSNs), while studying content propagation.  We model the propagation of a post/content of interest by an appropriate multi-type branching process. The branching process allows one to predict the emergence of global macro properties (e.g., the spread of a post in the network) from the laws and parameters that determine local interactions. The local interactions largely depend upon the timeline (an inverse stack capable of holding many posts and one dedicated to each user) structure and the number of friends (i.e., connections) of users, etc. We explore the use of multi-type branching processes to analyze the viral properties of the post, e.g., to derive the expected number of shares, the probability of virality of the content, etc.

In OSNs the new posts push down the existing contents in timelines, which can greatly influence content propagation; our analysis considers this influence. We find that one leads to draw erroneous conclusions when the timeline (TL) structure is ignored: a) for instance, even less attractive posts are shown to get viral; b) ignoring TL structure also indicates erroneous growth rates. More importantly, one cannot capture some interesting paradigm shifts/phase transitions, e.g., virality chances are not monotone with network activity parameter, as shown by analysis including TL influence.

%

In the last part, we integrate the online auctions into our viral marketing model. We study the optimization problem considering real-time bidding.   We again compared the study with and without considering the TL structure for varying activity levels of the network. We find that the analysis without TL structure fails to capture the relevant phase transitions, thereby making` the study incomplete.
\end{abstract}

\section{Introduction}
\label{intro}
The advent of the Internet has transformed the advertising industry in various ways. With the constant  year-on-year growth of the number of users, the global userbase of the Internet passed 3.5 billion mark in 2017,  constituting nearly half of the earth's population \cite{1IntStat}. This has made the Internet a powerful tool for organizations to interact with users and advertise their products/services in a personalized manner. In particular, Online Social Networks (OSNs) such as Facebook, Twitter, YouTube, etc play an instrumental role in the overall digital advertising of the products/services of various organizations. Users on these OSNs keep exchanging volumes of information/data in the form of images, blogs, texts, videos, etc. 
Due to   immense activities of the users in OSNs, the marketing/advertising companies promote their commercial contents by leveraging the strengths of these OSNs.  


In \textit{viral marketing,} the content providers (CPs)/advertisers create contents that are appealing to the users (e.g., giving offers, discounts, advertising in attractive manner). When  users find the service/product good enough, they involuntarily spread a word about it, triggering word-of-mouth.  Users share the content with their friends, and the information is thus spread through OSNs. 
In the abstract sense, information spreads like a virus from one person to another, and hence called viral marketing. However, the content propagation has additional complexities which must be incorporated in the model to accurately investigate the process/phenomenon. And  we study the same  in this paper.

\subsection{Motivation and scope of research: Timeline structure}
OSNs store volumes of information consumed by the users. At the user level,  
 these segments of information  (called posts) are  arranged based on their chronological order for the display to users \cite{1TLLit}. In other words, posts appear at different levels (each level holds one post) based on their newness  on each user's page in an OSN, for instance, News Feed in Facebook. We call this reverse chronological appearance of the posts a `timeline' (TL), one dedicated for each user. 
This order of storing contents on timelines (TLs)  and the related dynamics have great influence on content propagation.
   \textit{However, no attention is paid to the TL structure of the posts/contents appearing on a user's page in viral marketing literature.} We study the content propagation phenomenon over OSNs, considering the inherent TL structure.

 \begin{figure}[ht]
  \vspace*{-35mm} 
 \begin{minipage}{16cm}
\includegraphics[width = 16cm, height = 16cm]{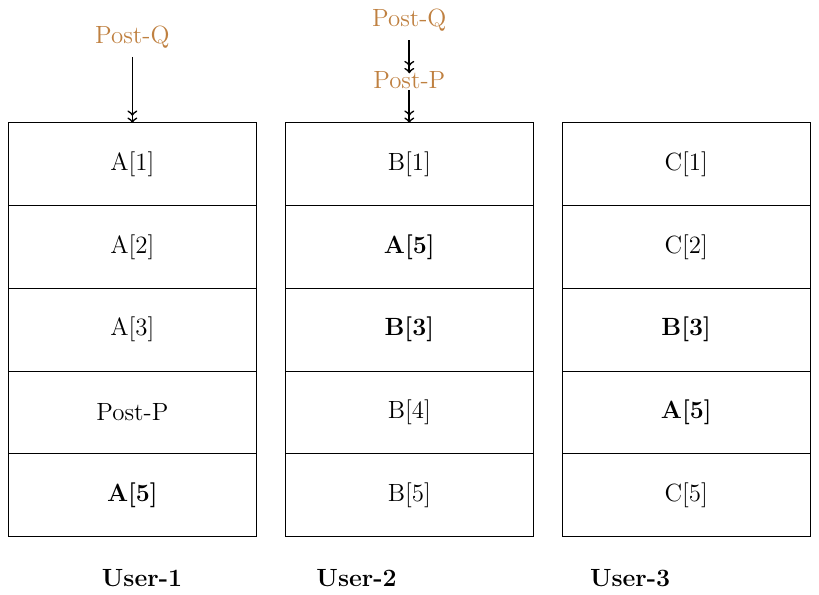} 
\vspace*{-5.5cm} 
\caption{TL structure \label{1Fig_TLIntro}} 
 \end{minipage}
\end{figure} 
 A typical example of TL structure   (for three users) is shown in  Figure \ref{1Fig_TLIntro}.
Referring to the Figure \ref{1Fig_TLIntro}, the natural question  to ask is:  how many users have a particular post of interest?, and the next immediate one is: at what level does that post reside (i.e., the position)? For instance, all the three users have the post `A[5]' on their TLs, but at different levels. It is clear  that the posts positioned on the top of the TLs receive more attention/visibility compared to the ones at lower levels. Further, the arrival of new contents keeps shifting/pushing down the existing contents of  a TL.
Consequently, a particular content of interest may reach lower levels before the user visits its TL, and  may miss the user's attention. Thus, the content of interest can potentially be missed by such shift transitions. Technically, a user can scroll through indefinite number of posts at various levels. However, it is known that users' attention is limited to the first few levels \cite{1Scroll}. 
We consider this aspect while analysing content propagation in the viral marketing scenario. We observe (theoretically as well as numerically) that this aspect makes a huge difference in the conclusions related to such a study.   In addition,  the TL structure also influences content propagation because of the following other aspects:  a) multiple posts with the same user (the attention gets divided);  and b) decreasing interests towards reading contents at lower levels of TLs etc.   We find that these key aspects must be incorporated into the model to accurately study the content propagation. We discover that  without these key elements into consideration, one leads to draw erroneous conclusions. 

Further when content of competing content providers circulates through the same social network, at the same point of time, there would be lot more influences on the propagation of the  content.
These influences are more complicated  with  timeline structure; for example, user might neglect the content of a low influential content provider when it has  (simultaneously) competing content.  Or alternatively one might be interested in forwarding all, or the user might be interested only in forwarding the post that appears first.  It might be possible that content of a content provider gets viral, but not that of the others etc.  These aspects require study of   decomposable   branching processes  and the same along with the propagation of the competing content is considered in Part-II \cite{part2} of this work.

\subsubsection*{Our approach and contribution: 
 Branching processes and viral marketing}
Content propagation over OSNs follows a number of models based on factors such as empirical evidences, the structure of an OSN, etc.
There has been an extensive literature on the content propagation over the OSNs, and an important approach for modeling the dynamics of content propagation has been the branching processes (see \cite{1BranchVM,1YU}, etc). Authors in \cite{1BranchNonMarkov} studied information diffusion in the real viral marketing campaigns (involving 31000 individuals) and showed that the branching processes explain the dynamics of information diffusion.  The branching processes (BPs) are adequate to incorporate the characteristics of content diffusion (e.g., phase transition-epidemic threshold) and provide explicit expressions for many important performance measures.   As an example, authors in \cite{1BranchVM} provided a discrete time branching model to predict the spread of a campaign. They estimated, using the theory of branching processes, the campaign's performance via various measures such as the number of forwarded e-mails, the number of viral e-mails, etc. as a function of system parameters. 
Other studies (e.g., \cite{1VirBranch,RumorSpread1}),  reinforce that the branching process can well fit the content propagation trajectories collected from real data.

In a branching process, a parent produces identically and independently distributed (IID) offsprings. When one models content propagation over an OSN as a branching process, any parent should produce identical number of offspring  (identical to the other parents) and independent of the offspring produced previously. Moreover, the parents keep producing offspring even when the population explodes. This is possible only when the OSN has infinite population. 
One can assume that the OSNs with huge user base have infinite population (unbounded number of users). Further, when users have identically distributed number of friends, then the BPs can model the content propagation over OSNs. This simplifies modeling and analysis. We use \underline{Multi-Type} Branching Processes (MTBPs) (e.g., \cite{1AthreyaBook,1AthreyaPaper}) to  model  the influence of TL structure on  content propagation in OSNs.  We further extract the realistic features of content sharing in a typical OSN and incorporate them appropriately in our model. The branching processes can mimic most of the phenomenon that influences the content propagation. For example, one can model the effects of multiple posts being forwarded to the same friend,  and multiple forwards of the same post, etc.

A post on a higher level in TL has better chances of being read by the user. Posts of appealing nature, e.g., containing irresistible offers, have a great chance of being in circulation, and we call it the post quality factor. Posts of similar nature appearing at lower levels on the TL have smaller chances of appreciation, etc.  To study all these factors, one needs to differentiate the TLs that have the `post'  at different levels, and this is possible only through multitype BPs. The following are some elements of our approach:
\begin{itemize}
\item Using the well-known results of multi-type branching processes (MTBPs), we obtain closed-form expressions for some performance measures which provide insights for the performance of the campaign, e.g., visibility of the content, virality\footnote{A post is said to be viral if a large number of posts are already shared and if the number of users with posts is exploding with time.}, etc.  For the other measures, we either have approximate (time asymptotic) closed-form expressions or simple fixed point equations whose solution provides the required measures. 
\item  We study the influence of various network (structural) parameters on the content propagation. We also study the effects of posts sliding down the  TLs due to network activity.
As the  mean number of friends (network activity) increases, one can expect contents to spread more rapidly (monotonous behavior). Contrary to that, we discover non-monotonous behavior in the virality of a content. This phenomenon is fundamentally due to TL structure.

\end{itemize} 
  

In the last part, we integrate online auctions into our viral marketing model. We
study an optimization problem considering real-time bidding.  We compare the study considering the  TL structure to that without considering the TL structure for varying activity levels of the network. Our observations
are   similar: drastically different  optimizers and the analysis without
  TL structure fails to capture
the relevant phase transitions.

\subsection{Related work}
The huge growth in the activity on the Internet has generated wide interest in understanding content propagation on the Internet. Previously, peer to peer (P2P) networks (e.g., \cite{1xu,1YU})
 have played an important role. While of late there has been a lot of interest in content propagation over OSNs (e.g., 
 \cite{1BranchVM,1BranchNonMarkov,1Efficient,RumorSpread1}, etc). The P2P networks pull the required information from their peers, while in viral marketing the information is pushed for marketing purposes.
 In viral marketing, one needs to keep pushing information by passing it on to seed nodes to keep the flow going on.  We consider content propagation over OSNs like Facebook, Twitter, etc where the information (called post/message) is again pushed, but involuntarily.
 Here the post/content is forwarded to few initial seeds, and the post gets viral  based on the interest generated among the users and the extensive sharing.

There is a vast literature that studies the propagation of content 
over OSNs.  Many models discretize the time and study content propagation across the discrete time slots  (e.g., \cite{1Efficient}). As argued in \cite{1CPonGraph} and references therein, a continuous time version (events occurs at continuously distributed random time instances)  is a  better model  (\cite{1BranchVM,1BranchNonMarkov} ) and we consider the same. In majority of the works which primarily use graph-theoretic models, 
 the information is spread at maximum to one user at any message forward event (e.g., \cite{1Efficient,1CPonGraph,RumorSpread1} etc.).  However, when a user visits a OSN (e.g., Facebook, Twitter, etc.),  it typically forwards multiple posts and typically (each post) to multiple friends. 
 Authors in (\cite{1BranchVM,1BranchNonMarkov}, etc) study viral marketing problem,  where the marketing message is pushed    continuously via emails,  banner advertisements,
 or search engines, etc. This scenario allows multiple forwards of the same post, and is analysed using BPs. 
However, they do not consider the influence of other posts using the same medium,  and the other effects of TLs.   
As already mentioned, these aspects majorly influence the analysis. 
 Branching processes have been used in analysing various types of networks,
such as,  polling systems (\cite{1R1}) which have been used to model
local area networks and  P2P networks (\cite{1R2}), etc. We use branching processes
not only to study the time evolution of the contents of interest (extinction and viral growth)
but also to provide a spatio-temporal description of the process.
We model the evolution of the number of timelines that have a given 
content at a given level of the timeline (e.g., top of the timeline).

\section{System description} 
We consider a giant OSN, e.g., Facebook, Twitter, VKontakte, etc. In this paper, we track a content of interest corresponding to one specific content provider. In the second part of this work, we consider the contents of  multiple (competing) content providers\footnote{The extension of this work is submitted separately for details see Part-II \cite{part2} }. Users use these networks to connect to other users to share photos, news, events/activities taking place around them, commercial contents, etc. We briefly refer to these pieces of information as a post. Recall that these posts appear at different levels (on the screen) based on their newness, for instance, News Feed on Facebook. When a user visits\footnote{The users 'visit' OSNs at random intervals of time and in each `visit' it browses some/all new posts.} the OSN, it reads the posts on its timeline and shares a post, upon finding it appealing/useful, with some of its friends (users connected to him). In this sharing process, the post appears on the top level of the timelines of those friends with whom the post was shared. This brings about a change in the appearance of contents on the timelines  of recipients of the post. Basically, the existing contents of these TLs shift one level down each. And a user can share as many posts as it wants. The number of shares of a particular post by a particular user depends upon: a) the distribution of its number of friends;  and b) the extent to which the user liked the post.
And extensive sharing of the post amongst the users potentially makes the post viral. It is evident that the sharing of a post depends on how engaging the \textit{content provider} (CP) designs its post.  There are some more aspects which influence the content propagation. For example, users may become reluctant to read/share the contents on the lower levels of their TLs. When they see multiple posts of similar nature, they may appreciate few posts while the remaining ones receive reduced attention. 
 We study all those aspects and the dynamics created by the actions (e.g., like, share, etc)  of the users, which  have a major impact on the propagation of the commercial content.
  
\subsection*{Continuous time branching processes}
  The continuous time branching processes (CTBPs)  are often good candidates for modeling viral marketing models. We describe these processes briefly as follows. 
  Let  $X(0)$ be the number of initial particles in a CTBP. Each of these particles stays alive for an exponentially distributed time with parameter say $\lambda$ and then dies.  The `death' times of the particles are independent of the others. Hence, the first death occurs after  exponentially distributed time with parameter $X(0)\lambda$. Upon its death, it produces a random number (say $\zeta$) of offspring, which join the existing population. The number of particles immediately after the death of the first particle changes to $X(0) - 1 + \zeta$. The `death' times are again exponentially distributed  (by the memoryless property), and the process continues. It is well-known that (under certain assumptions) the BPs have certain dichotomy:  a) either the population gets extinct  (death of last individual TL), or b) the population grows exponentially fast with time;  and c) there is no third way.

 \subsection{Dynamics of  content propagation and branching process}
 The content propagation in a typical OSN is as follows. Let us say we are interested in the propagation of post-\textbf{P} when the process starts with $X(0)$ number of seed TLs. We track the post-\textbf{P} till first $N$ levels of TLs. It is important to note that $X(0)$ remain unread before their respective users become aware of the contents on the TLs, i.e., before they visit their TLs. We call this TLs as \textit{ number of unread TLs }(NU-TLs). If a user, among  $X(0)$, visiting its TL finds  post-\textbf{P} attractive, it reads the post and may share the same with a random number of its friends.  And post-\textbf{P} would be placed on the top level of the recipient TLs. As shown in  Figure \ref{1FigShift}, the recipient TL has  post-\textbf{P} on the top, and remaining posts shift down one level each.
 \begin{figure}[h]
 \begin{center}
\begin{tikzpicture}[scale=1]
\def \n {5}
\draw[->>, thick] (-2.5,6.8)--(-2.5,6) node[pos = -0.3]  { \textcolor{black}{\text{post-\textbf{P}}}};
\foreach \y/\x in {1/\textbf{A[5]},2/post-Q,3/A[3],4/A[2],5/A[1],}
	{ \draw (-4,\y) rectangle (-1,\y+1) node[pos=.5] {\x}; 
	} 
	\draw[->>, red,thick]   (-1,3.5) -- (0.5,3.5);
\end{tikzpicture}
\begin{tikzpicture}[scale=1]
\def \n {5}
\foreach \y/\x in {1/post-Q,2/A[3],3/A[2],4/A[1],5/post-\textbf{P}}
	{ \draw (-4,\y) rectangle (-1,\y+1) node[pos=.5] {\x}; 
	} 
\end{tikzpicture}
 \end{center}
\hspace*{-3cm}
\caption{  Shifting of the contents on a TL} \label{1FigShift}
\end{figure}

If some more posts are shared again with some of these recipient TLs, the contents further shift down of the corresponding TLs.  For instance, when one more post is shared after the post-\textbf{P} with the same shared user, the post-\textbf{P} resides on the second level of the corresponding TL.

We first argue that the continuous time version of the branching process fits the content propagation better than the discrete counterpart. In a CTBP,  any \underline {one} of the existing particles `dies'  after exponentially distributed time while in a discrete time version all the particles of a generation `die' together. When the number of copies of CP-post grows fast (i.e., when the post is viral), the time period between two subsequent changes decreases rapidly as time progresses. This is also well captured by CTBP, which mimics the content dynamics better.
 As the underlying OSN is huge, one can say that the visit times of users are virtually independent of the each other. We assume memory-less visit times, i.e., the users visit their TLs at intervals that are exponentially distributed as in a CTBP. The sharing process generates a random number, say $\zeta$, of new TLs holding post-\textbf{P}. If the user does not read or share the post after visiting its TL, then $\zeta = 0$. If sharing process is independent and identical across all the users, the new TLs $\zeta$ so generated resemble IID  offspring in a CTBP and the effective NU-TLs with post-${\bf P}$ may appear like the particles of a CTBP. 
 When one of the users of these NU-TLs (including the new ones)   visits its TL  and starts sharing the  post-\textbf{P} (as before),   then the content propagation   dynamics again resemble a CTBP. 
 
 However, the CTBP  described above does not capture some aspects related to post-propagation process.   Post-${\bf P}$ can disappear from some of the TLs, before the corresponding user's visits. To be precise, the post-${\bf P}$ would disappear from a TL with $(N-l+1)$ or more shares, if initially post-${\bf P}$ were at level $l$. For example, the post '\textbf{A[5]}' is lost by the arrival of post-\textbf{P}  in Figure \ref{1FigShift}. 
In all, the propagation of content in an OSN  is influenced by two factors:  a) the evolution of  TLs with post--${\bf P}$ when some other posts  are shared with them  (contents on the TL shift down); and b) the sharing dynamics of post--${\bf P}$ between different TLs.  
 

If we consider a CTBP  with a single type of population,  all the particles will have the same death rate and offspring distribution (e.g., \cite{1AthreyaBook}).  However, the disappearance of post-${\bf P}$ from a TL depends upon the level at which the post is available. 
Further, we will see that many more aspects of the dynamics depend upon the level at which the post-${\bf P}$ resides. 
Thus clearly, the single type CTBP is not sufficient, and we require a  multi-type  continuous time branching process (MTBP).  
An MTBP   describes the population dynamics in the scenarios with a finite number of population-types. All the particles belonging to one type have the same death rate and offspring distribution; however, these parameters could be different across different types.  
To model the rich behaviour of the propagation dynamics,  we will require (details in later sections) a particle of a certain type to produce offspring of other types. This modeling feature is readily available with MTBPs. We will show that the propagation dynamics can be well modeled by an appropriate MTBP, where for any $l \le N$, all the TLs with post-${\bf P}$ in level $l$ form one type of population. We use the following feature of the branching processes: it suffices to study the evolution of the population with one initial/seed particle. To be more specific, the analysis starting with multiple seeds can be derived using the analysis with one seed particle (details are in later sections).  
\\
\textbf{Assumptions:}
 We track the post of the CP and study the time evolution of the post over  TLs till first $N$ levels. We assume a TL with posts of the CPs is not written\footnote{We say a TL is written when a friend of it shares a post which changes its content.} with the post of the same CP again. In a huge social network, it reasonable to assume that the probability of the same post being shared again with any user is very small. Also, one can find applications that satisfy such assumptions.  As an example, consider few organizations that plan to advertise their products using a coupon system.   Also, consider that these coupons can be shared with friends.  But a user with one or two such coupons can not be shared with another coupon at a later point of time.  To avoid multiple shares to the same user, there is a control mechanism. Any user sharing the coupons with its friends, needs to declare the recipients in a list which disables the share of coupons to the same recipients.

\section{Single content provider model} 
\label{sec_single_cp_analysis}
We consider a single CP and refer to its post as the CP-post.   The TLs containing CP-post may have it  at any level from one to $N.$
These TLs also contain the other posts, and the movement of these posts can also affect the propagation of the CP-post. And our focus would be on CP-post. We say a user is of type $l$, if its TL contains the CP-post on level $l$  and the top $l-1$ levels do not contain the CP-post.
 Let $X_l(t)$ represent the number of \underline{unread} TLs (NU-TLs) of  type $l$  at time $t$.   We study the time evolution of $\{ X_l (t) \}_{l} $.
We will show below that the $N$-valued vector process $\textbf{X}(t) := \{X_1(t),X_2(t),\cdots,X_N(t)\}$ is an MTBP  under suitable conditions. 

\subsection{Modeling details}
\noindent{\bf Birth-death process via shift and share transitions:}  To model content propagation process by an appropriate  branching process, one needs to specify the `death' of
an existing parent  (a TL with `unread' CP-post in our case) and the distribution of its offspring. A user of type $l$ is said to `die' either when its TL is written by another user or when the user itself wakes up (visits its TL) and shares the post
with some of its friends. In the former event, exactly one user of type $(l+1)$ (if $l < N$) is `born' while the latter event gives birth to a random number of offspring of types $1$ or $2$ or $\cdots N$. 

 If {\it$i-1$ (with $i \le N$) posts are shared with the same user after the CP-post, then the CP-post is available on the $i$-th level and we will have a type $i$ offspring.} Assume  that a user produces offspring of type $i$ with probability $\rho_i$ and that
\underline{  $\rho_1 > 0.$}  Note that $\sum_i \rho_i=1$. In general, users have lethargy to view/read all the posts. We represent this via a level based reading probability, $r_l$, which represents the probability that a typical user reads the post on the level $l$. It is reasonable to assume \underline{$r_1 \ge r_2 \cdots \ge r_N.$} We have two types of transitions that modify the MTBP, which we call {\it shift and share transitions.} In the share transition, a user first reads the CP-post and based on the interest generated, it shares CP-post with a random number of friends. The Figure \ref{1Fig_TL_ShareCh} below describes the share transition.
 \begin{figure}[h]
\begin{tikzpicture}[scale =1.1]
\node[ rounded corners] at (.7,2.8){ TL};
\draw[ xscale = 0.5, yscale =0.5, thick, fill= blue!10,rounded corners] (0,0) rectangle (3,5)  node[pos=.5]  { CP-post };
\draw[ ->, very thick ] (1.5,1.25) to ++(1,0);
\foreach \x in {1,...,5}  \draw[xscale = 0.5, yscale =0.5] (0,\x) -- (3,\x ) (-1,5-\x +.5 ) node {\tiny{Level}-\x};
\draw[thick, fill= blue!10,rounded corners] (2.5,0) rectangle ++(2,2.5) node[align=center,pos=.5] { \footnotesize{Reads CP-post} \\  w.p.  $r_3$ \\ (at level-3)};
\draw[ ->,very thick ] (4.5,1.25) to ++(1,0);
\draw[thick, fill= blue!10,rounded corners] (5.5,0) rectangle ++(2,2.5) node[align=center,pos=.5] {Interest \\ generated \\  w.p. $\eta$};
\draw[ ->, very thick ] (7.5,1.25) to ++(1,0);
\draw[thick, fill= blue!10,rounded corners] (8.5,0) rectangle ++(1.6,2.5) node[align=center,pos=.5] { \normalsize{CP-post } \\ shared  \\ \normalsize{with $\zeta$} \\  \# friends};
\draw[ ->, very thick ] (10.1,1.25) to ++(.9,0);
\draw (10.5,1.4) node[align = right, below] {\small{$\rho_3$}};
\draw[thick, fill= blue!10,rounded corners] (11,0) rectangle ++(2,2.5)  node[pos=.5]  { {\footnotesize{ \#cp-post-3++}}  };
\draw (12,2.2) node {{\footnotesize{ \#cp-post-1++}}  } (12,1.7) node { } (12,.7) node {   } (12,.2) node { {\footnotesize{ \#cp-post-5++}}  };
\foreach \x in {1,...,5}  \draw[xscale = 0.5, yscale =0.5] (22.6,\x) -- (25.6,\x );
\path (9.95,1.1) node(x) {}  (11.1,2.5) node(a) {} (11.1,1.7) node(b) {}  (11.1,.75) node(d) {}  (11.1,0) node(e) {};
\draw [ ->, very thick ]  (x) -- (a)  node[near end,below] {\small{$\rho_1$}}; 
\draw [ ->, very thick ]  (x)   -- (e) node[near end,above] {\small{$\rho_5$}}; 
\end{tikzpicture}
\caption{Share transition}\label{1Fig_TL_ShareCh}
\end{figure}

 In the shift transition, user with CP-post is written by other users, and the position of CP-post shifts down.
 
\noindent{\bf  CP-post propagation dynamics:}  Let  $\mathcal{G}_1$ represent the subset of users with CP-post at some level, while  $\mathcal{G}_2$ contains the other users. 
We assume the OSN  (and hence $\mathcal{G}_2$) has infinitely many users and note  
 $\mathcal{G}_1$ at time $t$ has, 
 \begin{eqnarray}
X(t) &:=& \sum_{l\le N} X_l (t), \label{Eqn_tot_num}
 \end{eqnarray}
 number of users.  Group ${\cal G}_2$ has an infinite number of users/agents, and this remains the same irrespective of the size of ${\cal G}_1$, which is finite at any finite time. 
Thus, the transitions between ${\cal G}_2$ and  ${\cal G}_1$ are more significant, and one can neglect the transitions within  ${\cal G}_1$.  It is obvious that we are not interested in transitions within  ${\cal G}_2$ (users without CP-post). We thus model the action of these groups in  the following consolidated manner: 
\begin{itemize}
\item \textbf{ Share transition:} Any user from $ \mathcal{G}_1 $ wakes up after $exp(\nu)$ time (exponentially distributed with parameter $\nu$)  to visit its TL and writes to a random (IID) number of users of $\mathcal{G}_2$ (refer to Figure \ref{1FigShift}).
\item \textbf{ Shift transition:}  The TL of  any user of  ${\cal G}_1$ is written by  one of the users of ${\cal G}_2$, and the time intervals between two successive writes are exponentially distributed with  
parameter $\lambda$ (refer to Figure \ref{1Fig_TL_ShareCh}).
\end{itemize} 

The state of the network, ${\bf X} (t)$, changes when the first of the above-mentioned events occurs. At time $t$, we have $X(t)$ (see equation (\ref{Eqn_tot_num})) number of users in $\mathcal{G}_1$ and thus (first) one of them wakes up according to exponential distribution with
parameter $X(t)\nu$. 
Similarly, the first TL/user of the group ${\cal G}_1$ is written with a post after exponential time with parameter $X(t)\lambda$. 
Thus, the state ${\bf X} (t)$, changes    after exponential time with
parameter $X(t)\lambda + X(t)\nu$.
Thus, the rate of transitions at any time is proportional to $X(t)$, the number of NU-TLs at that time, and hence, the rate of transitions increase sharply as time progresses,  
when the post gets viral. Considering all the modeling aspects,  the   IID offspring generated by one $l$-type user are summarized as below (w.p. means with probability):
\begin{eqnarray}
\xi_l = 
\left\{
\begin{array}{llll}
{\bf e}_{l+1}\mathbb{1}_{l < N} & \text{ w.p.}\ \ \theta  := \frac{\lambda}{\lambda + \nu} \mbox{ and }  \\ 
\zeta {\bf  e}_i & \text{ w.p.} \ \  (1-\theta ) r_l\rho_i \hspace{0.4cm} \forall i \leq   N \\
0 & \text{ w.p. } \ \ (1-\theta)(1-r_l).
\end{array}
\right.   
\label{Eqn_offspring}
\end{eqnarray}
where $\textbf{e}_l$ represents standard unit vector of size $N$ with one in the $l$-th position, $\mathbb{1}_A$ represents the indicator, $\zeta$ is the random number of friends to whom the post is shared   and $r_l$ is the probability the user reads/views a post on level $l$. Figure \ref{1Fig_TL_SSCh} demonstrates the transitions.
 \begin{figure}[ht]
\begin{tikzpicture}[scale =1.4]
\node[ rounded corners] at (1.1,2.8){Timeline of a user};
\draw[ xscale = 0.5, yscale =0.5, thick, fill= blue!10,rounded corners] (0,0) rectangle (5,5)  node[pos=.5]  { CP-post (at level-3)};
\foreach \x in {1,...,5}  \draw[xscale = 0.5, yscale =0.5] (0,\x) -- (5,\x ) (-1,5-\x +.5 ) node {Level-\x};
\draw[->, very thick] (2.5,1.25)  to[out = 80,in = 140] (6,4.5);
\draw  (3,3) node[ above,sloped] { \rotatebox{55}  {Shift transition w.p. $\theta$}}; 
\draw[ xscale = 0.5, yscale =0.5, thick, fill= blue!10,rounded corners] (10,4) rectangle (15,9)  (12.5,5.4) node  { \footnotesize{CP-post now at level-4}};
\foreach \x in {1,...,5}  \draw[xscale = 0.5, yscale =0.5] (10,\x+3) -- (15,\x +3) (9,9-\x +.5 ) node {Level-\x};
\draw[->, very thick] (2.5,1.25) -- (6,-.5);
\draw    (4,0) node{ \rotatebox{-20} {Share transition w.p. $1-\theta$}};
\draw[thick, fill= blue!10,rounded corners] (6,-2) rectangle ++(2.5,2.5) node[align=center,pos=.5] {CP-post shared with  \\ $\zeta$ no. of friends \\ \\ \# CP-post++};
\end{tikzpicture}
\caption{Propagation of CP-post: transitions}\label{1Fig_TL_SSCh}
\end{figure}

 Recall that users (offspring) of type $i$ are produced with probability $\rho_i$ during the share transitions. From equation (\ref{Eqn_offspring}) the offspring distribution  is  identical at all 
  time  instances $t$, $\zeta$ can be assumed  independent across users, and  hence  $\xi_l$ are IID offspring from any type $l$ user.
Further, all the transitions occur after memoryless exponential times, and hence ${\bf X} (t)$ is an MTBP with  $N$- types  (e.g. \cite{1AthreyaPaper}). 
\\
\\
\noindent{\bf PGFs and post quality factor:}
Let $f_F( s, \beta)$ be the probability generating function (PGF) of 
 the number of friends, $\Friends$, of a  typical user, parametrized by   $\beta$. 
For example,$f_F(s, \beta) = \exp (\beta (s-1) ) $ stands for Poisson distributed $\Friends $,   $f_F(s,  \beta) = (1-\beta)/(1- \beta s) \ $ stands
for geometric  $\Friends $.  Let $m = f'_F (1, \beta)$ represent the corresponding mean.
A user shares the post with some/all of its friends ($\zeta$ of equation (\ref{Eqn_offspring})) based on how engaging the post is. Let the post quality factor $\eta $ quantify the extent of the CP-post engagement on a (continuous) scale of 0 to 1 where $\eta = 0$ means the worst and $\eta =1$ is the best quality. {\it We assume that the mean  of  the number of shares  is proportional to this quality factor.}  In other words,   $m (\eta) = m \eta  $ represents the post quality dependent mean of the random shares.  Let $f(s, \eta, \beta)$ represent the   PGF of $\zeta$.  For example, for Poisson  friends, the PGF and the expected value of $\zeta$  are given respectively by:

$$f (s, \eta, \beta) = f_F(s, \eta \beta) = \exp (\beta \eta (s-1) ) \mbox{ for any $s$ and } m (\eta) = \eta \beta. $$

For Geometric friends, one may assume the post quality dependent parameter 
$$ \beta_\eta =  {(1-\beta)}/{(1-\beta+\beta\eta)}, \mbox{  which ensures  } m(\eta) = \eta \beta .$$ And then the PGF of $\zeta$ is given by
$
f(s,\eta,\beta) = f_F(s, \beta_\eta) = {(1-\beta_\eta)}/{(1-\beta_\eta s)}.$
One can derive such PGFs for other distributions of $\Friends.$  \textit{Interestingly enough, we find that most  of the analysis does not depend upon the distribution of \  $\Friends$
but only on its expected value.}

 Let $\textbf{s} :=  (s_1, \cdots,s_N )$ and  $\bpgf (\textbf{s}, \eta ):= \sum_{i = 1}^N f(s_i,  \eta, \beta) \rho_i.$
The  post quality factor dependent  PGF,  of the offspring distribution of the overall branching process,  is given by (see equation (\ref{Eqn_offspring})):
\vspace*{-0.2cm}
 \begin{eqnarray}
\label{Eqn_hl}
\boxed{ h_l(\textbf{s})  =  \theta\left(s_{l+1}\mathbb{1}_{l<N} + \mathbb{1}_{l=N} \right) + (1-\theta)r_l \bpgf (\textbf{s}, \eta ) +(1-\theta)(1-r_l)}.
\end{eqnarray}

 \subsection{ Generator matrix}
The key ingredient required for analysis of any  MTBP  is its generator matrix. We begin with the generator for  MTBP that represents the evolution of unread TLs with CP-post.
We refer to this process briefly as \underline{TL-CTBP},   timeline continuous time branching process. The generator matrix, $A $, is given by $A = (a_{lk})_{N\times N}$, where $ a_{lk}= a_l\left( {\partial h_l({\bf s})}/{\partial s_k} \Big | {_{ {\bf s} = {\bf 1} }} - \mathbb{1}_{\{l=k\}} \right)$  and $a_l$ represents the transition rate of a  type-$l$ particle (see \cite{1AthreyaPaper} for details). 
For our case, from previous discussions   $a_l = \lambda + \nu$ for all $l$. 
Further, using equation (\ref{Eqn_hl}),   the  matrix $A$  for our single CP case is given by (with $c := (1-\theta) m\eta$, $ c_l = c\rho_l$)

\begin{eqnarray}
A = (\lambda + \nu)
\label{Eqn_genmatrix}
 \left[ \begin{array}{ccccccc}
 c_1 r_1 -1 & c_2r_1 + \theta & \cdots  & c_{N-1}r_1  & c_Nr_1 \\
 c_1 r_2   & c_2r_2 -1 & \cdots & c_{N-1}r_2  & c_Nr_2 \\
  & \vdots    \\
 c_1 r_{N-1}  & c_2r_{N-1}   & \cdots &  c_{N-1}r_{N-1} -1 &  c_{N} r_{N-1} + \theta \\
 c_1 r_N  &  c_2r_N  &  \cdots &  c_{N-1}r_N  & c_{N} r_N -1 \\
 \end{array} 
 \right] .
\end{eqnarray}
The largest eigenvalue and the corresponding eigenvectors of the above generator matrix are instrumental in obtaining the analysis of TL-CTBP (\cite{1AthreyaPaper})
and the following lemma establishes important properties about the same. 
We also prove that the resulting TL-CTBP is positive regular\footnote{A matrix   $B$  is called positive regular (irreducible) if there exists an $n$ such that the matrix $B^n$ has all strict positive entries. 
A BP is positive regular when its mean matrix is positive regular.
With $A$ as generator, the positive regularity is guaranteed if $e^A$ is positive regular (e.g. \cite{1AthreyaPaper}).
}, which is an important property that establishes the simultaneous survival/extinction of all the types of TLs.
 \begin{lemma}
 \label{Lemma_positive_regular}
  i) When $0< \theta < 1$, the matrix $e^{At}$ for any $t > 0$  is positive regular.  \\
ii) \ Let $\alpha $ be the maximal real eigenvalue of the generator matrix $A$. This eigen value lies in the real interval, 
 i.e.,  $\alpha \in \bigg(\textbf{r.c}-1, \ \textbf{r.c} -1+\theta \bigg) (\lambda+\nu)$, where  inner product $\textbf{r.c} := \sum_{i=1}^N  r_ic_i.$
When the reading probabilities have special form $r_l = d_1 d_2^l$ (for some  $0 \le d_1, d_2 \le 1$), then  

\begin{eqnarray*}
\alpha \to \left(\textbf{r.c} -1+\theta d_2\right)(\lambda + \nu) \ \text{as } N \to \infty.  
\end{eqnarray*}
iii)  The left and right eigenvectors ${\bf u} $, ${\bf v}$ corresponding to $\alpha$ satisfy the following equations  $ c_1\textbf{r.u}  = \sigma u_1$ and $  c_1\textbf{r.v}  =  \sigma v_N$ where $\sigma := \alpha/( \lambda+ \nu) + 1$. We have
 
\begin{eqnarray}
u_l &=& \sum_{i = 0}^{l-1} \frac{\rho_{l-i}}{\rho_1} \left(\frac{\theta}{\sigma}\right)^{i} u_1,   \ 2 \le l \le N \mbox{ \ and \ }  v_l = \sum_{i = 0}^ {N-l} \left(\frac{\theta}{\sigma}\right)^i \frac{r_{l+i}}{r_N} v_N  \ \ \ 1 \le l \le N-1.   \nonumber 
 \end{eqnarray}
 \textbf{Proof:} The proof is given in Appendix.  \eop 
\end{lemma}

At the finest details, we now developed a full-fledged MTBP that models the content propagation.  The multitype continuous time branching processes (MTBPs) are well studied in the literature (e.g., \cite{1AthreyaBook,1AthreyaPaper}).  The analysis of MTBP largely depends upon its generator matrix.  Lemma \ref{Lemma_positive_regular} describes the characteristics of the generator matrix specific to our model.  It yields in positive regularity of   TL-CTBP, i.e., the generator matrix $A$ (\ref{Eqn_genmatrix}) is positive regular. The largest eigenvalue $\alpha$ of   $A$ characterizes the growth rate of NU-TLs. We later see that left and right eigenvectors, \textbf{u} and  \textbf{v} (corresponding to $\alpha$) characterize the visibility of the CP-post.  Using the characterizations of Lemma \ref{Lemma_positive_regular} and the rich theory of MTBPs,  we derive various performance measures specific to this content propagation.

 The CP would be interested in many related performance measures as a function of post quality factor and we consider the same in the next section. 
\subsection{Performance analysis }
\label{sec_perf}
   If the CP invests sufficiently in preparing the content/post and ensures a good quality, the post can get viral.  It is important to note that the overall evolution of post depends on the number of seed TLs with CP-post. It is sufficient to consider one seed TL to derive various performance measures corresponding. This is because of properties of the branching processes: the analysis of process is quite similar when started with multiple seed TLs (i.e., growth rate is same).  The central questions in a branching process which are also relevant to our content propagation process include:  a) What is the extinction probability, i.e., the probability with which the entire population gets extinct?;  b) What is the rate at which the population grows?;  and c) What is the total progeny? etc.  We apply the well-known results addressing the above questions to our context and derive some performance measures.  We employ fixed point techniques to obtain the other performance measures. We begin with the probability of extinction.

\subsubsection{Extinction probabilities}
 Depending upon the context of the problem, for instance, an awareness campaign, the CP may be interested in knowing the chances of dissemination of its information to a large population, i.e., the chance of virality of its post. This probability can be obtained directly using the extinction probability of the corresponding MTBP, as explained below.  The CP-post is said to be extinct when it disappears completely off the OSN, i.e., none of the $N$-length TLs contain the CP-post eventually (as time progresses).
Let $q_l$ be the  probability with  which the process gets extinct when TL-CTBP starts with one TL of type $l$,
$$
 q_l :=  P \big (\textbf{X}(t) = \textbf{0} \  \text{for some} \ t>0 \big |\textbf{X}(0)=\textbf{e}_l \big ).
$$
Let ${\bf q} := \{q_1, q_2, \cdots, q_N\}$ represent the vector of extinction probabilities. 

Under  positive regularity conditions of Lemma \ref{Lemma_positive_regular}.(i) when a BP  is not extinct, the population  grows exponentially fast to infinity  (see \cite{1Harris,1AthreyaPaper} , etc).  This fact is established for our TL-CTBP in Theorem \ref{Theorem_virality}, provided in the later subsections.   {\it Thus, we have a dichotomy: the post gets viral with the exponential rate when it is not extinct and dies off completely otherwise.}  And hence the extinction probability equals one minus the probability of virality.
 \begin{lemma}
 \label{Lemma_extn}
   Assume  $0 < \theta < 1$ and  $ E[ \Friends log \Friends ] < \infty$  with $\Friends\log(\Friends) := 0$ when $\Friends = 0$. Then clearly $ E[ \zeta log \zeta ] < \infty $ for any post quality factor $\eta$.  
Hence we have the following:\\
 (i)  If $\alpha \le 0$,  extinction occurs   w.p.1, i.e., $\textbf{q} = {\bf 1} = (1, \cdots, 1)$;\\
(ii) If $\alpha > 0$, then\footnote{Vector ${\bf q} < {\bf s}$ if  $q_i < s_i$ for all components $i$.} $\textbf{q} < \textbf{1}$, i.e., the post gets viral with positive probability irrespective of the type of the   seed TL.
In this case the extinction probability vector ${\bf q}$ is the unique solution of  
 the equation, 
$h(\textbf{s}) = \textbf{s}$,  and lies
 in the interior of $[0,1]^N$. 
\end{lemma}
\noindent{\bf Proof} It follows from \cite[Theorems 1-2]{1AthreyaPaper} and by Lemma \ref{Lemma_positive_regular}. \hfill{$\blacksquare$}
\\
\\
It is easy to verify that the hypotheses of this lemma are easily satisfied by many distributions. For example, Poisson, Geometric etc. satisfy  $E[ \Friends log \Friends ] < \infty$.

By Lemma \ref{Lemma_extn}.(ii) the extinction probabilities are obtained by
  solving  $h({\bf s}) = {\bf s}$. The extinction probability can be obtained by conditioning on events and is given as  below, when the process starts a type-$l$ TL: 
 \begin{eqnarray}
\label{Eqn_extn_fixedpoint}
 q_l = \theta \Big (q_{l+1} \mathbb{1}_{\{l < N\}}  +  \mathbb{1}_{\{l = N\}} \Big) + (1-\theta)r_l \ \bpgf({\bf q},\eta) + (1-\theta)(1-r_l).
  \end{eqnarray}
 The above simplifies to:
  \begin{eqnarray}
q_{N-l}   =   (q_N -1) \sum_{i = 0}^l \theta^{l-i}\frac{r_{N-i}}{r_N}  + 1 \mbox{  for any } 1 \le l < N, \label{ExtProbSingleCP1}
\end{eqnarray}and the solution of the above provides the extinction probabilities.

\noindent{\bf Virality Threshold:}   
By Lemma \ref{Lemma_extn}.(ii)  the CP-post {\it gets viral, i.e., the  TL-CTBP survives and explodes   with non-zero probability}, when   $\alpha  > 0.$   When $N$ is sufficiently large, by Lemma \ref{Lemma_positive_regular}.(ii) and Lemma \ref{Lemma_extn}.(ii),
\begin{eqnarray}
\alpha \approx  (  m \eta (1-\theta) \bm{\rho.r} - 1    + \theta d_2 )   (\lambda + \nu) = (m \eta \bm{\rho.r}  - 1)  \nu  - (1- d_2) \lambda.
\label{Eqn_ViralTh}
\end{eqnarray}
It is well-known that the BPs survive with positive probability if the largest eigenvalue of the generator matrix, $A$, is positive (supercritical process). We have an (almost) equivalent of the same, i.e.; the TL-CTBP can survive when $ m\eta(1-\theta)\bm{\rho.r}   > 1 -\theta d_2$   (see (\ref{Eqn_ViralTh})),
for a BP   pitted against the shifting process. The virality threshold, denoted by $\bar{\eta}$, is defined    in terms of network parameters and is given by 
\begin{eqnarray}
\bar{\eta} & > & \frac{1 -\theta d_2}{m (1-\theta)\bm{\rho.r}}. \label{Eqn_EtaBar}
\end{eqnarray}
Thus, the virality chances are influenced by post quality  $\eta$, shift factor (1-$\theta$),  by the types of posts produced as given by $\bm{\rho}$, the mean number of friends $m$   and  the reading probabilities ${\bf r}$. In effect, the virality chances are influenced by factor, $(1-\theta )\eta \bm{\rho.r}.$ 
\\
\\
\textbf{ No-TL Case: What if all the effects of the TLs were neglected?}

Majority of the works (e.g., \cite{1BranchNonMarkov,1BranchVM}) considers study of content propagation without considering TL structure, and as mentioned before, this is an incomplete study. We would like to compare our conclusions with the case when the effects of TLs are neglected.  When users do not posses TL structure, there will be: 
\begin{enumerate}
\item  No notion of post residing at various levels, i.e., all posts reside at one level only, and so $N = 1$;   consequently its reading probability is one ($r_1 = 1$),  
and  further,  $\rho_1 =1, \rho_i = 0 \ \forall i >1$.
\item No notion of shifting effect, consequently $\theta = 0$ which is equivalent to saying \  $\lambda = 0.$
\end{enumerate}
The remaining modeling details of the content propagation are the same as before. In the view of  this, it is evident that the content propagates according to  a single type continuous time Markov branching process.  Thus, the analysis of this case boils down to a special case of the TL model (with $\lambda = 0, N = 1, \rho_1 = 1, r_l = 1 \forall \ l$ ). 
For this special case, from  equation  (\ref{Eqn_ViralTh}),  the rate of growth say $\alpha_{No-TL}$ is given by: 
  $$
\alpha_{No-TL}  \approx  (m \eta    - 1)  \nu. 
$$
 Observe that the post  gets viral when $m \eta > 1$, as is well understood in branching and viral  marketing literature (\cite{1BranchVM}). However, as mentioned before this neglects the key aspects of content propagation\textemdash effects of TLs. It is  accompanied by an erroneous conclusion  that the virality chances are  influenced  by $m$ and $\eta$ only. While in reality there is additional influence, which  is summarized by factor $(1-\theta )\eta \bm{\rho.r}$. 

In  \underline{No-TL} case, the extinction probability is obtained is given by solving  $ q = f(q,\eta) $ (substituting the parameter values in equation (\ref{Eqn_extn_fixedpoint})). It again becomes evident that the effect of post residing at various levels is disappeared.
The extinction probability is the same, whether it is started with one CP-post on level 1 or level 9.  This is again a wrong interpretation and the solutions of the equation  (\ref{ExtProbSingleCP1})/(\ref{Eqn_extn_fixedpoint})  provide the correct extinction probabilities which considers the influence of TLs. 
\\
\\
\textbf{ Influence of the Network Connectivity on Extinction Probability:}
We call an OSN sparsely connected when a sizable portion of the users have less number of friends (random), i.e., when they have a smaller mean number of friends, $m=E[\mathbb{F}]$. Whereas in a densely/highly connected OSN, a sizable number of users have a large number of friends and hence $m$ is large. 
We study the impact of network connectivity on the extinction probability from sparsely connected OSN to densely connected OSN.
When the mean   $m$  increases, the network becomes more active as the sharing of different posts becomes more pronounced. The TLs are flooded with different posts rapidly, so do the TLs containing post-{\bf P},  and one might anticipate an increase in its virality chances.  However, these TLs also receive the other posts rapidly, resulting in rapid shifts to their contents. 
Thus, with an increase in $m$,  the   $\lambda$ increases, and so does $\theta$.  
 \vspace*{-2.5cm}
 \begin{figure}[ht]
\begin{center}
\includegraphics[width = 10cm, height = 12cm]{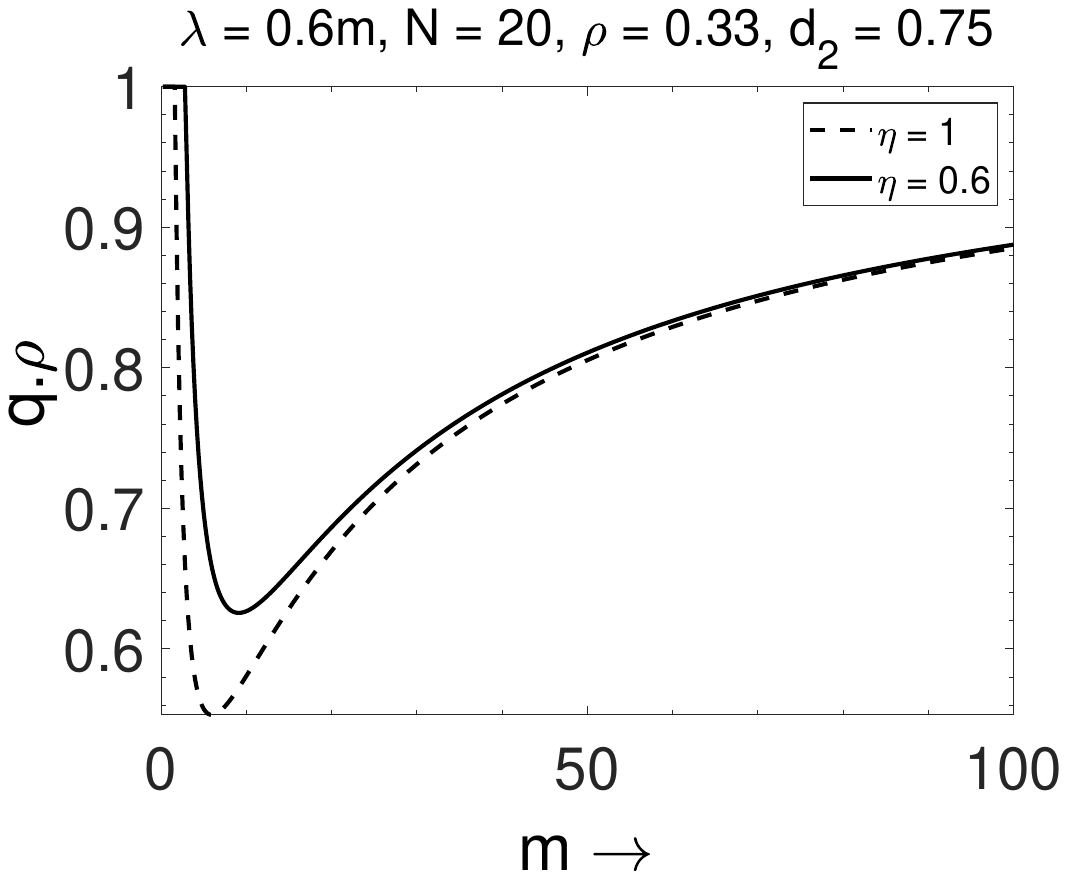} 
\vspace*{-3cm}
\caption{Extinction   vs $m = E[\mathbb{F}]$ \label{1Fig_Theta}} 
\end{center}
\end{figure}  
 \vspace*{-.5cm}
We observe  an interesting phenomenon in Figure \ref{1Fig_Theta},  with respect to the virality chances  $1-{q}_\rho := \sum_l q_l \rho_l$,  when $\lambda$ is set proportional to mean $m$.
To begin with, the virality chances $1-q_\rho$ improve ($q_\rho$ decreases) with mean $m$, as anticipated.  However, if one increases $m$ further,  we notice an increase in $q_\rho$. 
Basically, increased $m$ implies more shares of
 post-{\bf P} to new users \textit{but it also implies post-{\bf P} is missed more often}.      This phenomenon is mainly observed because of timeline structure: when TL structure is neglected, any user will view all the posts with equal interest irrespective of their levels. And consequently, one would not have noticed the effect of $m$  on extinction probability (as in No-TL case).   There seems to be {\it an optimal number of mean friends, which is best suited for post propagation.} 

\subsubsection{Time evolution of the NU-TLs}
The number of unread timelines (NU-TLs), at various time instances, may serve as an indicator of the reach of the CP-post. The reach of CP-post is another yardstick of the campaign effectiveness. In this section, we obtain the time evolution of NU-TLs. 
We have the following theorem which is instrumental in obtaining the expected number of NU-TLs in viral scenario:
   \begin{theorem} 
   \label{Theorem_virality}
 Let $(\Omega, \mathcal{F}, \mathbb{P})$ be an appropriate probability space and  let $\{\mathcal{F}_t \}$ be the natural filtration for 
TL-CTBP     $\textbf{X}(.)$, i.e.,  for each $t,$ 
     $\mathcal{F}_t $ is the $\sigma$-algebra generated by $\{\textbf{X}(t'); \ t' \leq t\}$.
The  process $ \{  \textbf{v}.\textbf{X}(t)e^{-\alpha t}; \ t\ge 0\}$, with ${\bf v}, \alpha$ as   in Lemma \ref{Lemma_positive_regular},    
       is a non negative martingale  (with natural filtration) 

\begin{eqnarray}
\hspace{13mm} \mbox{ and }
 \vspace{14mm}
\lim_{t\to\infty}\textbf{X}(t, \omega)e^{-\alpha t} = W(\omega)\textbf{u} \hspace{0.6cm} \text{ for   almost all   } \omega, 
 \label{LimitingSingleCP}
\end{eqnarray}  
 where $W$ is a non negative random variable that satisfies\footnote{ {\it We use $E_l$ and $P_l$ to represent the conditional expectation and  probability respectively when TL-CTBP    starts with one $l$-type TL.}}: $ P_l(W=0)=q_l$, $E_l[W] = v_l$ for each $l$, with \ ${\bf u.v}=1$.  
   \end{theorem}
   {\bf Proof:} Under the assumptions of Lemma \ref{Lemma_extn}, the TL-CTBP satisfies the hypotheses of  Theorem 1 of \cite{1AthreyaPaper}. \hfill{$\blacksquare$}
   
 The CP-post gets extinct on the sample paths with $W =0$ in equation (\ref{LimitingSingleCP}) (see \cite{1AthreyaPaper} for details, also observe 
 that $ P_l(W=0)=q_l$  from the above theorem). It  gets viral  in the complementary paths, i.e., when  $W > 0$,
 as is also evident from the limit\footnote{Note that for large $t$, the NU-TLs  $\textbf{X}(t, \omega) \approx W(\omega)\textbf{u}  e^{\alpha t} $, which grows exponentially fast when $W(\omega) >0$.} given by equation (\ref{LimitingSingleCP}). 

On the viral paths, we have two important measures: 1) the growth rate $\alpha$, and 2) the visibility of the post. The growth rate characterizes the rate at which the post spreads through the OSN.  From (\ref{LimitingSingleCP}),  the TLs grow exponentially fast with time at the rate $\alpha$ (given by equation (\ref{Eqn_ViralTh})), i.e., as according to $e^{\alpha t}$.    
And the other measure,  \textit{the visibility of the post}  can be determined by the number of \textit{potential} users that can read the post and thereby get influenced to buy the product/service. Recall that users attention is limited to the first few number of levels. Clearly,  the visibility of the post depends on the level at which it resides on the NU-TLs. The more the number of TLs having post on higher levels, the more the visibility. 
The number of potential users viewing the post on the level $l$ is approximately $r_l u_l e^{\alpha t}$ (large $t$) where $u_i$ is the $i$-th component of vector \textbf{u}. 
We define the visibility of the post at level say  $l$ as the  \textit{fraction} of NU-TLs holding the post at level $l$ after a long time $t$ which is given as   $u_l/\sum_i u_i$. 

 We also obtain the time evolution of the \textit{expected value of NU-TLs}. This result is obtained as a corollary of Theorem \ref{Theorem_virality}.
\begin{corollary} \label{cor_1}
When \ $\alpha > 0$\ and starting with one \ type-$i$ \ seed TL, \ \ $\sum_{l = 1}^N E_i\left[ X_l(t)\right] = e^{\alpha t}v_i \sum_{l} u_l.$ Further, when $r_i = d_1d_2^i, \ \rho_i = {\tilde \rho}\rho^i$ (with $\sum_i \rho_i=1$ and $0 < \rho \le 1$) for all $i$, we have  
\begin{eqnarray}
\sum_{l = 1}^N E\left[ X_l(t)\right] =  \varrho e^{\alpha t} d_2^{i-1} 
\mbox{ with } \varrho  :=  (1-d_2 \rho) \left(\frac{1}{1-\rho} - \frac{\theta}{\rho }  \frac{1}{\sigma -\theta} \right) \frac{(\sigma - \theta d_2)(\sigma \rho -\theta)}{(\sigma -\theta) (\rho -\theta)}.
\end{eqnarray}    
\end{corollary}
\textbf{Proof:}
Using the fact that $v.X(t)e^{-\alpha t}$ is a martingale and $\textbf{u.v} = 1$, one can write the following
\begin{eqnarray*}
 E\left[\textbf{v.X}(t)e^{-\alpha t}\right] &=& E[v.X(0)e^{-\alpha \times 0}] = v_i; \ \  \textbf{u}. E\left[\textbf{v.X}(t)e^{-\alpha t}\right] = \textbf{u}v_i  \hspace*{0.9cm} \Big(\mbox{recall} \ \ \textbf{X}_i(0) = 1 \Big)
 \\
 \textbf{u}. E\left[\textbf{v.X}(t)e^{-\alpha t}\right] & = &   \textbf{u.v}E\left[\textbf{X}(t)e^{-\alpha t}\right]= E\left[\textbf{X}(t)e^{-\alpha t}\right] = \textbf{u}v_i \hspace*{1cm} \because \textbf{u.v} = 1. 
\end{eqnarray*}
Thus $E\left[\textbf{X}(t)\right] = \textbf{u}v_i e^{\alpha t}$. Further, by taking the sum of individual component of the expected value of the random vector 
\begin{eqnarray*}
\sum_{l = 1}^N E\left[ X_l(t)\right] =  E\left[\sum_{l = 1}^N X_l(t)\right] = e^{\alpha t}v_i \sum_{l} u_l \hspace{1cm } \because N \  \mbox{is finite.}
\end{eqnarray*}
Substituting the value of $ v_i \sum_{l} u_l $ from equation (\ref{Eqn_vi_sum_u}) in Appendix,  we get the desired result
\begin{eqnarray*}
\sum_{l = 1}^N E\left[ X_l(t)\right] =  e^{\alpha t} d_2^{i-1} (1-d_2 \rho) \left(\frac{1}{1-\rho} - \frac{\theta}{\rho }  \frac{1}{\sigma -\theta} \right) \frac{(\sigma - \theta d_2)(\sigma \rho -\theta)}{(\sigma -\theta) (\rho -\theta)} = \varrho e^{\alpha t} d_2^{i-1} .
\end{eqnarray*}  
  \hfill{$\blacksquare$}
\subsubsection{Time evolution of the number of shares}
We derive another important performance measure,   the expected number of shares of the post,  before a given time  $t$. This measure gives the total spread of the post, i.e., the total number of shares a post gets in the given time-frame  (e.g., the number of shares in Facebook).  
It is basically the total number of distinct TLs (i.e., users)  that received a copy of the Post before time $t$.  
\textit{It is important to observe here that `number of shares'   is different from the  well known `total progeny'\footnote{The total number of offspring produced so far, by the BP.} of the underlying BP.} \ The `number of shares' is   due to offspring generated by share transition only, while the `total progeny'  is due to both `share' as well as `shift' transition offspring. 
We discuss the number of shares in viral ($\textbf{q} < \textbf{1}$) and non viral (sure extinction) scenario.
\\  
 \textbf{Viral scenarios:}    We employ probability generation based technique to obtain the time evolution of number of shares.  Let $Y(t)$  be  the accumulated number of  shares  till time $t$ and  let $Y =\lim_{t \to \infty} Y(t)$ (can also be infinity) be the eventual number of shares. The following Lemma captures the time evolution of number of shares.
  \begin{lemma}
  \label{Lemma_Shares} Let ${\bf y} (t) := [y_1 (t) \cdots, y_N(t)]$ with $y_l (t) := E_l [Y(t)] = E[ Y(t) | {\bf X}(0) = {\bf e}_l]$, the expected number of shares till time $t$ when started with one $l$-type TL for each $l$. If $\alpha > 0$, we have
  \begin{flalign}
  \label{Eqn_yt1}
  {\bf y} (t)&=
e^{At} \Big( \textbf{1} + (\lambda + \nu) A^{-1} \textbf{k} \Big) - (\lambda + \nu)A^{-1}\textbf{k}  
\mbox{ where } \textbf{k}  =  \  [1-\theta, 1-\theta, \cdots, 1-\theta,1 ]^T.
\end{flalign}
\textbf{Proof:} The proof is given in Appendix. \eop
\end{lemma} 

Thus,  the expected number of shares grow exponentially fast,  with time, for viral scenarios.  Further, the growth rate $\alpha$ (see eqn. (\ref{Eqn_ViralTh})) is the same as that for the unread posts.
From  (\ref{Eqn_yt2}), 
  for large $t$,  
the expected shares when started with one type-$l$ particle is:
  \begin{eqnarray}
  \label{Eqn_yt}
  y_l (t)   \approx  e_{l,0} e^{\alpha t}  \mbox{ with }  e_{l,0} =     v_l \sum_{i=1}^N u_i \left (1 + \frac{\nu}{\alpha}\right ) + v_l \frac{\lambda}{\alpha} u_N.
 \end{eqnarray}  
 \textbf{Non viral scenarios:} 
When population gets extinct with probability one, the expected  number of  total shares is finite.
One can directly obtain the expected number of shares by conditioning on the first transition event as follows:

\begin{eqnarray}
\label{Eqn_yl}
y_l   := E_l[Y] =  \theta y_{l+1} \mathbb{1}_{\{l < N\}} + (1-\theta) r_l \left(m\eta + m \eta \textbf{y.}\bm{\rho}\right); \ \ \  \mbox{for all } l\le N. 
\end{eqnarray}
\noindent  
On recursively simplifying the above system of equations backward, we obtain the following for any $l \le N$
\begin{eqnarray}
y_l   =  (1-\theta)m\eta(1+\textbf{y.}\bm{\rho}) \sum_{i=0}^{N-l} \theta^{N-l-i} r_{N-i}. \label{FiniteSharesSingleCP} 
\end{eqnarray} 
Summing the above over $l$ after multiplying with $\rho_l$, we obtain:
\begin{eqnarray*}
\textbf{y.}\bm{\rho} = \sum_{l=1}^N \rho_l y_l  =  (1-\theta)m\eta(1+\textbf{y.}\bm{\rho})  \sum_l \rho_l \sum_{i=0}^{N-l} \theta^{N-l-i} r_{N-i} .
\end{eqnarray*}
Thus, the FP equation for $\textbf{y.}\bm{\rho}$ is  linear and
  hence we have a unique FP solution for $\textbf{y.}\bm{\rho} $
   whenever
$$(1-\theta)m\eta  \sum_l \rho_l \sum_{i=0}^{N-l} \theta^{N-l-i} r_{N-i} <  1.$$ 
 If
$(1-\theta)m\eta{\bf r.\bm{\rho}}-1+\theta= {\bf r.c}-1+\theta < 0,$ from  Lemma \ref{Lemma_positive_regular}.(ii) $\alpha < 0$ and the process would be extinct w.p. one. In this scenario:
\begin{eqnarray*}
(1-\theta)m\eta  \sum_l \rho_l \sum_{i=0}^{N-l} \theta^{N-l-i} r_{N-i}   
&\le& (1-\theta)m\eta  \sum_l \rho_l  r_l \sum_{i=0}^{N-l} \theta^{N-l-i}
\\
&  & \hspace*{-2cm} =
(1-\theta)m\eta  \sum_l \rho_l  r_l \frac{1- \theta^{N-l-1}}{1-\theta} 
= m \eta {\bf r.\bm{\rho}} < 1
\end{eqnarray*} 
because $r_1 \ge r_2 \cdots \ge r_N.$ We can similarly show using the limit of the eigenvalue $\alpha $ of Lemma \ref{Lemma_positive_regular}, that when the process is extinct w.p. one, the above condition is  always satisfied  asymptotically.  To be more precise the condition is   satisfied for all  $N$ bigger than a threshold ${\bar N}$, whenever the process is extinct w.p. one.  

We thus have the following  unique FP for $\textbf{y.}\bm{\rho}$ under the conditions discussed above:
\begin{eqnarray}
\label{Eqn_single_CP_yrho}
\textbf{y.}\bm{\rho} = \frac{(1-\theta)m\eta  \sum_l \rho_l \sum_{i=0}^{N-l} \theta^{N-l-i} r_{N-i} }{ 1-  (1-\theta)m\eta \sum_l \rho_l \sum_{i=0}^{N-l} \theta^{N-l-i} r_{N-i}}.
\end{eqnarray}

One can substitute the above in equation (\ref{FiniteSharesSingleCP}) to obtain $y_l$ for all $l$:
\begin{eqnarray}
y_l   =  \frac{ (1-\theta)m\eta   \sum_{i=0}^{N-l} \theta^{N-l-i} r_{N-i} }{ 1-  (1-\theta)m\eta \sum_l \rho_l \sum_{i=0}^{N-l} \theta^{N-l-i} r_{N-i}}  . \label{FiniteSharesSingleCP1} 
\end{eqnarray} 
Also, it is easy to verify that the FP is unique by uniqueness of the FP solutions for $\textbf{y.}\bm{\rho}$.
 \\
Note that in No-TL case, the number of shares is computed using
 \begin{eqnarray}
 y = E[Y] = m\eta(1+y) \implies y = \frac{m\eta}{1-m\eta} \label{Eqn_NoTLSharesExt}
 \end{eqnarray}
which is again inaccurate.

{\bf Special case:} Say  $r_i = d_1d_2^i$,  $\rho_i = {\tilde \rho}\rho^i$ (with $\sum_i \rho_i=1$ and $0 < \rho \le 1$) for all $i$,
one can easily  simplify the above. 
We have the following
\begin{eqnarray*}
\sum_l \rho_l \sum_{i=0}^{N-l} \theta^{N-l-i} r_{N-i}
  &= &  d_1 {\tilde \rho}  \sum_l \rho^l \sum_{i=0}^{N-l} \theta^{N-l-i} d_2^{N-i}   
   =     d_1 {\tilde \rho}  \sum_l \rho^l d_2^l \sum_{i=0}^{N-l} \theta^{N-l-i} d_2^{N-l-i} 
   \\  
     & = &   d_1 {\tilde \rho}  \sum_l \rho^l d_2^l \sum_{i=0}^{N-l} \theta^{ i} d_2^{i}   =   d_1 {\tilde \rho}  \sum_l \rho^l d_2^l  \frac{ (1- (\theta d_2)^{N-l+1} ) }{1-\theta d_2} 
      \\
      & =& \frac{ d_1 {\tilde \rho} }{1-\theta d_2} 
      \left ( \rho d_2  \sum_{l=0}^{N-1} \rho^l d_2^l - d_2   (\theta d_2)^{N}  \rho  \sum_{l=0}^{N-1} \rho^l  \theta^{-l}   \right) 
      \\
       &= &   
      \frac{ d_1 {\tilde \rho} }{1-\theta d_2} 
      \left ( \rho d_2 \frac{1- (\rho d_2)^{N}}{1-\rho d_2} -    (d_2)^{N+1}  \rho   \frac{ \theta^N - \rho^N}{\theta(\theta - \rho)}    \right).
\end{eqnarray*} 
Substituting this in equation (\ref{Eqn_single_CP_yrho}) and  under the limit $N \to \infty$, we obtain the following compact expression (where 
$\tilde{\rho}  = (1-\rho)/\rho$  now in the limit):
\begin{eqnarray}
\label{Eqn_yrho}
\textbf{y.}\bm{\rho} \approx  \frac{O_{mean}}
{1- O_{mean}} \mbox{ with }  O_{mean} :=(1-\theta)m\eta \frac{ (1- \rho) d_1d_2 }{(1-\theta d_2)(1-\rho d_2)}.
\end{eqnarray}

\subsection{Validation of the number of shares}
 We validate our theoretical expression for the expected number of shares by Monte Carlo simulation based on a real dataset; \textit{Stanford Large Network Dataset Collection} (SNAP) dataset as provided in ego-Facebook, Social Networks section \cite{1SNAP}. The dataset consists of friends' list of 4039 Facebook users and undirected connections among them. The sum of the number of friends of all these 4039 users (undirected connections) stands at 88234.   To judiciously validate the theoretical finding, we add new users to the existing dataset as it has insufficient users originally. Basically, we split the friends of the nodes that have higher degree of connections into multiple sets. We then created new users and made  undirected connections  by randomly choosing the nodes from each of the above-mentioned sets. (We now have a total of 20109 users.) 

\begin{figure}[H]
\begin{minipage}{10cm}
\begin{center}
\includegraphics[width = 9.5cm, height = 10cm]{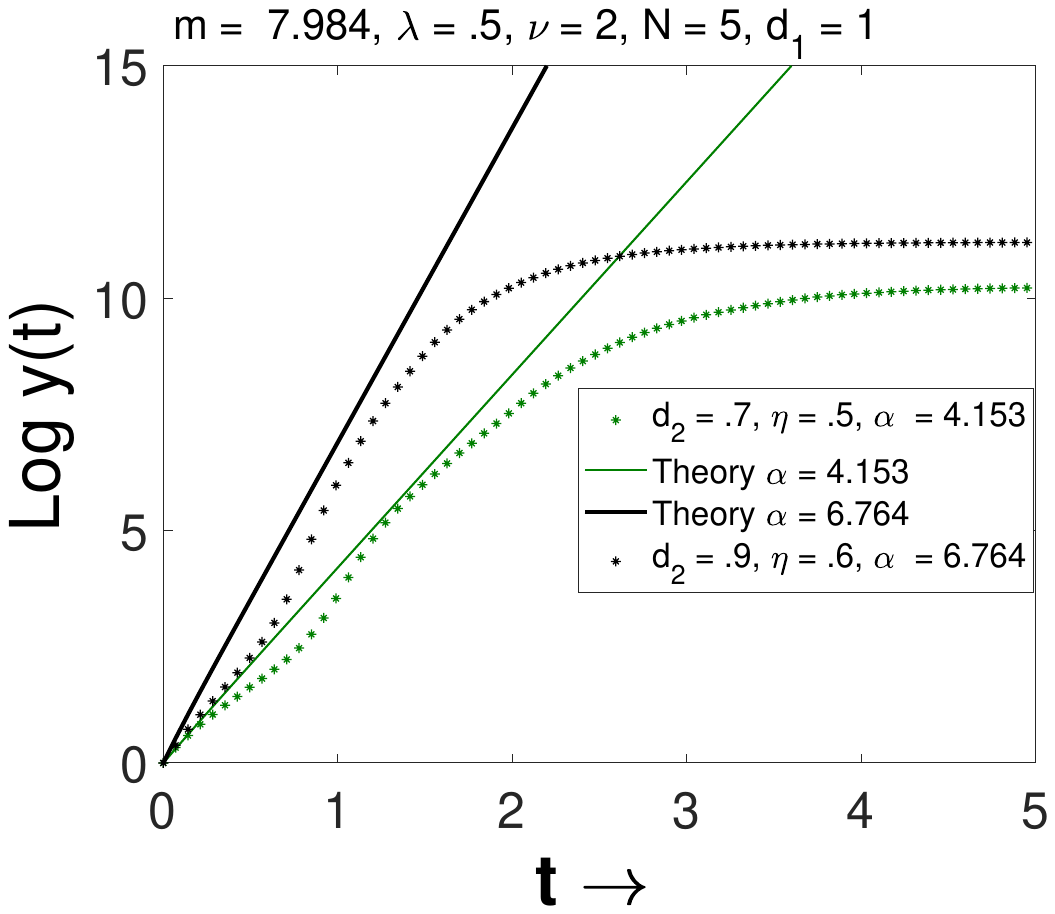} 
\vspace*{-3cm}
\caption{Time evolution of the expected shares: Theory vs Simulation \label{1SharesTheoryVsSimFB}} 
\end{center}
\end{minipage}
\end{figure}

 We emulate our content propagation model on the above dataset as follows. We represent each user by a TL comprising five levels (N = 5). The starting type-1 seed TL reads the CP-post with probability $r_1$ shares it with a random number of friends from its friends' list (as in the dataset) while influenced by the post quality factor $\eta$. We incorporated all the other details, e.g., shifting, the lifetime of a TL, etc into the simulation.  We obtain the number of shares in each sample path (realization) at fixed points in time. We then computed the average number of shares generated in such 8000 sample paths at each of the fixed time instances $y_1(t)$, i.e., the time evolution of the expected number of shares. 
 We plot the time evolution of the expected number of shares obtained theoretically and via simulation in Figure \ref{1SharesTheoryVsSimFB}).  And we compute the number of shares  theoretically using the same set of values. For the sake of convenience, we use the natural log scale on the $y-$axis.
 \vspace*{-1cm}
 
  As the number of users are finite (dataset), the trajectory of $\log y_1(t)$  begins to saturate as time elapses in Monte Carlo simulation. While theoretically, the expected number of shares continues to grow indefinitely.  We  see in  Figure \ref{1SharesTheoryVsSimFB}) that the theoretical trajectory  of $\log y_1(t)$ matches well to that of the simulation based trajectory till saturation.

\section{Viral marketing and real time bidding} \label{Ch_oxn}
The performance measures obtained in the previous sections can be useful in many advertisement/campaign related objectives such as brand awareness, search engine optimization, maximizing the number of clicks to a post/advertisement (ad), etc. 
In this section, we will study online auctioning for advertisements in viral marketing using the performance measures as obtained in the previous sections.

The publishers of OSNs sell the advertisement inventory/space to various content providers (CPs) via auction mechanism commonly known as \textit{real time bidding} (\cite{1BidEst}). For example, Facebook auctions billions of advertisement space inventory every day, and the advertisements  (ads) of the winners are served. Real-time bidding enables the CPs to automatically submit their bids in real time, and the advertisement of the highest worth (based on bid amount and its performance) is thus served. By virtue of auctioning, a natural competition occurs among the CPs for winning auctions. A content provider (CP) has to win the auction to get sufficient number of seed (initial) timelines. 
The virality/sharing of the post further depends upon the quality of the advertisement/post (recall the post quality factor $\eta$). On summarizing, the CP has to invest in two aspects: a) the bid amount to win the auction, and b) the amount spent to the design of the post ($\eta$).
Recall that designing of a post could include providing authentic information about your services/products, or providing quality content,  or giving offers, etc.  Inappropriately tailored post can make users lose interest in the post, and thereby reducing the virality chances.

Content providers (CPs) typically have wide-ranging objectives while advertising on OSNs. For example, a CP may be in interested in enhancing the brand awareness of its products. Brand awareness plays a central role in users' decision making for a purchase. Such an objective is achieved if the brand promotional post gets viral. 
Recall, we say a post gets viral if it spreads on a massive scale via its sharing among the users.  Given that a post gets viral, a CP may be interested in knowing how fast the post spreads, i.e., the rate of virality. Other objectives, a CP may be interested in, include:  maximizing the number of clicks on its post, improving its reputation, increasing its presence in the marketplace, etc. 

In previous sections, we derived some of these performance measures. For example, we obtained the time evolution of the number of shares and NU-TLs which characterize the rate of virality. We also obtained the expression for the probability of virality. On the other hand, in non-viral (sure extinction) scenarios, we computed the expected number of total shares before extinction. We provided  explicit expressions for some of the performance measures as a function of controllable parameters while others are represented as the solutions of appropriate  FP (fixed point) equations. One can use these measures to study a relevant optimization problem taking auctions into account. In particular, and without loss of generality, we take the expected shares/NU-TLs as an indicative of the performance of CP's posts.

\subsection{ Optimal budget allocation}  
In the single CP model\footnote{When we consider the study of competing content in Part-II, the  CPs further compete over relative visibility of their own content (details in Part-II \cite{part2}). }, the CP (indirectly) competes with other  CPs only for advertisement inventory space (i.e., for winning initial seeds). This is because the other CPs are advertising unrelated content. 
We consider the details related to winning auctions, and then the resultant rewards derived by single CP. The CP has to first win the auction, and then its post will propagate via the forwards/shares as discussed before. Recall that these shares/forwards generate revenue to the CP.
 Therefore, it becomes important for the (concerned) CP to know the bid distribution of the other CPs/advertisers. In particular, we need the highest bid of the advertisers participating in the auction. Authors in \cite{1BidEst} show that the maximum bid value follows the log-normal distribution with  parameters mean $\mu_{b}$ and variance \ $\sigma_{b}^2$.  \ Let $\mathbf{B}$ denote the distribution of the maximum bid values, it follows form \cite{1BidEst}
 \begin{eqnarray}
 \log \mathbf{B} \sim N(\mu_b, \sigma^2_b)  \hspace{0.3cm} \mbox{where \ }  N(\centerdot,\centerdot) \ \mbox{ \ is Gaussian distribution.} \label{Eqn_BidDist}
\end{eqnarray}  

 As mentioned before, we take the expected value of NU-TLs as one among several choices of performance measures to study the optimization problem. More specifically,  we take the sum of the expected number of users with CP-post at various levels, $\sum_{l=1}^N E[X_l(t)]$ for some large $t$, as an indicator of CP's revenue. Note that when the characteristics of underlying social network (e.g., sparsely connected) are such that the probability of extinction is \textit{one}, non-viral scenario, the CP gets  \textit{zero} reward as the NU-TLs become zero after some time.  Whereas in the viral scenarios (i.e.,  $\textbf{q} < \textbf{1}$), we have $ \bar{\eta} \le \eta \le 1$ (see (\ref{Eqn_EtaBar})) and the CP gets  $\sum_{l=1}^N E[X_l(t)]$ provided that it wins the auction. 

Authors in \cite{1xeta} state that the winner of the auction is decided based on the bid amount and the corresponding quality of the post/advertisement collectively. In other words, the CP wins the auction when the bid amount $x$ and $\eta$ collectively exceeds the bid distribution, i.e., $ x \eta > \mathbf{B} $. Thus, the probability of winning the bid is $P( \mathbf{B} < x \eta) $, which is the cumulative density function (CDF) of log-normal distribution. Given that the CP wins the auction, its content is placed at the top-level of one TL, i.e., we begin content propagation with one seed TL of type-1. By the time the seed user visits its TL, the post might have shifted down or might disappear completely off the TL. 
 
 In all, the CP invests: 1) $x$ for the bid amount to win the auction, and 2)  $\kappa_1 \eta $ for preparing the post and $\kappa_1 >0$.  Let us say the CP wants to maximize its utility, denoted by $\mathbf{C}\left(x,\eta \right) $ where
  \[   
\mathbf{C}\left(x,\eta \right)  =  
     \begin{cases}
       \left(\log E\Big(\sum_{l} X_l(t)\Big)  -  \kappa_2(x +\kappa_1 \eta)\right)P( \mathbf{B} < x \eta), & \text{if } \ \bar{\eta} \le \eta \le 1
       \\
       0, & \text{else},
     \end{cases}
\]
\\
where the weightage $\kappa_2$ captures  trade-off between the reward 
$\log E\Big(\sum_{l} X_l(t)\Big)$  and the overall cost $x +\kappa_1 \eta$. 

The close-form expression of CDF of log-normal distribution with erf as the error function
\begin{eqnarray*}
P( \mathbf{B} < x \eta) = \frac{1}{2} + \frac{1}{2} \mbox{erf}\left(\frac{\log x \eta - \mu_b}{\sqrt{2} \sigma_b}\right) = \frac{1}{2} + 
 sign \left(f(x\eta)\right)  \frac{1}{\sqrt{\pi}}  \int_{0}^{ \abs{f(x \eta)}} e^{-z^2} dz,
 \end{eqnarray*}
 \\
  where $sign(a) = +1 \ $ if $a \ge 0$ and -1 otherwise; and $ f(x\eta)  = \frac{\log x \eta - \mu_b}{\sqrt{2} \sigma_b}$. 
 Using Corollary \ref{cor_1}, we rewrite it as 
 \[   
 \mathbf{C}\left(x,\eta \right)  =  
     \begin{cases}
       \left(\log \Big( e^{\alpha t}v_i \sum_{l} u_l \Big) - \kappa_2( x +\kappa_1 \eta)\right) \left(\frac{1}{2} + 
 sign\left(f(x\eta)\right)  \frac{1}{\sqrt{\pi}} \int_{0}^{\abs{f(x \eta)}} e^{-z^2} dz \right), \\ &  \hspace{-10mm}\text{if } \   \bar{\eta} \le \eta \le 1
       \\
       0, & \text{else.}
     \end{cases}
\]
Thus, the optimization problem is stated as: 
\begin{eqnarray}
O1: \ \ \max_{x,\eta}\mathbf{C}\left(x,\eta \right) \ \ s.t.  \ x \ge 0 \ \mbox{and} \ 0 \le \eta \le 1.
\end{eqnarray}
In some scenarios, the CP is constrained by   limited budget. Given a budget amount say $\bar{B}$, how to  allocate/divide the same  into bid amount $x$ and  $\eta$  related cost,  such that the revenue is maximized;  in other words we consider constrained optimization (revenue maximization) problem  under  the budget constraint $\mathcal{B}(x,\eta) = x + \kappa_1 \eta \le \bar{B}$.  This leads to the formulation of a\textit{ variant of the above stated optimization problem}:

\begin{eqnarray}
O2: \ \ \max_{x,\eta} \ \log E\Big(\sum_{l} X_l(t)\Big)P(  \mathbf{B}\leq x \eta) \ \ s.t. \ \ \bar{\eta} \le \eta \le 1, \ x \ge 0, \ \ x + \kappa_1 \eta \le \bar{B}.
\end{eqnarray}
Optimizers of the above problem give the best allocation of the  available  budget to the following factors: 1) winning the auction,  and 2) maintaining post quality such that overall spending does not exceed $\bar{B}$.
\begin{proposition} \label{Prop_BgtTightness}
Any pair of optimizers, $(x^*, \ \eta^*)$, of \ $O2$ \ satisfy \  $x^* + \kappa_1 \eta^* = \bar{B}.$
\end{proposition}
{\bf Proof} The proof is given in Appendix.  \eop

Due to the complex nature of the underlying objective functions, it is hard to  analyse both of the optimization problems analytically.  In particular, we are interested in obtaining the optimizers and study their variations with different system parameters in both the optimization problems. Let us say $C^*$ and $C^*_{con}$ be optimal objective values of $O1$ and $O2.$ We compare and contrast the optimizers and objective values of $O1$ and $O2$ in the plots below.

\begin{figure}[H]
\begin{minipage}{9cm}
\begin{center}
\includegraphics[width = 6.5cm, height = 9cm]{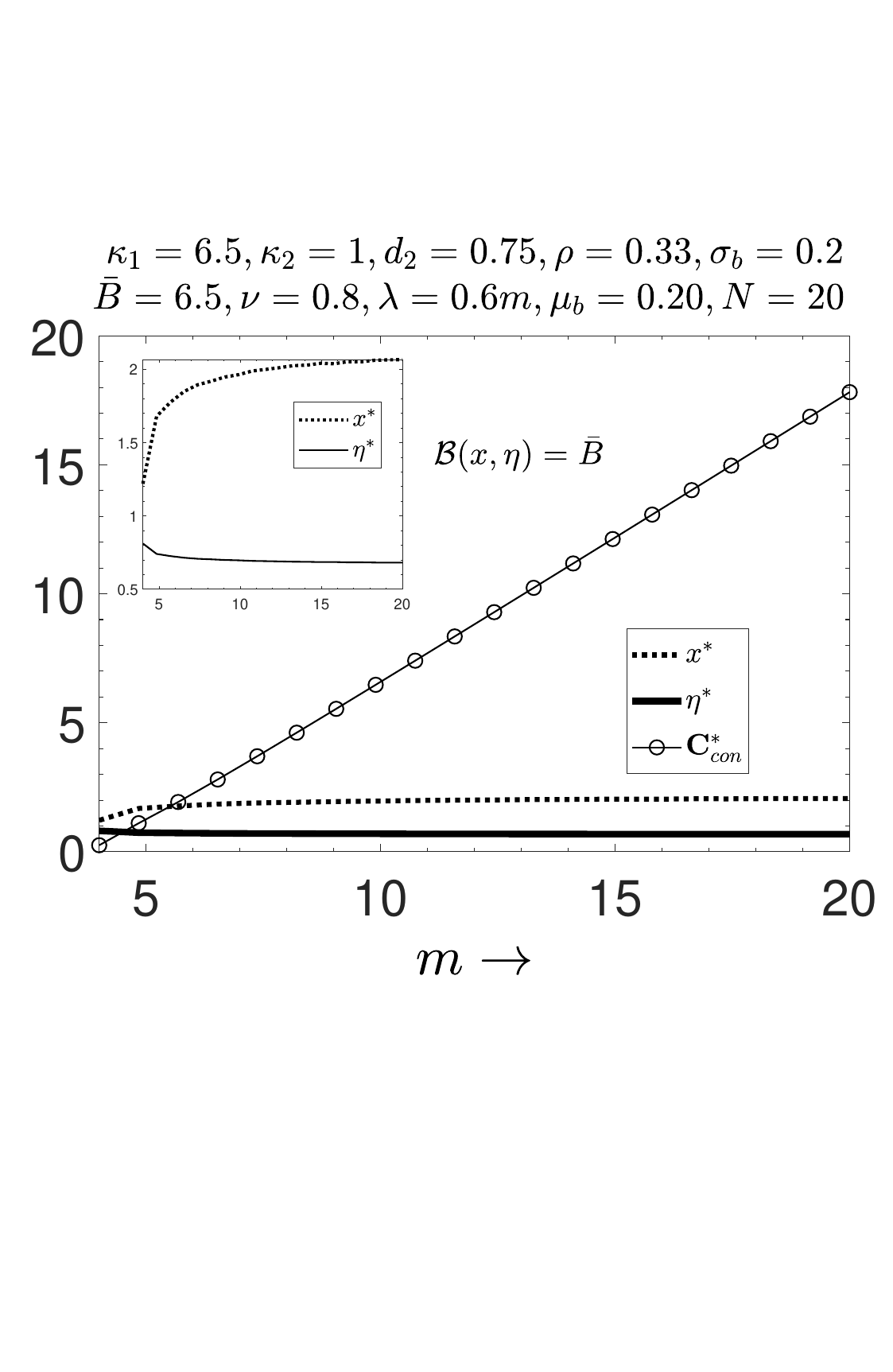} 
\vspace*{-2.6cm}
\caption{$O2$ with TL \&  $\lambda \propto m$  \label{1Fig_ConMeanLam_6}} 
\end{center}
\end{minipage}
\hspace*{-2cm}
\begin{minipage}{9cm}
\begin{center}
\includegraphics[width = 7.5cm, height = 9cm]{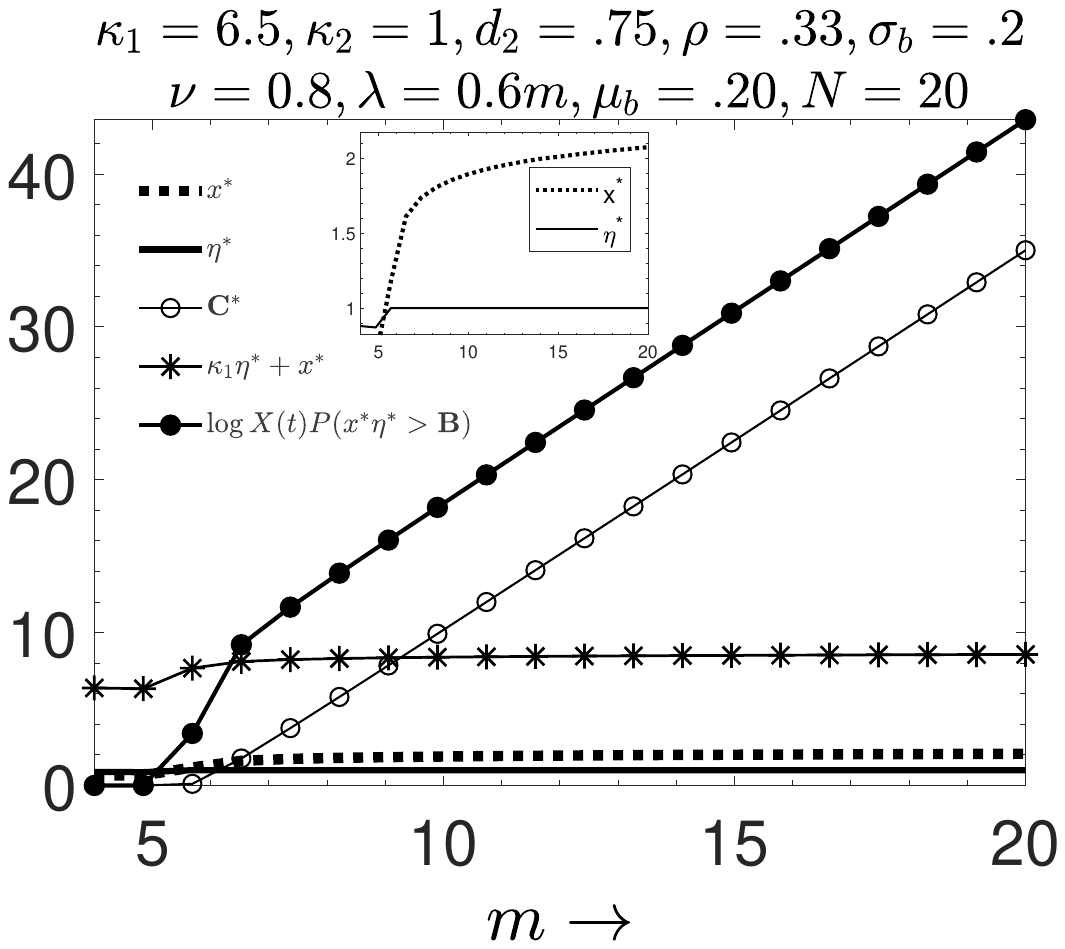}
\vspace*{-2.6cm}
\caption{$O1$ with TL \& $\lambda \propto m$  \label{1Fig_UnMeanLam_6}} 
\end{center} 
\end{minipage}
\end{figure}

Figures \ref{1Fig_ConMeanLam_6}  and \ref{1Fig_UnMeanLam_6}  depict the CP's  spending 1) on the bid amount  for  winning the auction,  $x$; and 2) on the post quality factor $\eta$ in order to maximize its utility.  When the CP has a limited budget, i.e., as in optimization problem $O2$,  we see in Figure \ref{1Fig_ConMeanLam_6} that  $x^*$ increases whereas $\eta^*$ decreases as the mean number of friends increases ($m$). Eventually, both settle at their respective constant values, i.e., $x^* \approx 2.08, \ \eta^* \approx 0.69$. This pattern is attributed to two factors: 1) the cost factor for $\eta$, i.e., $\kappa_1$ is comparable to $\bar{B}$; and 2) the increasing mean number of friends accounts for the steady decrease of  $\eta$ to 0.69. In other words, as $m$ increases, the post can get viral with smaller $\eta$ (see equation (\ref{Eqn_EtaBar}), and hence, the CP tends to proportionally invest more in winning bid. This kind of trend is seen only when $x$ and $\eta$ are taken together as in budget constraint. While in optimization problem $O1$, due to the absence of budget constraint, we immediately see in Figure \ref{1Fig_UnMeanLam_6} that $x^*$ increases unrestricted and $\eta^*$ also increases to its maximum value one as the network activity increases (measure by $m$).  Hence, we see higher optimal objective values attained in $O1$ compared to those of $O2.$ Basically, the CP can utilize the increasing connectivity of network (i.e., $m $ increases) by investing more in $x$ and $\eta$. Note that the trend is different for the extinction probability as in Figure \ref{1Fig_Theta}. There we saw that as connectivity of network increases, the virality chances decreases. However, as seen now, the virality chances reduces but the expected shares still improves.

\begin{figure}[H]
\begin{minipage}{9cm}
\begin{center}
\includegraphics[width = 7.5cm, height = 10cm]{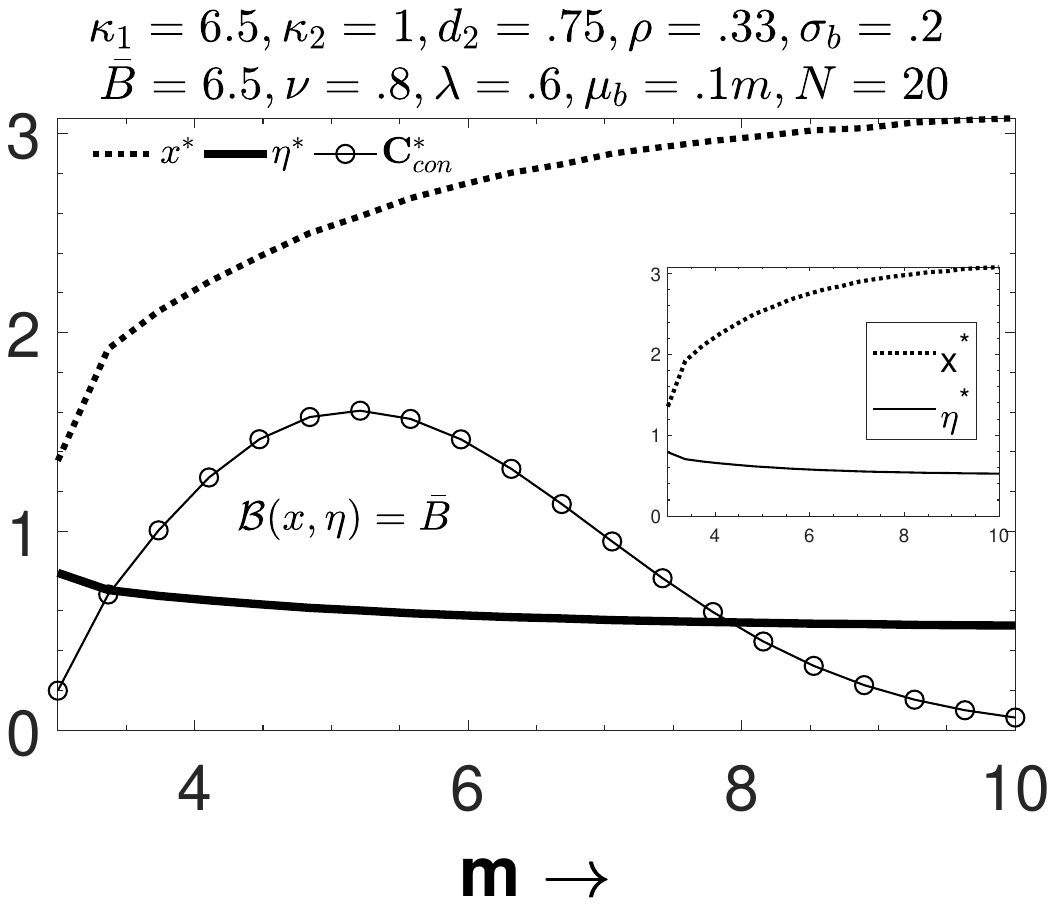}
\vspace*{-3cm}
\caption{$O2$ with TL \& $\mu_b \propto m$  \label{Fig_ConMeanMu0P1}} 
\end{center}
\end{minipage}
\hspace*{-1.8cm}
\begin{minipage}{9cm}
\begin{center}
\includegraphics[width = 8cm, height = 10cm]{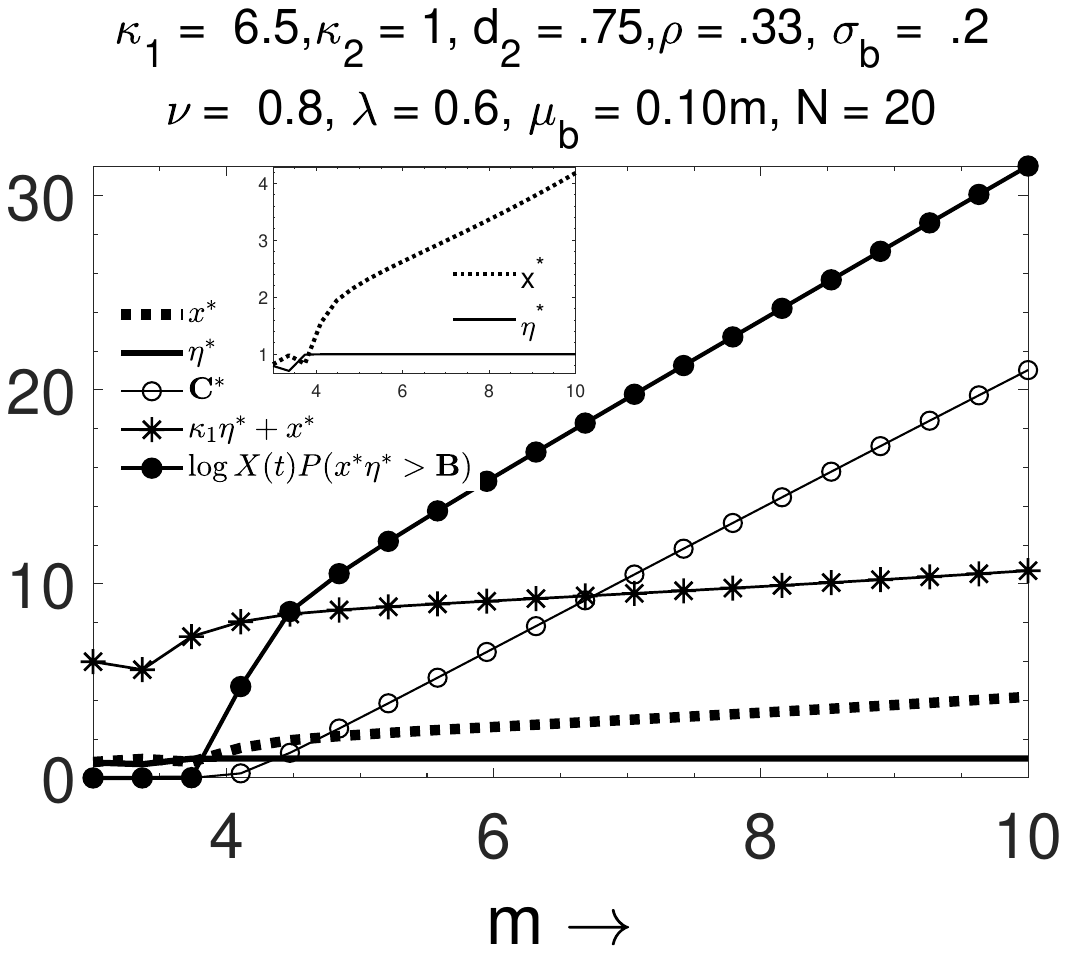} 
\vspace*{-3cm}
\caption{$O1$ with TL \& $\mu_b \propto m $  \label{Fig_UnMeanMu0P1}}. 
\end{center}
\end{minipage}
\end{figure}

While in $O2$, the CP steadily increases  allocation in $x$, and hence, proportionally invests lesser in $\eta$. And  both eventually settle to the constant values.
\vspace*{-1.5cm}
When mean number of friends increases, it is natural that the CPs would bid more as it would be easier for a content to get viral (recall $\alpha \propto m$). In other words,  an increase in $m$  implies an increase in $\mu_b$ ( $\mu_b \propto m$).
 When the mean of bid distribution increases, it gets difficult to win the auction  in constrained problem $O2$ (see Figure \ref{Fig_ConMeanMu0P1}). Consequently, the CP has to invest more in winning the auction, which comes at the cost of reducing the post quality $\eta$ ($x$ can not increase unrestricted due to the budget constraint). Further, as explained earlier the increasing mean number accounts for the steady decrease in $\eta^*$, and thereby, $x^*$ increases and both of them converge to the fixed values as can be seen in Figure \ref{Fig_ConMeanMu0P1}. Note that the objective value $C_{con}^*$, in this case, decreases after $m \approx 5$ because the allocation to $x$ is considerably higher than that seen in Figure \ref{1Fig_ConMeanLam_6}. Again in Figure \ref{Fig_UnMeanMu0P1}, without the  budget constraint (in $O1$), we do not see the trend. The optimal value increases with an increase in $m$, as $x^*$ can take unrestricted values. 
\vspace*{-1cm}

\begin{figure}[H]
\begin{minipage}{9cm}
\begin{center}
\includegraphics[width = 8cm, height = 10cm]{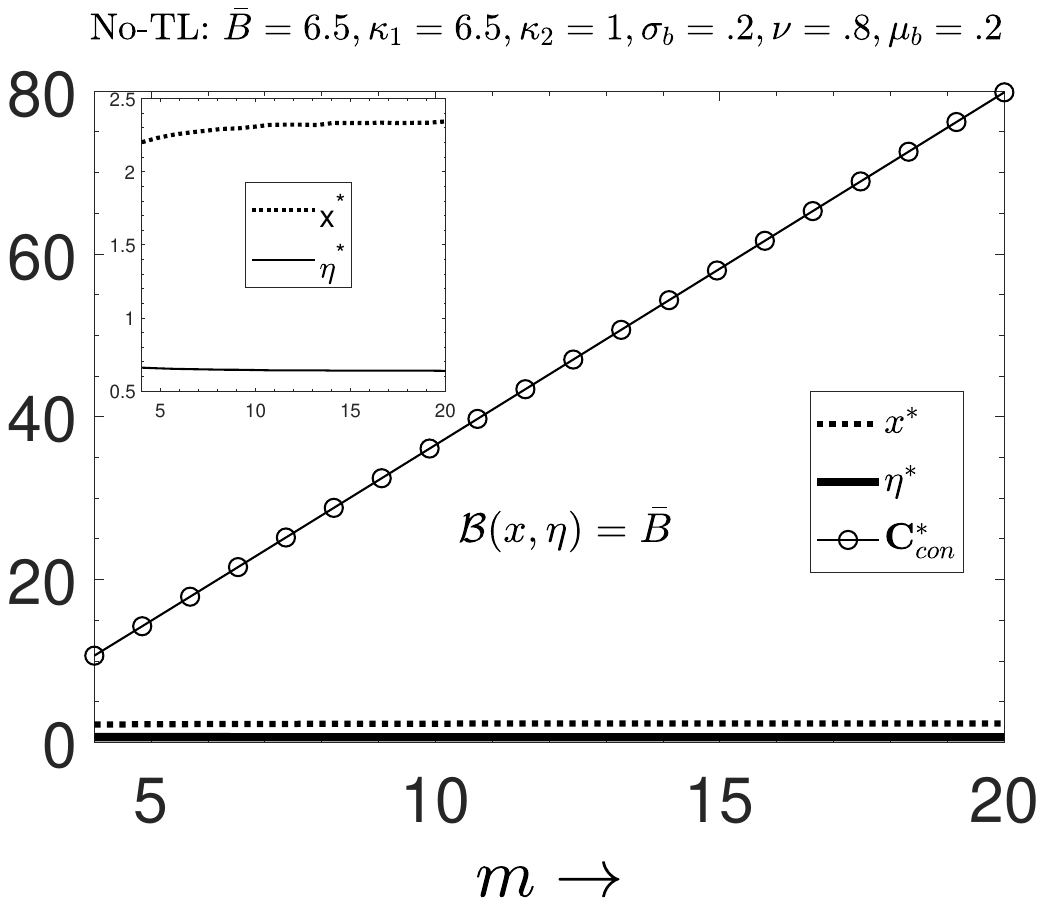}
\vspace*{-3cm}
\caption{Problem $O2$ without TL  \label{1Fig_ConNoTL}} 
\end{center}
\end{minipage}
\hspace*{-2cm}
\begin{minipage}{9cm}
\begin{center}
\includegraphics[width = 7cm, height = 10cm]{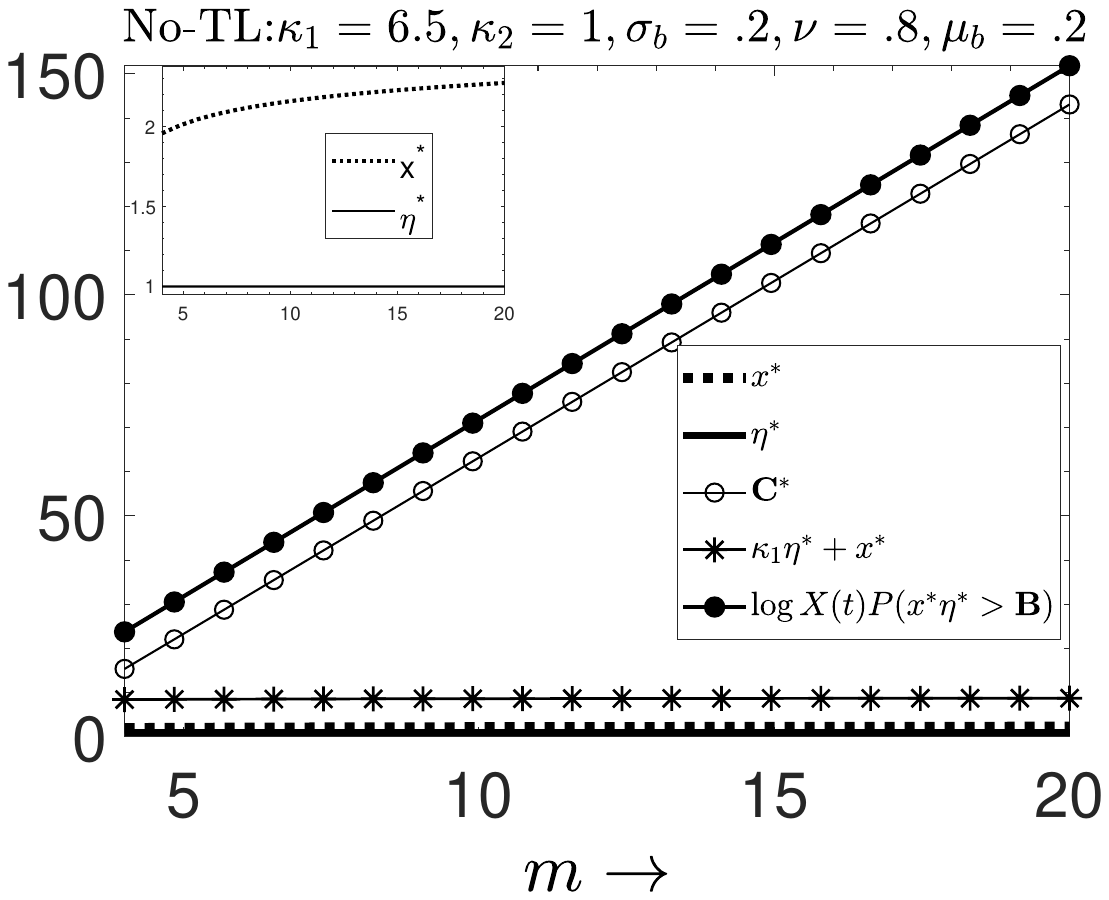} 
\vspace*{-3cm}
\caption{Problem $O1$ without TL  \label{1Fig_UnNoTL}} 
\end{center}
\end{minipage}
\end{figure}

\textbf{Impact of timeline structure on optimizers:}
Earlier we studied the impact of TL structure on the post propagation. We now see  through Figures  \ref{1Fig_ConNoTL}  and \ref{1Fig_UnNoTL}  how neglecting the TL structure influences the optimizers.  In No-TL case, the optimal values in both versions $O1$ and $O2$  are higher than that of their respective timeline scenarios,  overestimating the realistic optimal value. Also, the realizable optimal objective value (TL case) further gets compromised when it is accompanied by adopting No-TL case optimizers. They may be sub-optimal for TL scenario for, e.g., in the context of $O2$ problem with No-TL case concluding $x^* \approx 2.27, \ \eta^* \approx 0.64$ (see Figure \ref{1Fig_UnNoTL}) to be optimizers is fallacious (the actual optimizers in the TL  case as in Figures \ref{1Fig_ConMeanLam_6} and \ref{1Fig_UnMeanLam_6} are $x^* \approx 2.08, \ \eta^* \approx 0.69$). Thus, ignoring TL can cause a CP to make sub-optimal decisions and may indicate  false trends.
 \vspace*{-1.5cm}

\section*{Conclusions}
We studied the propagation of a post of interest over a huge OSN.  We modeled the propagation of the post, considering the timeline (TL) structure, by an appropriate multi-type branching process. We found that the underlying branching process exhibits a certain dichotomy: either the post gets extinct or gets viral. We obtained various performance measures such as the time evolution of the number of unread posts, the expected number of shares, the probability of virality, etc. We showed that the expected number of shares grow at the same rate as the number of unread posts. We compare our results with the results that one would obtain without considering TL structure.   We discovered that without considering the TL structure, one leads to draw erroneous conclusions. For instance, we found that a study without TLs shows that even less attractive posts can get viral. It also indicates erroneous growth rates. 
More importantly, we also observe that without TL effects, one cannot capture some interesting paradigm shifts/phase transitions in certain behavioral patterns. For example, as the network becomes more active, one anticipates that it is more beneficial to engage in the network. The studies which do not incorporate these effects of TL lead to this erroneous conclusion; and argue that the virality chances increase monotonically as the mean number of friends increases ($m$). We demonstrated that virality chances do not increase monotonically with the number of friends. After a certain value of $m$, it decreases for some intermittently active networks (medium $m$ values).   To be more specific, for some range of parameters,  less active networks are preferable to more active networks.

Lastly, we integrated online auctions into our viral marketing model. We studied the optimization  problem considering the online auctions. We again compared the study with and without considering TL structure for varying activity levels of the network. We observe that the analysis without considering TL structure fails to capture phase transitions, thereby making the overall study incomplete.
Our study provides a framework using which, one can estimate important performance measures related to content propagation over online social networks, which further can be used in solving relevant optimization/game theoretic problems. 

\bibliographystyle{apacite}

\section*{Appendix}

{\large{\bf Proof of Lemma \ref{Lemma_positive_regular}:} }
i)  The matrix $e^{At}$ for any $t > 0$ is positive regular iff $e^A$ is (\cite{1AthreyaPaper}), because $A+I$ has only non-negative entries. Thus it is sufficient to prove $e^{A}$ is positive regular. Without loss of generality we can drop the multiplier $\lambda + \nu.$
Then the matrix $A$ can be written in the following way $A = A_1 + A_2 $, where 
\begin{eqnarray*}
 A_1 = \begin{bmatrix}
 c_1 r_1  & c_2r_1 + \theta & \cdot  & c_{N-1}r_1  & c_Nr_1 \\
 c_1 r_2   & c_2r_2 & \cdot & c_{N-1}r_2  & c_Nr_2 \\
  \cdot & \cdot & \cdot  & \cdot  & \cdot  \\
 c_1 r_{N-1}  & c_2r_{N-1}   & \cdot &  c_{N-1}r_{N-1} &  c_{N} r_{N-1} + \theta \\
 c_1 r_N  &  c_2r_N  &  \cdot &  c_{N-1}r_N  & c_{N} r_N 
\end{bmatrix} 
 \end{eqnarray*}
and  $A_2 = Diag(-1)$ is the diagonal matrix with `-1' on all the diagonals. Thus, $e^{A} = e^{A_1} e^{A_2} =  e^{-1} e^{A_1}$  since the matrices commute. For any $i$, one can express
\begin{eqnarray}
 e^{A_i} & = & I + A_i + \frac{A_i^2}{2!}  + \frac{A_i^3}{3!} + \cdots,  \label{ExpMatrix}
\end{eqnarray}
where $I$ is the identity matrix.
Also $e^{A_2} = e^{-1} I$ commutes with  $e^{A_1}.$
 A matrix is positive regular if there exists an $n$ such that $A^n$ has all positive entries. 
If $c_l > 0$ and $r_l > 0$ for all $l$, then $A_1$ is trivially positive regular  and hence $e^A$ is also positive regular. 

Consider a general case, where some of the constants can be zero, in particular consider the case with  $c_l = 0 \ \forall l > 1$ and $c_1 > 0$.   For this case: 
\begin{eqnarray*}
 A_1 = \begin{bmatrix}
 c_1 r_1  &  \theta & 0 & 0 & \cdot  & 0  & 0 & 0  & 0  \\
 c_1 r_2   &0  &\theta & 0 & \cdot & 0 & 0 & 0  & 0  \\
 c_1 r_3  & 0 & 0 &   \theta & \cdot  & 0  & 0 & 0  & 0  \\
 c_1 r_4   &0  &0  & 0 & \cdot & 0 & 0 & 0  & 0  \\
  \cdot & \cdot & \cdot  & \cdot  & \cdot  & \cdot  & \cdot & \cdot  & \cdot\\
     c_1 r_{N-3}  & 0 & 0 & 0 & \cdot & 0   &  \theta & 0 & 0 \\
   c_1 r_{N-2}  & 0 & 0 & 0 & \cdot & 0  & 0  &  \theta & 0 \\
 c_1 r_{N-1}  & 0 & 0 & 0 & \cdot & 0  & 0  &  0& \theta \\
 c_1 r_N  &  0 & 0 & 0 & \cdot & 0  & 0  & 0 & 0
\end{bmatrix}
\end{eqnarray*}

Then it is clear that  
\begin{eqnarray}
 A_1^2 = \begin{bmatrix}
 c_1^2 r_1^2 +\theta c_1 r_2  & \theta c_1r_1  &\theta^2 & \cdot  & 0& 0 \\
 c^2_1r_1 r_2  +\theta c_1 r_3 & \theta c_1 r_2 & 0& \cdot & 0& 0 \\
  \vdots & \vdots & \vdots  & \vdots  & \vdots  \\
c^2_1r_1 r_{N-2} +\theta c_1 r_{N-1}  & \theta c_1r_{N-2}  & 0 & \cdot & 0 &   \theta^2 \\
c^2_1r_1 r_{N-1} +\theta c_1 r_N  & \theta c_1r_{N-1}  & 0 & \cdot & 0 &  0 \\
 c_1^2 r_1 r_N  & \theta  c_1r_N & 0 &  \cdot &  0& 0 
\end{bmatrix}
\label{Eqn_A1_sqr}
 \end{eqnarray}
 The third power $A_1^3 = A_1^2 A_1$ will have first three columns positive because the first two columns in $A_1^2$ is have strict positive terms and the first $2\times 3$ sub matrix of $A_1$ 
\begin{eqnarray*}
 \begin{bmatrix} 
 c_1 r_1  &  \theta & 0    \\
 c_1 r_2   &0  &\theta  
 \end{bmatrix}
\end{eqnarray*}
 has at least one positive entry  in every column. 
 Continuing this way one can verify that $A_1^N$  has all positive entries by induction.    Basically once  $A_1^n$ has first $n$ columns with only  positive entries, because the first $n\times (n+1)$ sub-matrix of $A_1$ has atleast one positive entry in  every column, the matrix $A_1^{n+1}  = (A_1^n) \times A_1$ will have its first $n+1$ columns with only  positive entries.
   Further $A^n$ has only non negative entries for any $n \in \mathbb{N} $. From (\ref{ExpMatrix}) it is direct that $e^{A_1}$ is positive regular and so is $e^{-1}e^{A_1}$. 

 For the general case, when only some of $\{c_l\}$ are non-zero since terms are non-negative, the positive regularity follows from the above case and expansion (\ref{ExpMatrix}).   The result is true as long as   $c_1 > 0$.
 
{\bf   Proof of parts (ii)-(iii):}
 We proved that $e^{A} $ is positive regular. By Frobenius-Perron theory of positive regular matrices: a)there exists an eigenvalue, call it $e^\alpha$, of the matrix $e^{A}$ whose algebraic and geometric multiplicities are one and which dominates all the other eigenvalues in the absolute sense. In fact, $\alpha$ would be a real eigenvalue of matrix $A$, and  it  dominates the real components of all other eigenvalues of the matrix $A$; 
    b) there exists a left eigenvector $\textbf{u}$ and a right  eigenvector $\textbf{v}$, both  with all positive components,  corresponding to ${\alpha }$.  
   Fix one such set of left and right  eigenvectors $\textbf{u}$, $\textbf{v}$.  

Note that the eigenvectors of matrices $A$ and $e^{A}$ are the same.  
Any left eigenvector of $\alpha$, in particular  $\textbf{u}$,  satisfies   $\textbf{u}A = \alpha\textbf{u} $ and hence we get the following system of equations
relating  $\textbf{u}$ and $\alpha$
\begin{eqnarray}
(\lambda + \nu) c_1\textbf{r.u} - (\lambda + \nu)u_1  = \alpha u_1  \mbox{ or in other words }  c_1\textbf{r.u}   &=& \frac{\alpha +  \lambda + \nu}{\lambda + \nu} u_1,  \mbox{ and similarly}
\nonumber \\
c_l\textbf{r.u}  + \theta u_{l-1} & =&\frac{\alpha +  \lambda + \nu}{\lambda + \nu} u_l, \ \ l \geq 2 \label{SystemEqSIngleCP}.
\end{eqnarray}
Simplifying the above we obtain the following relation among various components of  left eigenvector  $\textbf{u}$:  for any $l \le N$
\begin{eqnarray}
\label{Eqn_ul_pattern}
u_l = \sum_{i = 0}^{l-1} \frac{\rho_{l-i}}{\rho_1} \left(\frac{\theta}{\sigma}\right)^{i} u_1; \ \ \ \sum_{i=1}^N  u_i =  \sum_{l=1}^N\frac{\rho_l}{\rho_1} \sum_{i=0}^{N-l} \left(\frac{\theta}{\sigma}\right)^{i}u_1 \mbox{ where } \sigma := \frac{\alpha +  \lambda + \nu}{ \lambda + \nu}.
 \end{eqnarray}  
Following exactly the same procedure, we obtain the relation among various components of   right  eigenvector  $\textbf{v}$ which are
\begin{eqnarray}
\label{Eqn_vl_pattern}
v_l = \sum_{i = 0}^{N-l} \frac{r_{l+i}}{r_N} \left(\frac{\theta}{\sigma}\right)^{i} v_N \ \forall l = 1,2, \cdots, N-1.
 \end{eqnarray}
This completes the proof of part (iii).

 Fix $\textbf{u}$, $\textbf{v}$  as before, and consider  the following linear  function of $\sigma'$:
\begin{eqnarray}
  P(\sigma') := \left(\textbf{r.c}\right) \textbf{r.u} + \theta\sum_{i = 1}^{N-1} r_{i+1} u_i  - \sigma' \textbf{r.u} \label{SingleCPEqaution}
\end{eqnarray} 
where $\textbf{r.c} := \sum_{i = 1}^N  r_i c_i$ etc.  
Multiplying either  side of  the equation (\ref{SystemEqSIngleCP})  with $r_l$ and then summing over $l$  we notice that  $\sigma $ is  a zero of $P(.)$.  
In other words,  eigenvalue $\alpha = (\sigma^*- 1)(\lambda+\nu)$, where $\sigma^*$  is a   zero of $P(.)$ 
Because   $u_i >0$ for all $l$,  $\textbf{r.u} >0$ and similarly  $\textbf{r.c} > 0$. Thus  $\sigma$  is the only   zero of $P(.)$. 
It is clear that
$$P\left( \textbf{r.c} \right)= \theta\sum_{i = 1}^{N-1} r_{i+1} u_i > 0. $$ 
Since $r_i$s are monotonic, i.e., because $r_1 \ge r_2  \ge \cdots  \ge r_N$, 
 $$P\left(\textbf{r.c} + \theta  \right)=   \theta\sum_{i = 1}^{N-1} r_{i+1} u_i  - \theta \textbf{r.u}  < 0. $$  
 Thus, the only zero of  $P(.)$ lies in the  open interval interval $\Big (\textbf{r.c}, \ \ \textbf{r.c} + \theta \Big )$.  Thus 
 $\alpha \in \Big( \textbf{r.c}-1, \ \ \textbf{r.c} + \theta -1\Big )(\lambda+\nu).$
 
 Consider the special case with  $r_l = d_1d_2^l$, where $d_1$ and $ d_2 \le 1$ are constants, then clearly the only root of equation (\ref{SingleCPEqaution}) $\sigma$ equals 
\begin{eqnarray*}
\sigma = \textbf{r.c}+ \theta d_2\dfrac{\sum_{i = 1}^{N-1} r_{i} u_i}{\textbf{r.u}  } =   \textbf{r.c}+ \theta d_2\left(1-\frac{r_N u_N}{\textbf{r.u} }\right) . \end{eqnarray*}  
Now we  study the convergence of $\sigma$ as $N \to \infty$.
It is obvious that the eigenvectors/eigenvalues corresponding to different $N$ would be different.  We would normalize them by choosing the eigenvector ${\bf u}$
 with $u_1 = 1$ for any $N$.  
 With such a choice,  it is clear from (\ref{Eqn_ul_pattern}) that   $u_N$ remains bounded even when we let $N \to \infty$. 
 Thus as   $N \to \infty$
\begin{eqnarray*}
\sigma =  \textbf{r.c}+ \theta d_2\left(1-\frac{r_N u_N}{\textbf{r.u} }\right) \to \textbf{r.c}+ \theta d_2 \ \ \text{as \ } N \to \infty \ \ \ \because \left(r_N \to 0\right).
\end{eqnarray*}  
Thus,  as the number of TL levels increase the largest eigenvalue, $\alpha$ of matrix $A$ converges to $\left(\textbf{r.c}+ \theta d_2 -1 \right)(\lambda + \nu) $. 
\eop
 \\
 \\ \\
 \
{\large{\textbf{ Computation of $v_l \sum_l u_l:$}}}
Referring to Theorem \ref{Theorem_virality}, the left and right eigenvectors of the matrix $A$  are
 \begin{eqnarray*}
u_l = \sum_{i = 0}^{l-1} \frac{\rho_{l-i}}{\rho_1} \Big(\frac{\theta}{\sigma}\Big)^i u_1, \ \ v_l = \sum_{i = 0}^{N-l} \frac{r_{l+i}}{r_N} \Big(\frac{\theta}{\sigma}\Big)^i v_N; 
\ \  \ l \geq 2 \ \  \mbox{and \ } \sigma =\frac{\alpha}{(\lambda + \nu)} +1.
\end{eqnarray*} 
When $\rho_i = \tilde{\rho} \rho^i, r_i = d_1d_2^i$ with  $0 < d_1,d_2, \rho < 1.$ On substituting these values, we obtain 
\begin{eqnarray*}
u_l & = & \frac{1 -  \Big(\frac{\theta}{\sigma \rho}\Big)^l }{1-  \frac{\theta}{\sigma \rho }} \rho^{l-1} u_1
 = \frac{\rho^{l-1} -  \frac{1}{\rho}\Big(\frac{\theta}{\sigma }\Big)^l }{1-  \frac{\theta}{\sigma \rho}}  u_1;
\\
\sum_{l =1}^{N} u_l & = &
 \frac{u_1}{ \Big( 1- \frac{\theta}{\sigma \rho} \Big)}   \left(\frac{1-\rho^N}{1-\rho} - \frac{\theta}{\sigma \times \rho}  \frac{1 -\big(\frac{\theta}{\sigma}\big)^N}{1-\frac{\theta}{\sigma}} \right) \nonumber
\\[.2cm]
v_l & = &  \frac{1- \Big(\frac{\theta d_2}{\sigma}\Big)^{N-l+1} }  {1 - \frac{\theta d_2}{\sigma }}d_2^{l-N}  v_N 
\ = \   \frac{\frac{d_2^l}{d_2^N}- \Big(\frac{\theta }{\sigma}\Big)^{N -l+1} d_2^{N-l+1 + l-N}}  {1 - \frac{\theta d_2}{\sigma }}  v_N \nonumber 
\\
& =&  \frac{\frac{d_2^l}{d_2^N} - \Big(\frac{\theta }{\sigma}\Big)^{N+1} d_2 \Big(\frac{\sigma}{\theta}\Big)^l}  {1 - \frac{\theta d_2}{\sigma }}  v_N. \nonumber
\end{eqnarray*}

\begin{eqnarray}
 v_l \sum_l u_l & = & \frac{\frac{d_2^l}{d_2^N}- \Big(\frac{\theta }{\sigma}\Big)^{N+1} d_2 \Big(\frac{\sigma}{\theta}\Big)^l}  {1 - \frac{\theta d_2}{\sigma }}  v_N 
 \frac{u_1}{ \Big( 1- \frac{\theta}{\sigma \rho} \Big)}   \left(\frac{1-\rho^N}{1-\rho} - \frac{\theta}{\sigma \times \rho}  \frac{1 -\big(\frac{\theta}{\sigma}\big)^N}{1-\frac{\theta}{\sigma}} \right) \nonumber
  \\[0.2cm]
 & = & \frac{\frac{d_2^l}{d_2^N}- \Big(\frac{\theta }{\sigma}\Big)^{N+1} d_2 \Big(\frac{\sigma}{\theta}\Big)^l}
  {\Big(1 - \frac{\theta d_2}{\sigma }\Big)\Big(1 - \frac{\theta }{\sigma }\Big)} \left(\frac{1-\rho^N}{1-\rho} - \frac{\theta}{\sigma \times \rho}  \frac{1 -\big(\frac{\theta}{\sigma}\big)^N}{1-\frac{\theta}{\sigma}} \right) u_1 v_N. \label{Eqn_KaamKi}
    \end{eqnarray}
We require the value of $u_1v_N$ towards obtaining $v_l \sum_l u_l.$ For this, we will use  the fact that $\textbf{u.v} = 1$. So,
\begin{eqnarray*}
\sum_{l=1}^N u_l v_l & = & \sum_{l=1}^N \frac{\rho^{l-1} -  \frac{1}{\rho}\Big(\frac{\theta}{\sigma }\Big)^l }{1-  \frac{\theta}{\sigma\rho }}  \frac{\frac{d_2^l}{d_2^N}- \Big(\frac{\theta }{\sigma}\Big)^{N+1} d_2 \Big(\frac{\sigma}{\theta}\Big)^l}  {1 - \frac{\theta d_2}{\sigma }} u_1v_N = 1.
\end{eqnarray*}
Observe that
\begin{eqnarray*}
\sum_{l=1}^N  \left(\rho^{l-1} -  \frac{1}{\rho}\Big(\frac{\theta}{\sigma }\Big)^l \right) \times \left( \frac{d_2^l}{d_2^N}- \Big(\frac{\theta }{\sigma}\Big)^{N+1} d_2 \Big(\frac{\sigma}{\theta}\Big)^l\right)
&  & 
\\
&  & \hspace*{-7cm} = \ \frac{d_2}{d_2^N} \sum_{l=1}^N (d_2 \rho)^{l-1}  - \ \ \frac{1}{\rho d_2^N} \sum_{l=1}^N
\Big(\frac{\theta d_2}{\sigma }\Big)^l
- \sum_{l=1}^N \Big(\frac{\theta}{\sigma }\Big)^{N+1}\frac{d_2}{\rho} \Big(\frac{\sigma \rho}{\theta} \Big)^l + \sum_{l=1}^N \frac{d_2}{\rho} \Big(\frac{\theta}{\sigma }\Big)^{N+1}
\\
&  & \hspace*{-7cm}  = \frac{d_2}{d_2^N} \frac{1-(d_2\rho)^N}{1-d_2\rho} 
- \frac{1}{\rho d_2^N}\frac{\theta d_2}{\sigma } \frac{1-\Big(\frac{\theta d_2}{\sigma }\Big)^N}{1-\frac{\theta d_2}{\sigma }}
- \frac{d_2}{\rho} \Big(\frac{\theta}{\sigma }\Big)^{N+1} \frac{\sigma \rho}{\theta  } \frac{1-\Big(\frac{\sigma \rho}{\theta  }\Big)^N}{1-\frac{\sigma \rho}{\theta  }} + \frac{d_2}{\rho} \Big(\frac{\theta}{\sigma }\Big)^{N+1} N.
\end{eqnarray*} 
Substituting this value back to $\sum_{l=1}^N u_l v_l$, we get
{\footnotesize{ 
\begin{eqnarray}
&&\frac{u_1v_N}{\Big(1-  \frac{\theta}{\sigma\rho}\Big)\Big(1-  \frac{\theta d_2}{\sigma}\Big)} 
\left(\frac{d_2}{d_2^N} \frac{1-(d_2\rho)^N}{1-d_2\rho} 
- \frac{1}{\rho d_2^N}\frac{\theta d_2}{\sigma } \frac{1-\Big(\frac{\theta d_2}{\sigma }\Big)^N}{1-\frac{\theta d_2}{\sigma }}
- \frac{d_2}{\rho} \Big(\frac{\theta}{\sigma }\Big)^{N+1} \frac{\sigma \rho}{\theta  } \frac{1-\Big(\frac{\sigma \rho}{\theta  }\Big)^N}{1-\frac{\sigma \rho}{\theta  }} + \frac{d_2}{\rho} \Big(\frac{\theta}{\sigma }\Big)^{N+1} N \right) =1 \nonumber
\\
&& u_1v_N \ = \dfrac{  \Big(1-  \frac{\theta}{\sigma\rho}\Big)\Big(1-  \frac{\theta d_2}{\sigma}\Big)}
{\left(\frac{d_2}{d_2^N} \frac{1-(d_2\rho)^N}{1-d_2\rho} 
- \frac{1}{\rho d_2^N}\frac{\theta d_2}{\sigma } \frac{1-\Big(\frac{\theta d_2}{\sigma }\Big)^N}{1-\frac{\theta d_2}{\sigma }}
- \frac{d_2}{\rho} \Big(\frac{\theta}{\sigma }\Big)^{N+1} \frac{\sigma \rho}{\theta  } \frac{1-\Big(\frac{\sigma \rho}{\theta  }\Big)^N}{1-\frac{\sigma \rho}{\theta  }} + \frac{d_2}{\rho} \Big(\frac{\theta}{\sigma }\Big)^{N+1} N \right) } \label{Eqnu1vNVal}
\end{eqnarray}}}

  Now substituting  the value of $u_1 v_N $ in the equation \ref{Eqn_KaamKi}
  
\begin{eqnarray*}
 v_l \sum_l u_l & = &  \frac{\frac{d_2^l}{d_2^N}- \Big(\frac{\theta }{\sigma}\Big)^{N+1} d_2 \Big(\frac{\sigma}{\theta}\Big)^l}
  {\Big(1 - \frac{\theta d_2}{\sigma }\Big)\Big(1 - \frac{\theta }{\sigma }\Big)} \left(\frac{1-\rho^N}{1-\rho} - \frac{\theta}{\sigma \times \rho}  \frac{1 -\big(\frac{\theta}{\sigma}\big)^N}{1-\frac{\theta}{\sigma}} \right) 
  \\
  && \hspace{1cm} \times \dfrac{  \Big(1-  \frac{\theta}{\sigma\rho}\Big)\Big(1-  \frac{\theta d_2}{\sigma}\Big)}
{\left(\frac{d_2}{d_2^N} \frac{1-(d_2\rho)^N}{1-d_2\rho} 
- \frac{1}{\rho d_2^N}\frac{\theta d_2}{\sigma } \frac{1-\Big(\frac{\theta d_2}{\sigma }\Big)^N}{1-\frac{\theta d_2}{\sigma }}
- \frac{d_2}{\rho} \Big(\frac{\theta}{\sigma }\Big)^{N+1} \frac{\sigma \rho}{\theta  } \frac{1-\Big(\frac{\sigma \rho}{\theta  }\Big)^N}{1-\frac{\sigma \rho}{\theta  }} + \frac{d_2}{\rho} \Big(\frac{\theta}{\sigma }\Big)^{N+1} N \right) }
\end{eqnarray*}
{\normalsize{
\begin{eqnarray*}
& = &
\frac{\frac{d_2^l}{d_2^N}- \Big(\frac{\theta }{\sigma}\Big)^{N+1} d_2 \Big(\frac{\sigma}{\theta}\Big)^l}
  {\Big(1 - \frac{\theta }{\sigma }\Big)} \left(\frac{1-\rho^N}{1-\rho} - \frac{\theta}{\sigma \rho}  \frac{1 -\big(\frac{\theta}{\sigma}\big)^N}{1-\frac{\theta}{\sigma}} \right) 
  \\
  && \hspace{2cm} \times
\ \ \dfrac{  \Big(1-  \frac{\theta}{\sigma\rho}\Big)}
{\left(\frac{d_2}{d_2^N} \frac{1-(d_2\rho)^N}{1-d_2\rho} 
- \frac{1}{\rho d_2^N}\frac{\theta d_2}{\sigma } \frac{1-\Big(\frac{\theta d_2}{\sigma }\Big)^N}{1-\frac{\theta d_2}{\sigma }}
- \frac{d_2}{\rho} \Big(\frac{\theta}{\sigma }\Big)^{N+1} \frac{\sigma \rho}{\theta  } \frac{1-\Big(\frac{\sigma \rho}{\theta  }\Big)^N}{1-\frac{\sigma \rho}{\theta  }} + \frac{d_2}{\rho} \Big(\frac{\theta}{\sigma }\Big)^{N+1} N \right) }
\\
& = & \frac{d_2^l - \Big(\frac{\theta d_2}{\sigma}\Big)^{N+1}  \Big(\frac{\sigma}{\theta}\Big)^l}
  {\Big(1 - \frac{\theta }{\sigma }\Big)} \left(\frac{1-\rho^N}{1-\rho} - \frac{\theta}{\sigma \rho}  \frac{1 -\big(\frac{\theta}{\sigma}\big)^N}{1-\frac{\theta}{\sigma}} \right) 
    \\
  && \hspace{2cm} \times
\ \ \dfrac{  \Big(1-  \frac{\theta}{\sigma\rho}\Big)}
{\left(d_2 \frac{1-(d_2\rho)^N}{1-d_2\rho} 
- \frac{\theta d_2}{\sigma \rho } \frac{1-\Big(\frac{\theta d_2}{\sigma }\Big)^N}{1-\frac{\theta d_2}{\sigma }}
- \frac{1}{\rho} \Big(\frac{\theta d_2}{\sigma }\Big)^{N+1} \frac{\sigma \rho}{\theta  } \frac{1-\Big(\frac{\sigma \rho}{\theta  }\Big)^N}{1-\frac{\sigma \rho}{\theta  }} + \frac{1}{\rho} \Big(\frac{\theta d_2}{\sigma }\Big)^{N+1} N \right) }.
\end{eqnarray*} }}  
On simplifying and using property of limits:
{\normalsize{
\begin{eqnarray}
\lim_{ N \to \infty}   v_l \sum_{l =1}^{N} u_l  & = & \lim_{ N \to \infty}  \frac{d_2^l - \Big(\frac{\theta d_2}{\sigma}\Big)^{N+1}  \Big(\frac{\sigma}{\theta}\Big)^l}
  {\Big(1 - \frac{\theta }{\sigma }\Big)} \left(\frac{1-\rho^N}{1-\rho} - \frac{\theta}{\sigma  \rho}  \frac{1 -\big(\frac{\theta}{\sigma}\big)^N}{1-\frac{\theta}{\sigma}} \right) \nonumber
\\
&&\hspace*{-2.5cm} \times \ \lim_{ N \to \infty}  \dfrac{  \Big(1-  \frac{\theta}{\sigma\rho}\Big)}
{\left(d_2 \frac{1-(d_2\rho)^N}{1-d_2\rho} 
- \frac{\theta d_2}{\sigma \rho } \frac{1-\Big(\frac{\theta d_2}{\sigma }\Big)^N}{1-\frac{\theta d_2}{\sigma }}
-   d_2 \frac{ \Big(\frac{\theta d_2}{\sigma }\Big)^{N} -\big( d_2 \rho\big)^N}{1-\frac{\sigma \rho}{\theta  }} + \frac{1}{\rho} \Big(\frac{\theta d_2}{\sigma }\Big)^{N+1} N \right) }. \label{lims}
\end{eqnarray}}}
Note that in viral scenario $\alpha > 0$ and  $\sigma > 1$ and hence $\theta/\sigma < 1$, so \begin{eqnarray}
\lim_{ N \to \infty}   \Big(\frac{\theta }{\sigma}\Big)^{N+1} = 0, \ \  \ \lim_{ N \to \infty}\Big(\frac{\theta d_2}{\sigma }\Big)^{N+1}   = 0, \ \ \ \lim_{ N \to \infty}  N \Big(\frac{\theta d_2}{\sigma}\Big)^{N+1}  = 0 \ \because \ \  d_2 <1. 
\end{eqnarray}
 In what follows, for any fixed $l$, we have
{\footnotesize{
\begin{eqnarray}
\lim_{ N \to \infty}  \frac{d_2^l - \Big(\frac{\theta d_2}{\sigma}\Big)^{N+1}  \Big(\frac{\sigma}{\theta}\Big)^l}
  {\Big(1 - \frac{\theta }{\sigma }\Big)} \left(\frac{1-\rho^N}{1-\rho} - \frac{\theta}{\sigma  \rho}  \frac{1 -\big(\frac{\theta}{\sigma}\big)^N}{1-\frac{\theta}{\sigma}} \right) & = &
 \frac{d_2^l}
  {\Big(1 - \frac{\theta }{\sigma }\Big)} \left(\frac{1}{1-\rho} - \frac{\theta}{\sigma  \rho}  \frac{1 }{1-\frac{\theta}{\sigma}} \right) \hspace*{1cm}\label{EqnPhliLim}
\\
\lim_{ N \to \infty}  \dfrac{  \Big(1-  \frac{\theta}{\sigma\rho}\Big)}
{\left(d_2 \frac{1-(d_2\rho)^N}{1-d_2\rho} 
- \frac{\theta d_2}{\sigma \rho } \frac{1-\Big(\frac{\theta d_2}{\sigma }\Big)^N}{1-\frac{\theta d_2}{\sigma }}
-   d_2 \frac{ \Big(\frac{\theta d_2}{\sigma }\Big)^{N} -\big( d_2 \rho\big)^N}{1-\frac{\sigma \rho}{\theta  }} + \frac{1}{\rho} \Big(\frac{\theta d_2}{\sigma }\Big)^{N+1} N \right) } & = &
 \dfrac{  \Big(1-  \frac{\theta}{\sigma\rho}\Big)}
{\left( \frac{d_2}{1-d_2\rho} 
- \frac{\theta d_2}{\sigma \rho } \frac{1}{1-\frac{\theta d_2}{\sigma }} \right)}  \label{EqnDusriLim}
\end{eqnarray}}}
Substituting these limits (equations \ref{EqnPhliLim} and \ref{EqnDusriLim}) in equation \ref{lims}, we get
\begin{eqnarray*}
\lim_{ N \to \infty}   v_l \sum_{l =1}^{N} u_l  & = & \frac{d_2^l}
  {\Big(1 - \frac{\theta }{\sigma }\Big)} \left(\frac{1}{1-\rho} - \frac{\theta}{\sigma  \rho}  \frac{1 }{1-\frac{\theta}{\sigma}} \right)  \dfrac{  \Big(1-  \frac{\theta}{\sigma\rho}\Big)}
{\left( \frac{d_2}{1-d_2\rho} 
- \frac{\theta d_2}{\sigma \rho } \frac{1}{1-\frac{\theta d_2}{\sigma }} \right)} 
\\
& = &
\frac{d_2^{l-1}}
  {\Big(1 - \frac{\theta }{\sigma }\Big)} \left(\frac{1}{1-\rho} - \frac{\theta}{  \rho}  \frac{1 }{\sigma-\theta} \right)  \dfrac{  \Big(1-  \frac{\theta}{\sigma\rho}\Big)}
{\left( \frac{1}{1-d_2\rho} 
- \frac{\theta }{\rho } \frac{1}{\sigma -\theta d_2} \right)} 
\end{eqnarray*}
Thus, we have
\begin{eqnarray}
 v_l \sum_l u_l & = &   
d_2^{l-1} (1-d_2 \rho) \left(\frac{1}{1-\rho} - \frac{\theta}{\rho }  \frac{1}{\sigma -\theta} \right) \frac{(\sigma - \theta d_2)(\sigma \rho -\theta)}{(\sigma -\theta) (\rho -\theta)} \label{Eqn_vi_sum_u}
\end{eqnarray}
\eop
\\
\\ 
{\large{\bf Proof of Lemma \ref{Lemma_Shares} }}:
  Let $\bf{{\bf j}_{\bf x}}  = \{  j_{x_1},j_{x_2},\cdots,j_{x_N} \} $ be the number of TLs of type $1,2,\cdots,N$ respectively, and $y$ be the total number of shares. It is easy to observe that  $y \ge \sum_{i} j_{x_i}$. We write it in short form as $y \ge {\bf j}_{\bf x}$.  Define ${\bf s}_{\bf x}^{{\bf j}_{\bf x}} := \Pi_{i}  s_{x_i}^{j_{x_i}}.$ Then the PGF of TL-CTBP,  when started with one type-1 particle,  can be written as\footnote{ $P_{(\textbf{e}_1, 1) \to ({\bf j}_{\bf x}, y)} \big(t\big)$ is the probability that state $(\bf{e}_1, 1) $  (one type-1 particle and 1 total progeny/number of shares)  after time $t$ gets  transformed to population vector ${\bf j}_{\bf x}$ and number of shares $y$ } 
\begin{eqnarray*}
F_1(\textbf{s}, t) = \sum_{{\bf j}_{\bf x} = 0}^{\infty}\sum_{ \substack{y \ge  {\bf j}_{\bf x} }}^{\infty} P_{(\textbf{e}_1, 1) \to ({\bf j}_{\bf x}, y)} \big(t\big){\bf s}_{\bf x}^{{\bf j}_{\bf x}}s_{y}^{y} \ \ \text{and, }
\\
\frac{\delta F_1(\textbf{s}, t)}{\delta t}= \sum_{{\bf j}_{\bf x} = 0}^{\infty}\sum_{ \substack{y \ge  {\bf j}_{\bf x} }}^{\infty} P'_{(\textbf{e}_1, 1) \to ({\bf j}_{\bf x}, y)} \big(t\big){\bf s}_{\bf x}^{{\bf j}_{\bf x}} s_{y}^{y}.
\end{eqnarray*}
This is obtained by conditioning on the events of the first transition. Note that the populations generated by two parents evolve independently of each other and the procedure is similar to the standard procedure used in these kind of computations (e.g., \cite{1AthreyaBook}).
Let $\pmb{\xi} = \big( \xi_1, \xi_2, \cdots, \xi_N \big)$ represent the offspring produced by one parent of type-1 and let $\sxi := \sum_{i} \xi_i$. 

By backward equation, we have $P'_{1k}(t) = \sum_{j}q_{1j}P_{jk}(t)$; in our case it is
\begin{eqnarray*}
\frac{\delta F_1(\textbf{s}, t)}{\delta t} & = & \big(\lambda +\nu \big) \Big((1-\theta) r_1 \sum_{\pmb{\xi} }\sum_{{\bf j}_{\bf x} = 0}^{\infty}\sum_{ \substack{y \ge  {\bf j}_{\bf x} }}^{\infty} P_1(\pmb{\xi})P_{\big( \pmb{\xi}, \sxi +  1 \big)\to ({\bf j}_{\bf x}, y)} \Big( t \Big) {\bf s}_{\bf x}^{{\bf j}_{\bf x}}s_{y}^{y} 
+  \theta F_2(\textbf{s}, t)  
\\
& - & \sum_{{\bf j}_{\bf x} = 0}^{\infty}\sum_{ \substack{y \ge  {\bf j}_{\bf x} }}^{\infty} P_{(\textbf{e}_1, 1) \to ({\bf j}_{\bf x}, y)} \big(t\big){\bf s}_{\bf x}^{{\bf j}_{\bf x}}s_{y}^{y} + (1-\theta)(1-r_1)s_{y} \Big)
\\
\frac{\delta F_1(\textbf{s}, t)}{\delta t} & = &  \big(\lambda +\nu \big)\Big( (1-\theta) r_1 \sum_{\pmb{\xi} } P_1(\pmb{\xi}) \Pi_{i =1}^{N}\Big(  \sum_{{\bf j}_{\bf x} = 0}^{\infty}\sum_{ \substack{y \ge  {\bf j}_{\bf x} }}^{\infty} P_{(\textbf{e}_i, 1) \to ({\bf j}_{\bf x}, y)} \big(t\big){\bf s}_{\bf x}^{{\bf j}_{\bf x}}s_{y}^{y} \Big)^{\xi_i} s_y
\\
& + &
  \theta F_2(\textbf{s}, t)   -  F_1(\textbf{s}, t)+(1-\theta)(1-r_1)s_y\Big)
\\
\frac{\delta F_1(\textbf{s}, t)}{\delta t} & = & \big(\lambda +\nu \big)\Big( (1-\theta) r_1 s_y f_1\Big(\textbf{F}(\textbf{s}, t) \Big) +\theta F_2(\textbf{s}, t) 
- F_1(\textbf{s}, t)+(1-\theta)(1-r_1)s_y \Big)
\end{eqnarray*}
where $\textbf{F}(\textbf{s}, t) := \{ F_1(\textbf{s}, t),F_2(\textbf{s},t), \cdots, F_N(\textbf{s},t)\}.$ Similarly we can write  for any $l$
\begin{eqnarray*}
\frac{\delta F_l(\textbf{s}, t)}{\delta t} & = & \big(\lambda +\nu \big)\Big( (1-\theta) r_l s_y f_l\Big(\textbf{F}(\textbf{s}, t) \Big) + \theta \Big( \mathbb{1}_{l<N} F_{l+1}(\textbf{s}, t)   +s_y \mathbb{1}_{l = N}\Big)
\\
&- & F_l(\textbf{s}, t)+(1-\theta)(1-r_l)s_y \Big).
\end{eqnarray*} 
 Let $ \dot{y}_{l}(t) = \frac{\delta^2 F_l(\textbf{s},t)}{\delta t\delta s_y}\vert_{\textbf{s}=	1} \ \forall l =\{1,2,\cdots, N\}$ represent the time derivative of number shares till time $t$ when started with a type $l$ progenitor. We have the following expression 
\begin{eqnarray*}
\dot{y}_{1}(t) & = & \Big(\lambda +\nu \big) \Big( (1-\theta) r_1 f_1(1) + (1-\theta) r_1 \sum_{i=1}^{N}\frac{\delta f_1(\textbf{F}(\textbf{s},t) }{\delta F_i(\textbf{s},t)} \frac{\delta F_i(\textbf{s},t)}{\delta s_y}\Big{\vert_{\textbf{s}=	1}} + (1-\theta)(1-r_1)1 
\\
& + & \theta  \frac{\delta F_2(\textbf{s},t)}{\delta s_y}\Big{\vert_{\textbf{s}=	1}}  
-  \frac{\delta F_1(\textbf{s},t)}{\delta s_y}\Big{\vert_{\textbf{s}=	1}}  \Big )
\\
 & = & \big(\lambda +\nu \big) \Big( (1-\theta) r_1  + (1-\theta)(1-r_1) + (1-\theta) r_1 m\eta\sum_{i=1}^{N} \rho_i  y_{i}(t) + \theta  y_{2}(t)  -  y_{1}(t) \Big)
\\
 & = & \big(\lambda +\nu \big) \Big( 1-\theta + r_1\sum_{i=1}^{N} c_i  y_{i}(t)  + \theta y_{2}(t) -  y_{1}(t) \Big). 
\end{eqnarray*}
Similarly, we can write the above for any $l$
\begin{equation}
\dot{y}_{l}(t) = \big(\lambda +\nu \big) \Big(1- \theta +  r_l\sum_{i=1}^{N} c_i  y_{i}(t) +  \theta  y_{l+1}(t)  \mathbb{1}_{l<N} - y_{l}(t) +  \theta \mathbb{1}_{l =N}\Big).
\end{equation}
The above can be written  in matrix form as
{\footnotesize{
\begin{eqnarray*}
\frac{1}{\lambda + \nu}  \left[
\begin{array}{ccccc}
\dot{y}_{1}(t) \\ \dot{y}_{2}(t) \\ \vdots \\ \dot{y}_{N-1 }(t) \\ \dot{y}_{N}(t)
\end{array}
\right]
&= &
 \left[ \begin{array}{ccccccc}
 c_1 r_1 -1 & c_2r_1 + \theta & \cdots  & c_{N-1}r_1  & c_Nr_1 \\
 c_1 r_2   & c_2r_2 -1 & \cdots & c_{N-1}r_2  & c_Nr_2 \\
  & \vdots    \\
 c_1 r_{N-1}  & c_2r_{N-1}   & \cdots &  c_{N-1}r_{N-1} -1 &  c_{N} r_{N-1} + \theta \\
 c_1 r_N  &  c_2r_N  &  \cdots &  c_{N-1}r_N  & c_{N} r_N -1 \\
 \end{array} 
 \right] \left[
\begin{array}{ccccc}
y_{1}(t) \\ y_{2}(t) \\ \vdots \\ y_{N-1 }(t) \\ y_{N}(t)
\end{array}
\right]
+
\hspace{0.1cm}
 \left[
\begin{array}{c}
1-\theta  \\ 1-\theta  \\ \vdots \\ 1-\theta  \\ 1
\end{array}
\right].
\end{eqnarray*}}}
 Solving the above set of equations, we obtain:
\begin{eqnarray}
\left[
\begin{array}{c}
y_{1}(t) \\ y_{2}(t) \\ \vdots \\ y_{N-1 }(t) \\ y_{N}(t)
\end{array}
\right] 
&= &
e^{At}\left[
\begin{array}{c}
y_{1}(0) \\ y_{2}(0) \\ \vdots \\ y_{N-1 }(0) \\ y_{N}(0)
\end{array}
\right]  + e^{At} \int_{0}^{t}e^{-As}  (\lambda + \nu)\left[
\begin{array}{c}
1-\theta  \\ 1-\theta  \\ \vdots \\ 1-\theta  \\ 1
\end{array}
\right] ds
\\ 
&= &
e^{At}\left[
\begin{array}{c}
y_{1}(0) \\ y_{2}(0) \\ \vdots \\ y_{N-1 }(0) \\ {y}_{N}(0)
\end{array}
\right]  + e^{At} A^{-1}\big( I -e^{-At} \big) (\lambda + \nu)\left[
\begin{array}{c}
1-\theta  \\ 1-\theta  \\ \vdots \\ 1-\theta  \\ 1
\end{array}
\right]  
\end{eqnarray}
 
With ${\bf y} (t) := \{y_{1}(t), y_{2}(t), \cdots, y_{N}(t)\}$, we can represent the above as:
\begin{eqnarray*}
{\bf y} (t) := \left[
\begin{array}{c}
y_{1}(t) \\ y_{2}(t) \\ \vdots \\ y_{N-1 }(t) \\ y_{N}(t)
\end{array}
\right] 
&= &
e^{At} \Big(  {\mbox{ $ \left [\begin{array}{cccc}
y_1(0)  \\
\vdots \\
y_N(0)  
\end{array} \right ] $  } }  + (\lambda + \nu) A^{-1} \textbf{k} \Big) - (\lambda + \nu)A^{-1}\textbf{k}
\end{eqnarray*}  
where $\textbf{k}  = [1-\theta, 1-\theta, \cdots, 1-\theta,1 ]^T.$
\begin{eqnarray*}
\textbf{y}(t) & = & e^{At} \Big( \textbf{1} + (\lambda + \nu) A^{-1} \textbf{k} \Big) - (\lambda + \nu)A^{-1}\textbf{k}
\end{eqnarray*}
 From \cite[equation (45)]{1AthreyaPaper}, $e^{At}$ can be approximated  for large $t$.  By which, we can write
  \begin{eqnarray}
  \label{Eqn_yt2}
  {\bf y} (t)
    & \approx&  e^{\alpha t}  {\bf v u'}\Big( \textbf{1} + (\lambda + \nu) A^{-1} \textbf{k} \Big) - (\lambda + \nu)A^{-1}\textbf{k}
    \\    
 & \approx&  e^{\alpha t}  \Big(\textbf{v}\sum_{i = 1}^N u_i + \frac{\lambda + \nu} {\alpha}\textbf{vu}'\textbf{k}\Big)  - (\lambda + \nu)A^{-1}\textbf{k} \nonumber
\\
 & \approx&  \textbf{v} e^{\alpha t} \left ( \sum_i  u_i   \Big ( 1 + \frac{\lambda + \nu} {\alpha} (1-\theta)\Big)  + \frac{\lambda + \nu} {\alpha} u_N   \right )
- (\lambda + \nu)A^{-1}\textbf{k} \nonumber
 \\
& \approx &    \textbf{v} e^{\alpha t} \sum_i  u_i     \left (  1 +  \frac{1-\theta} {{\bf r.c} -1 + \theta d_2}       \right )
- (\lambda + \nu)A^{-1}\textbf{k}.\nonumber
\end{eqnarray}  
\eop 
\\
\\
{\large{\textbf{Proof of Proposition \ref{Prop_BgtTightness}}:}}
We will prove that at optimality the budget constraint is tight, i.e., $ x + \kappa_1 \eta = \bar{B}$ using the Lagrangian relaxation method. To do so, we first change the inequality budget constraint, $ x + \kappa_1 \eta \le \bar{B} $ , to equality constraint as follows. Let $s^2$ (ensuring it to be $\ge 0$) be slack variable  such that $ x + \kappa_1 \eta + s^2 = \bar{B}.$  Similarly, we have 
\begin{eqnarray*}
\eta - s_1^2 =\bar{\eta}, \ \eta  + s_2^2 = 1, \ \ x -s_3^2 = 0 \ \ \mbox{for the constraints } \ \eta  \ge \bar{\eta}, \ \ \eta \le 1 \ \ \mbox{respectively}
\end{eqnarray*}
where $s_1^2, s_2^2, s_3^2$ are the slack/surplus variables.
The Lagrangian function $\mathit{L}(x,\eta, \Lambda) $ with Lagrangian multiplier $\Lambda, \Lambda_1,\Lambda_2,\Lambda_3$ is given as
{\footnotesize
\begin{eqnarray*}
\max_{x,\eta} \ \  \underbrace{\log E\Big(\sum_{l} X_l(t)\Big)P( Bid \leq x \eta) -\Lambda \left(\bar{B} -x - \kappa_1 \eta - s^2\right) - \Lambda_1 \left( \bar{\eta} -\eta + s_1^2  \right) - \Lambda_2 \left(1-\eta  - s_2^2 \right) - \Lambda_3(x-s_3^2)}_{\mathit{L}(x,\eta, \Lambda,\Lambda_1,\Lambda_2,\Lambda_3)}. 
\end{eqnarray*}}
The critical points of  $ \log E\Big(\sum_{l} X_l(t)\Big)P( Bid \leq x \eta)$ with the given constraint, say 
$${\cal G}:= -x - \kappa_1 \eta - s^2, \ \ {\cal G}_1:= \bar{\eta} - \eta + s_1^2,  \ \ {\cal G}_2:= 1-\eta  - s_2^2 , \ \ {\cal G}_3:= x-s_3^2 $$
 are obtained by solving the following system of simultaneous equations\footnote{see \url{http://users.wpi.edu/~pwdavis/Courses/MA1024B10/1024_Lagrange_multipliers.pdf}} \footnote{and \url{http://www.math.harvard.edu/archive/21a_spring_09/PDF/11-08-Lagrange-Multipliers.pdf}}

\begin{eqnarray}
&& \hspace{-2cm}\frac{\partial\log E\Big(\sum_{l} X_l(t)\Big)P( \mathbf{B} < x \eta) }{\partial x} = \Lambda\frac{\partial {\cal G}}{\partial x} + \Lambda_1\frac{\partial {\cal G}_1}{\partial x} +
\Lambda_2\frac{\partial {\cal G}_2}{\partial x} +\Lambda_3\frac{\partial {\cal G}_3}{\partial x} 
\nonumber
\\
&& \hspace{.2cm} \implies \ \ 
\log E\Big(\sum_{l} X_l(t)\Big) \frac{  e^{-f(x\eta)^2} }{\sqrt{2 \pi} \sigma_b x } 
= -\Lambda + \Lambda_3 \label{Eqn_1CPLam}
\end{eqnarray}
\vspace{-.51cm}
{\small{
\begin{eqnarray}
&&\frac{\partial\log E\Big(\sum_{l} X_l(t)\Big)P( \mathbf{B} \leq x \eta) }{\partial \eta}  =  \Lambda\frac{\partial {\cal G}}{\partial \eta} + \Lambda_1\frac{\partial {\cal G}_1}{\partial \eta} +
\Lambda_2\frac{\partial {\cal G}_2}{\partial \eta} +\Lambda_3\frac{\partial {\cal G}_3}{\partial \eta} \nonumber
\\
&& \hspace{-.2cm} \implies  
\log E\Big(\sum_{l} X_l(t)\Big) \frac{  e^{-f(x\eta)^2} }{\sqrt{2 \pi} \sigma_b \eta } 
+ P( \mathbf{B} \leq x \eta)  \frac{\partial\log E\Big(\sum_{l} X_l(t)\Big)}{\partial \eta} =  -\kappa_1\Lambda -\Lambda_1-\Lambda_2  \label{Eqn_1CPk1Lam}
\\
&&\hspace{-2cm} \frac{\partial\log E\Big(\sum_{l} X_l(t)\Big)P( \mathbf{B}\leq x \eta) }{\partial s}   =  \Lambda\frac{\partial {\cal G}}{\partial s} + \Lambda_1\frac{\partial {\cal G}_1}{\partial s} +
\Lambda_2\frac{\partial {\cal G}_2}{\partial s} +\Lambda_3\frac{\partial {\cal G}_3}{\partial s} \implies  0 = 2\Lambda s   \label{Eqn_1CPLamS}
\end{eqnarray}}}
and also $ x + \kappa_1 \eta +s^2 = \bar{B}.$  We now compute the gradient w.r.t. to $s_1, s_2, s_3$:\footnote{Note that we mainly require equations \ref{Eqn_1CPLamS} and \ref{Ls3} for this proof.} 
{\small{
\begin{eqnarray}
\frac{\partial\log E\Big(\sum_{l} X_l(t)\Big)P( \mathbf{B}\leq x \eta) }{\partial s_1}   = \Lambda\frac{\partial {\cal G}}{\partial s_1} + \Lambda_1\frac{\partial {\cal G}_1}{\partial s_1} +
\Lambda_2\frac{\partial {\cal G}_2}{\partial s_1} +\Lambda_3\frac{\partial {\cal G}_3}{\partial s_1} & \implies  0 = -2s_1 \Lambda_1 \label{Ls1}
\\
\frac{\partial\log E\Big(\sum_{l} X_l(t)\Big)P( \mathbf{B}\leq x \eta) }{\partial s_2}   = \Lambda\frac{\partial {\cal G}}{\partial s_2} + \Lambda_1\frac{\partial {\cal G}_1}{\partial s_2} +
\Lambda_2\frac{\partial {\cal G}_2}{\partial s_2} +\Lambda_3\frac{\partial {\cal G}_3}{\partial s_2}& \implies  0 = 2s_2 \Lambda_2 \label{Ls2}
\\
\frac{\partial\log E\Big(\sum_{l} X_l(t)\Big)P( \mathbf{B}\leq x \eta) }{\partial s_3}   = \Lambda\frac{\partial {\cal G}}{\partial s_3} + \Lambda_1\frac{\partial {\cal G}_1}{\partial s_3} +
\Lambda_2\frac{\partial {\cal G}_2}{\partial s_3} +\Lambda_3\frac{\partial {\cal G}_3}{\partial s_3} & \implies  0 = -2s_3 \Lambda_3 \label{Ls3}
\end{eqnarray}}}
\vspace*{-1cm}
\begin{eqnarray}
\mbox{and} \ \ \ \bar{\eta} - \eta + s_1^2 = 0, \ 1-\eta  - s_2^2 =0 , x-s_3^2 =0. \label{zzz}
\end{eqnarray}   
\\
Referring to equation \ref{Ls3}, we have either $s_3 =0$ or $\Lambda_3 = 0.$ If $s_3 = 0$, then $x =0$ (see equation \ref{zzz}); which is clearly not an optimal solution (zero objective value)  as the objective can be improve when $x>0.$
In particular, we do not need to compute $\Lambda_1, \Lambda_2, s_1, s_2, s_3$ for this proof. We only need to prove that $s = 0.$ For this,  observe that equation(\ref{Eqn_1CPLamS}) gives that either $\Lambda = 0$ or $s = 0$. However, $\Lambda \ne 0$ because $\log E\Big(\sum_{l} X_l(t)\Big) \frac{  e^{-f(x\eta)^2} }{\sqrt{2 \pi} \sigma_b x } $ is positive (recall $\Lambda_3 = 0$). Therefore, we must have $ s = 0$, which consequently brings out the tightness of budget constraint $x + \kappa_1 \eta  = \bar{B} $.
  Hence proved.
  \eop
  \\ 
 
\newpage
\vspace*{13cm}
\begin{center}
{\Large{\textbf{Part-II: Competitive  Viral Marketing Branching Processes in OSNs}}}
\end{center}
\newpage
\centerline{\Large \bf Competitive  Viral Marketing Branching Processes in OSNs}

\begin{abstract}

We study the content propagation of competing contents in Online Social Networks. We  model the propagation of competing posts/contents by an appropriate  branching process. The underlying branching process turns out to be decomposable.  Consequently, the evolution of the competing posts can be drastically different from each other. We utilize the existing theory of branching process and our newly developed results on decomposable branching process to study this problem. We obtain various performance measures such as the time evolution of the population of one of competing posts, extinction probabilities, etc. We also compare our results with the results that one would obtain without considering the timeline structure. We find that one leads to draw erroneous conclusions when the timeline structure is ignored. At last, we formulate a game theoretic framework to study the competition considering the online auctions. We numerically compute the Nash equilibria.
\\
\\
 \textbf{Keywords:} Viral marketing, Branching processes, Online social network, Game theory,  Martingales, Online auctions. 
\end{abstract}

\section{Introduction}
\label{intro}
In \textit{viral marketing,} the content providers (CPs)/advertisers create contents/posts that are appealing to the users. When a user finds a post about products/services attractive, it spreads a word about it. The post is transmitted from one user to its neighbour, which causes a chain reaction. By the extensive sharing/transmission of a post  the post spreads on a massive scale, then we say the post got viral and hence this process is called viral marketing. In Part-I (\cite{Ranbir2}) of this  work, we studied viral marketing branching process for the propagation of posts corresponding to a content provider (CP). In this paper, we will extend this study to investigate the propagation of posts  corresponding to competing content providers. 
 \\
\textbf{ \large{Online social network and timelines:}}
Online Social Networks (OSNs) store volumes of information about the users. An important feature of these OSNs is the \textit{timeline } (TL) structure of the appearance of the posts. Each post appears at a certain level based on its newness  on each user's page in an OSN, for instance, News Feed in Facebook. We call this reverse chronological appearance of the posts a `timeline' (TL). There is one TL dedicated for each user. \textit{As mentioned in  Part-I \cite{Ranbir2}, no attention is paid to the TL structure of the posts/contents appearing on a user's page in viral marketing literature.} We study the content propagation of competing CPs over OSNs, considering the inherent TL structure.

\begin{figure}[h]
 \begin{center}
\begin{tikzpicture}[scale=1]
\def \n {5}
\foreach \y/\x in {1/msg,2/post-\textbf{Q},3/post-\textbf{P},4/A[2],5/A[1]}
	{ \draw (-4,\y) rectangle (-1,\y+1) node[pos=.5] {\x}; 
	}
	 \node[above, xshift =-80] {user-1};                                                                    
\def \n {5}
\foreach \y/\x in {1/B[2],2/msg,3/B[1],4/post-\textbf{P},5/post-\textbf{Q}}
	{ \draw (3,\y) rectangle (0,\y+1) node[pos=.5] {\x}; 
	} 
\node[above, xshift =40] {user-2};                                                                    
 
 \def \n {5}
\foreach \y/\x in {1/post-\textbf{P},2/D[1],3/A[2] ,4/C[2],5/C[1]}
	{ \draw (7,\y) rectangle (4,\y+1) node[pos=.5] {\x}; 
	} 
\node[above, xshift =150] {user-3};                                                                    
\end{tikzpicture}
 \end{center}
\caption{Timeline structure} \label{Fig_TLIntro}
\end{figure}
A typical example of TL structure   (for three users) with competing content (say posts {\bf P} and {\bf Q}) is shown in  Figure \ref{Fig_TLIntro}.
The figure  shows the TL, consisting of different posts at different levels, for  three users. 
Users 1 and 2 have both the posts, while user 3 has only post-{\bf P}.
Our goal is to understand the propagation of these competing posts. 
Here a natural question  to ask is:  at given time, how many users have post-\textbf{P} or post-\textbf{Q}?;  and the next immediate one is at what level does that post reside (i.e., the position)? For instance, all the three users have the post `post-\textbf{P}' on their TLs, but at different levels. It is clear  that the posts positioned on the top of the TLs receive more attention/visibility compared to the ones at lower levels.
Further, when a user has two or more posts of  competing nature, it may pay more attention to the one at top level; for example, user 2 may pay more attention to post-{\bf Q}   while user 1
may pay relatively  more attention to post {\bf P}.  
 Note also that the arrival of new contents keeps shifting/pushing down the existing contents of  a TL. Thus, a particular content of interest may reach lower levels before the user visits\footnote{The users 'visit' OSNs at random intervals of time and in each `visit' it browses some/all the new posts.}  its TL, and the user may miss it. Technically, a user can scroll through indefinite number of posts. However, it is known that users' attention is limited to the first few levels \cite{Scroll}. As  in Part-I \cite{Ranbir2}, we consider this aspect in analysing content propagation, further,  considering the extra complications that arise because of competing content on TLs. Without these key elements, one leads to draw erroneous conclusions (see Part-I \cite{Ranbir2} for similar results in the context of propagation of a single post). 
\\
\\
\textbf{\large{ Methodology:}}
Similar to Part-I, we model the content propagation of (competing) posts as an appropriate branching process  (see Part-I \cite{Ranbir2} for more details).   
The branching processes can mimic most of the phenomenon that influences the content propagation;  one can model the effects of multiple posts being forwarded to the same friend,  and multiple forwards of the same post, etc. Further, we  require  a multitype branching process to model this propagation.  This is because   we require separate counts of TLs  with the given post (say post-${\bf P}$)  at each level  and  at any time instance. A post on a higher level in TL has better chances of being read by the user. Posts of appealing nature, e.g., containing irresistible offers, have a great chance of being in circulation, and we call it the post quality factor. Posts of similar nature appearing at lower levels on the TL have smaller chances of appreciation, etc.  To study all these factors, one needs to differentiate the TLs that have the `post'  at different levels, and this is possible only through multitype branching processes (BPs).  

Further to the above, we have more factors to consider  while studying the competing contents. One may have TLs with one particular content (e.g., post-${\bf P}$ or ${\bf Q}$),  or may have TLs with both the posts but at different levels.  
The propagation of a post (say post-${\bf P}$) is impacted by  that  of  the other post (e.g., post-${\bf Q}$).  This  impact is largely different than that considered in Part-I \cite{Ranbir2}  because of the competition between the two posts. When a new post is received by a TL with post-${\bf P}$: a)   post-${\bf P}$  can  only lose its  position in higher levels, when the new post is unrelated to post-${\bf P}$ (like in Part-I);  b) if the new post is (competing) post-${\bf Q}$,  in addition,    the post-${\bf P}$ may receive reduced attention  from  the user of TL.  In Part-I \cite{Ranbir2},  we studied the propagation of content belonging to single content provider,  we required non-decomposable (or irreducible) branching processes.  Whereas for the propagation of competing content,  we find that a decomposable branching process is needed.  The analysis of decomposable is far more complicated, and is a relatively less studied object in literature.

\subsection*{An overview of branching processes}
As our analysis uses the theory of branching processes extensively, we briefly present an overview of branching process literature.  Branching processes can be categorized on a number of factors,  for example, the discrete and continuous time branching processes (classification by time), single type and multi-type branching processes,  critical, super-critical or sub-critical, etc. Further, each category has subcategories giving rise to numerous variants of the branching processes (e.g., \cite{AthreyaBook}, \citep{Harris} etc).
 
In a multitype  Markovian  \textit{continuous time branching process} (CTBP), a particle lives for an exponentially distributed random time (e.g., \cite{AthreyaPaper}). It produces a random number of offspring of various types independent of the other particles and then dies. And this continues. The underlying generator  matrix, say $A$,  plays a vital role in carrying out the analysis of continuous time branching process (CTBP). When the $e^{At}$ matrix is   positive regular\footnote{Matrix $e^{At}$ is positive regular if each entry of it is strictly positive for some $t_0 >0.$  }, the underlying CTBP is classified as irreducible/non-decomposable. In this case, the branching processes exhibit a certain dichotomy (e.g., \cite{AthreyaBook,AthreyaPaper}):  a) either all the types survive together and grow exponentially (with time) with the same rate; b) or all the types get extinct after some time.  Additionally,  the largest eigenvalue of  $A$, say $\alpha$, determines the growth rate, extinction probability, etc.
  The CTBP is called subcritical, critical and supercritical    based on  whether $\alpha < 0, \alpha = 0$ and $\alpha > 0$ respectively. 
When $\alpha \le 0$,  the population gets extinct with probability one. Whereas when  $\alpha > 0$, the CTBP can survive on some sample paths. Further, on these sample paths, all types grow exponentially fast with the common rate $\alpha$  provided the CTBP is non-decomposable \cite{AthreyaPaper}.  

When the process is such that the particles of certain types do not produce offspring of certain other types,  we have a very different variety of branching process called as decomposable branching process. In this process, the types get partitioned into different classes, where the types across different classes may have different characteristics. These processes behave significantly different from  non-decomposable processes.  First and foremost, the dichotomy no longer holds, i.e., a particular class  (a group of types) may thrive/survive whereas another gets extinct.
Secondly, the types in different classes may have different growth rates, etc.

 It turns out that the branching processes (BPs) modeling the competing content is a  decomposable branching processes and this necessitated the extension of the theory of the corresponding  BPs.  In particular we  obtained the  time evolution of  expected performance   of a metric  that resembles  the well known    total progeny (or total population, i.e., the total number of population from the start, including the perished ones) of the BP.  This performance helps in estimating the expected `number of shares'. 

{ 
The decomposable processes are relatively less studied in comparison with the non-decomposable objects. There are strands of literature (e.g., \cite{Kesten,Haut}) that study the discrete time decomposable branching processes. But to the best of our knowledge, there is no dedicated literature on decomposable continuous time branching process (DCTBPs). 
Some analysis could probably be derived using the discrete-time results. However, there are many more questions that need to be answered directly for continuous time versions, and we consider the same.  

Authors in  \cite{Kesten} studied the multitype decomposable branching process in  discrete time framework. Their study mainly focused on investigating the growth rates of the particles belonging to different classes.  Authors in \cite{Haut} investigated the class-wise extinction probabilities, where the extinction of a class is shown to be the minimal non-negative solution of the extinction probability equation but with added constraints.   We  mainly study  the `total progeny' of different classes  in a decomposable branching processes, and  in the context of continuous time processes. The continuous time processes are useful in studying many practical problems such as cancer biology, viral marketing problem etc. In these problems, it becomes important to know the following: the growth rates of different classes, 
the growth rates of `total progeny' of different classes, the probability that  the particles of one specific class explodes, while, the others get extinct,  etc. 

\ignore{ We find that the results  required  for our analysis are missing in the current  literature related to decomposable branching process. By this way, we contribute to decomposable branching processes. In particular, we derive a performance measure similar to total progeny (which is well analyzed for irreducible branching processes),  for  decomposable branching processes in the continuous time framework. 

 We employ the well-known results of branching processes and our newly developed results on decomposable processes to study competing posts.  The decomposable branching processes are relatively less studied objects, particularly in the continuous
time framework.}
As usual practice in the theory of decomposable branching processes, we group
various types into different irreducible classes. These irreducible classes evolve according to the well-studied nondecomposable/irreducible branching processes  (when they start in their own class) and we study further evolution of them when they are inter connected according to a reducible generator matrix. }

\section{System description} 
We consider an OSN with large number of users, for example, Facebook, Twitter, etc.  Users use these networks to share pieces of information such as messages, photos, videos etc. We briefly refer to these pieces of information as  posts. The posts are stored in a reverse chronological order on inverse stacks which we refer to as timelines (TLs).  When a user visits  the OSN, it reads the posts on its timeline (TL)  and shares a post, upon finding it appealing/useful, with some of its friends. 
Due to this, the shared post appears on the top level of the timelines of those friends with whom the post was shared. This brings about a change in the appearance of contents on the timelines  of the recipients of the post. Basically, the existing contents of these TLs shift one level down. And a user can share as many posts as it wants. The number of shares of a particular post by a particular user depends upon: a) the distribution of  number of friends of the user;  b)  the level in the TL at which the post resides; and c) the extent to which the user liked the post etc. Basically,  the sharing of a post depends on how engaging the \textit{content provider} (CP) designs its post.
And extensive sharing of the post amongst the users/friends potentially makes the post viral. 

  There are many more aspects which influence the content propagation (see Part-I \cite{Ranbir2}). Users may become reluctant to read/share the contents on the lower levels of their TLs. When they see multiple posts of similar nature, they may appreciate few posts while the remaining receive reduced attention. 
 We study all these aspects and the dynamics created by the actions (e.g., like, share, etc)  of the users, which  have a major impact on the propagation of the commercial content.
 Further, in this paper,  we consider propagation of multiple posts which compete with each other; for example two competitors can spread simultaneously their advertisements  through the same social network and users response  to one of the posts  depends also upon the post of the competitor.  
 
 We consider multi-type continuous time branching process to model the propagation of competing content. We begin with the description of the relevant dynamics and that in an appropriate branching process.
 
%
%

 \subsection{Dynamics of  content propagation and branching process}
 The content propagation in a typical OSN is as follows. Let us say we are interested in the propagation of two competing posts, namely  post-\textbf{P} and post-\textbf{Q},  when the process starts with $X(0)$ number of seed TLs.
Some of these $X(0)$ TLs have only  post-\textbf{P}/post-\textbf{Q}, while some others have both the posts.  Further, the tagged (post-\textbf{P}/post-\textbf{Q}) posts can be residing anywhere on the first $N$-levels of the corresponding TLs;  
  we track these posts only till first $N$ levels of the TL.  Note that the posts of these  $X(0)$ TLs remain unread before their respective users visit their TLs.  Thus, we call these  as \textit{ number of unread TLs }(NU-TLs). If a user, among  $X(0)$, visiting its TL finds  post-\textbf{P}/post-\textbf{Q} or both  attractive, it reads the post(s) and may share the same with a random number of its friends.  And post-\textbf{P}/post-\textbf{Q} or both  would be placed on the top level(s)  of the recipient TLs (see Figure \ref{FigShift}). As shown in the figure, the recipient TL has  post-\textbf{P}   and post-\textbf{Q} on the top levels, and the existing posts shift down by two levels. Further, it is clear (in this example) that the  post-\textbf{Q}  is shared before the  post-\textbf{P}.
 \begin{figure}[h]
 \begin{center}
\begin{tikzpicture}[scale=1]
\def \n {5}
\draw[->>, thick] (-2.5,6.8)--(-2.5,6) node[pos = -0.3]  { \textcolor{black}{\text{post-\textbf{P}}}};
\foreach \y/\x in {1/A[4],2/A[3],3/A[2],4/A[1],5/post-\textbf{Q},}
	{ \draw (-4,\y) rectangle (-1,\y+1) node[pos=.5] {\x}; 
	} 
	\draw[->>, red,thick]   (-1,3.5) -- (0.5,3.5);
\end{tikzpicture}
\begin{tikzpicture}[scale=1]
\def \n {5}
\foreach \y/\x in {1/A[3],2/A[2],3/A[1],4/post-\textbf{Q},5/post-\textbf{P}}
	{ \draw (-4,\y) rectangle (-1,\y+1) node[pos=.5] {\x}; 
	} 
\end{tikzpicture}
 \end{center}
\hspace*{-3cm}
\caption{  Shifting of the contents on a TL, when  post-\textbf{P} and post-\textbf{Q} are shared with it} \label{FigShift}
\end{figure}

If some more posts are shared   with some of these recipient TLs, their contents further shift down.  For instance, in Figure \ref{FigShift},  if one more post is shared after  post-\textbf{P}, the post-\textbf{P}  would reside on the second level (and   post-\textbf{Q} on level 3) of the  TL. 
  
  { 
As argued in Part-I \cite{Ranbir2}, }the continuous time version of the branching process fits the content propagation better than the discrete counterpart. In a CTBP,  any \underline {one} of the existing particles `dies'  after exponentially distributed time while in a discrete time version all the particles of a generation `die' together;
the  users of the  TLs with  post-\textbf{P}/post-\textbf{Q}  visit their respective TLs at different instances of times.    As the underlying OSN is huge, one can say that the visit times of users are virtually independent of the each other. We assume that a user visits its TL after exponentially distributed time.
 
 When the number of copies of a CP-post grows fast (i.e., when the post is viral), the time period between two subsequent changes decreases rapidly as time progresses. This is also well captured by CTBP, which mimics the content dynamics better.
 
 The sharing process generates a random number, say $\zeta$, of new TLs holding post-\textbf{P} or post-\textbf{Q} or both. If the user does not read or share the post after visiting its TL, then $\zeta = 0$. If sharing process is independent and identical across all the users, the new TLs $\zeta$ so generated resemble IID  offspring in a CTBP and the effective NU-TLs with one or both of the posts  may appear like the particles of a CTBP. 
 When one of the users of these number of unread TL  (NU-TLs) visits its TL  and starts sharing the  post-${\bf P}$/post-\textbf{Q} or both (as before),   then the content propagation dynamics again resemble that in a CTBP. 
 
 However, the CTBP  described above does not capture some aspects related to the modeling of the post-propagation process.   Post-${\bf P}$/post-\textbf{Q} can disappear from some of the TLs, before the corresponding user's visit. For example, post-${\bf P}$ would disappear from a TL with $(N-l+1)$ or more shares (before user's visit), if initially post-${\bf P}$ were at level $l$.  Further, TLs with only one of the two posts are different from the ones that have both, etc.  
 
 Thus, we will need (continuous time branching processes) CTBPs  with multiple types of particles to model   this kind of content propagation. 
Further, with competing content, one of posts may get viral and the other may get extinct.  This kind of an effect is not seen in irreducible BPs. Thus, the BPs that model our process  cannot be irreducible, they will have to be  decomposable  {(e.g., \cite{dec1, Kesten}). }
%

\section{Modeling details}
 We consider \textit{two} competing content providers (CPs) and refer to them as CP-1 and CP-2 respectively. The competing CPs are operating in a similar kind of business, e.g., tourism industry, hotel/restaurant services, manufacturing businesses, etc. And they have  competing contents/posts. We track the posts of both the CPs till first N levels of the TLs.  The propagation of competing  content can   be  modelled by a \textit{multitype branching process} (MTBP). 
As already mentioned, we have multiple types of population, and they further can be classified into three classes, as explained below. 


\subsection{Different Types  of TLs}
\label{sec_Diff_TLs}
\subsubsection {Exclusive-types} There are two classes in this category and each class contains the users  holding the post of  one of the competing CPs only on their TLs.  Let $X_{l,0}(t)$ be the number of unread TLs having CP-1 post at level $l$ and these do  not contain CP-2 post at time $t$. And let $ \textbf{X}^1_{ex}(t) : = \{X_{1,0}(t),  \cdots, X_{N,0}(t)\} $  and  $\textbf{X}^2_{ex}(t) := \{X_{0,1}(t),   \cdots, X_{0,N}(t)\}$ denote the population vector of NU-TLs holding CP-1 and CP-2 posts respectively.  These types are exactly like those in Part-I (\cite{Ranbir2}). We describe the details with one exclusive type post (say post-${\bf P}$) and we refer to it as CP-post. We have two types of transitions namely shift and share transitions that modify the NU-TLs holding exclusively the CP-post.  Figure \ref{Fig_TL_SSCh} demonstrates the transitions.
 \begin{figure}[ht]
\begin{tikzpicture}[scale =1.4]
\node[ rounded corners] at (1.1,2.8){Timeline of a user};
\draw[ xscale = 0.5, yscale =0.5, thick, fill= blue!10,rounded corners] (0,0) rectangle (5,5)  node[pos=.5]  { CP-post (at level-3)};
\foreach \x in {1,...,5}  \draw[xscale = 0.5, yscale =0.5] (0,\x) -- (5,\x ) (-1,5-\x +.5 ) node {Level-\x};
\draw[->, very thick] (2.5,1.25)  to[out = 80,in = 140] (6,4.5);
\draw  (3,3) node[ above,sloped] { \rotatebox{55}  {Shift transition w.p. $\theta$}}; 
\draw[ xscale = 0.5, yscale =0.5, thick, fill= blue!10,rounded corners] (10,4) rectangle (15,9)  (12.5,5.4) node  { \footnotesize{CP-post now at level-4}};
\foreach \x in {1,...,5}  \draw[xscale = 0.5, yscale =0.5] (10,\x+3) -- (15,\x +3) (9,9-\x +.5 ) node {Level-\x};
\draw[->, very thick] (2.5,1.25) -- (6,-.5);
\draw    (4,0) node{ \rotatebox{-20} {Share transition w.p. $1-\theta$}};
\draw[thick, fill= blue!10,rounded corners] (6,-2) rectangle ++(2.5,2.5) node[align=center,pos=.5] {CP-post shared with  \\ $\zeta$ no. of friends \\ \\ \# CP-post++};
\end{tikzpicture}
\caption{Propagation of CP-post: transitions}\label{Fig_TL_SSCh}
\end{figure}

\noindent
{\bf Share Transitions:} 
In the share transition, a user first reads the CP-post and based on the interest generated, it shares CP-post with a random number of friends. 
When a user visits its TL, it reads some/all posts located on its TL and shares them with some of its friends. 
We illustrate the sharing transition in Figure \ref{Fig_TL_ShareCh}.
 \begin{figure}[h]
\begin{tikzpicture}[scale =1.1]
\node[ rounded corners] at (.7,2.8){ TL};
\draw[ xscale = 0.5, yscale =0.5, thick, fill= blue!10,rounded corners] (0,0) rectangle (3,5)  node[pos=.5]  { -  - } ;
\draw (0.7,1.7) node {{\normalsize{ CP-post}}  };
\draw[ ->, very thick ] (1.5,1.25) to ++(1,0);
\foreach \x in {1,...,5}  \draw[xscale = 0.5, yscale =0.5] (0,\x) -- (3,\x ) (-1,5-\x +.5 ) node {\tiny{Level}-\x};
\draw[thick, fill= blue!10,rounded corners] (2.5,0) rectangle ++(2,2.5) node[align=center,pos=.5] { \footnotesize{Reads CP-post} \\  w.p.  $r_2$ \\ (at level-2)};
\draw[ ->,very thick ] (4.5,1.25) to ++(1,0);
\draw[thick, fill= blue!10,rounded corners] (5.5,0) rectangle ++(2,2.5) node[align=center,pos=.5] {Interest \\ generated \\  w.p. $\eta$};
\draw[ ->, very thick ] (7.5,1.25) to ++(1,0);
\draw[thick, fill= blue!10,rounded corners] (8.5,0) rectangle ++(1.6,2.5) node[align=center,pos=.5] { \normalsize{CP-post } \\ shared  \\ \normalsize{with $\zeta$} \\  \# friends};
\draw[ ->, very thick ] (10.1,1.25) to ++(.9,0);
\draw (10.5,1.4) node[align = right, below] {\small{$\rho_3$}};
\draw[thick, fill= blue!10,rounded corners] (11,0) rectangle ++(2,2.5)  node[pos=.5]  { {\footnotesize{ \#CP-post-3++}}  };
\draw (12,2.2) node {{\footnotesize{ \#CP-post-1++}}  } (12,1.7) node { } (12,.7) node {   } (12,.2) node { {\footnotesize{ \#CP-post-5++}}  };
\foreach \x in {1,...,5}  \draw[xscale = 0.5, yscale =0.5] (22.6,\x) -- (25.6,\x );
\path (9.95,1.1) node(x) {}  (11.1,2.5) node(a) {} (11.1,1.7) node(b) {}  (11.1,.75) node(d) {}  (11.1,0) node(e) {};
\draw [ ->, very thick ]  (x) -- (a)  node[near end,below] {\small{$\rho_1$}}; 
\draw [ ->, very thick ]  (x)   -- (e) node[near end,above] {\small{$\rho_5$}}; 
\end{tikzpicture}
\caption{Share transition giving birth to exclusive types  }\label{Fig_TL_ShareCh}
\end{figure} 
The user   reads (and shares) the posts residing on different levels, with varying levels of interest based on many factors (as in Part-I \cite{Ranbir2}). Firstly, it reads posts on higher levels with higher probabilities than those on lower levels; the interest can also depend upon the influence of the content provider; and it can further depend upon the quality of the post etc.    We assign probability $r_i$ for reading the post at level-$i$, and note that $r_1 \ge r_2 \cdots \ge r_N$. If the user finds the post interesting, which is determined by the post quality factor $\eta$, it shares the post with random number of its friends. And if the user shares more posts along with the tagged post (e.g., post-${\bf P}$), the level of the  post  (in  the recipient TLs) changes accordingly. 
We consider this aspect into our model by defining $\rho_i$, where $\rho_i$ is the probability of sharing $i$ number of posts (see Figure \ref{Fig_TL_ShareCh}). Furthermore, a user can respond more actively to the post of a more influential CP. Let $w_j$  ($\ge 1$) be the influence factor of \CPj  with  $j=1$ or $2$, and we assume that  the post quality factor  of 
  \CPj is given by $w_j \eta_j$.
    Thus, the CP with high influence factor can obtain good results even with a  lower post quality.  To {\it simplify the notation, we use $\eta_j$ to represent $w_j \eta_j$} and with this  notation,  $\eta_j \in [0, 1/w_j]$.

\noindent
{\bf Shift Transition:}  In the shift transition, when the TL of user with CP-post is written by other users,   the position of CP-post shifts down (see Figure \ref{FigShift} and \ref{Fig_TL_ShareCh}).
 
\noindent{\bf  CP-post propagation dynamics:}  
Let  $\mathcal{G}_1$ represent the subset of users with CP-post\footnote{Similarly  $\mathcal{G}_2$  contains subset of  users with  post of CP-2.} at some level, while  $\mathcal{G}$ contains the other users of social network without any post of our interest (i.e., post-\textbf{P} and post-\textbf{Q}). 
We assume the OSN  (and hence $\mathcal{G}$) has infinitely many users and note  
 $\mathcal{G}_1$ at time $t$ has, 
 \begin{eqnarray}
X_{ex}^1(t) &:=& \sum_{l\le N} X_{l,0} (t), \label{Eqn_tot_num_ex}
 \end{eqnarray}
 number of users.  Group ${\cal G}$ has an infinite number of users/agents, and this remains the same irrespective of the size of ${\cal G}_1$ (and  ${\cal G}_2$), which is finite at any finite time. 
Thus, the transitions between ${\cal G}$ and  ${\cal G}_1$ are more significant, and one can neglect the transitions within  ${\cal G}_1$.  It is obvious that we are not interested in transitions within  ${\cal G}$ (users without CP-posts). We thus model the action of these groups in  the following consolidated manner: 
\begin{itemize}
\item In the share transition, any user from $ \mathcal{G}_1 $ wakes up after $exp(\nu)$ time (exponentially distributed with parameter $\nu$)  to visit its TL and writes to a random (IID) number of users of $\mathcal{G}$  (refer to Figure  \ref{Fig_TL_ShareCh}).
\item In the shift transition,  The TL of  any user of  ${\cal G}_1$ is written by  one of the users of ${\cal G}$, and the time intervals between two successive writes are exponentially distributed with  
parameter~$\lambda$.
\end{itemize} 

The  state\footnote{to be more precise,  the components of the entire system state, corresponding to the post of CP-1.} of the network, ${\bf X}_{ex}^1 (t)$, changes when the first of the above-mentioned events occurs. At time $t$, we have $X_{ex}^1(t)$ (see equation (\ref{Eqn_tot_num_ex})) number of users in $\mathcal{G}_1$ and thus (first) one of them wakes up according to exponential distribution with
parameter $X_{ex}^1(t)\nu$. 
Similarly, the first TL/user of the group ${\cal G}_1$ is written with a post after exponential time with parameter~$X_{ex}^1(t)\lambda$. 
Thus, the state ${\bf X}_{ex}^1 (t)$, changes    after exponential time with
parameter $X_{ex}^1(t)\lambda + X_{ex}^1(t)\nu$.
Thus, the rate of transitions at any time is proportional to $X_{ex}^1(t)$, the number of NU-TLs at that time, and hence, the rate of transitions increase sharply as time progresses,  
when the post gets viral. Considering all the modeling aspects,  the   IID offspring generated by one $(l,0)$-type user are summarized as below (w.p. means with probability):
\begin{eqnarray}
\xi_{l,0} = 
\left\{
\begin{array}{llll}
{\bf e}_{l+1}\mathbb{1}_{l < N} & \text{ w.p.}\ \ \theta  := \frac{\lambda}{\lambda + \nu} \mbox{ and }  \\ 
\zeta {\bf  e}_i & \text{ w.p.} \ \  (1-\theta ) r_l\rho_i \hspace{0.4cm} \forall i \leq   N \\
0 & \text{ w.p. } \ \ (1-\theta)(1-r_l).
\end{array}
\right.   
\label{Eqn_offspring_ex}
\end{eqnarray}
where $\textbf{e}_l$ represents standard unit vector of size $N$ with one in the $l$-th position, $\mathbb{1}_A$ represents the indicator, $\zeta$ is the random number of friends to whom the post is shared   and $r_l$ is the probability the user reads/views a post on level $l$. Recall that users (offspring) of exclusive type $(0,i)$ are produced with probability $\rho_i$ during the share transitions.

  From equation (\ref{Eqn_offspring_ex}) the offspring distribution  is  identical at all 
  time  instances $t$, $\zeta$ can be assumed  independent across users, and  hence  $\xi_l$ are IID offspring from any type $l$ user.
Further, all the transitions occur after memoryless exponential times, and hence ${\bf X}_{ex}^1 (t)$ by itself is  an MTBP with  $N$- types  (e.g. \cite{AthreyaPaper}), when one starts only this exclusive type TLs. 
\\
\\
\noindent{\bf PGFs and post quality factor:}
Let $f_F( s, \beta)$ be the probability generating function (PGF) of 
 the number of friends, $\Friends$, of a  typical user, parametrized by   $\beta$. 
For example,$f_F(s, \beta) = \exp (\beta (s-1) ) $ stands for Poisson distributed $\Friends $,   $f_F(s,  \beta) = (1-\beta)/(1- \beta s) \ $ stands
for geometric  $\Friends $.  Let $m = f'_F (1, \beta)$ represent the corresponding mean.
A user shares the post with some/all of its friends ($\zeta$ of equation (\ref{Eqn_offspring_ex})) based on how engaging the post is. As mentioned before,  the post quality factor $\eta $ quantifies the extent of the CP-post engagement on a (continuous) scale of 0 to $w_1$,  where $\eta = 0$ means the worst and $\eta =w_1$ is the best quality. {\it We assume that the mean  of  the number of shares  is proportional to this quality factor.}  In other words,   $m (\eta) = m \eta  $ represents the post quality dependent mean of the random shares.  Let $f(s, \eta, \beta)$ represent the   PGF of $\zeta$.  For example, for Poisson  friends, the PGF and the expected value of $\zeta$  are given respectively by:
$$f (s, \eta, \beta) = f_F(s, \eta \beta) = \exp (\beta \eta (s-1) ) \mbox{ for any $s$ and } m (\eta) = \eta \beta. $$

For Geometric friends, one may assume the post quality dependent parameter 
$$ \beta_\eta =  {(1-\beta)}/{(1-\beta+\beta\eta)}, \mbox{  which ensures  } m(\eta) = \eta \beta .$$ And then the PGF of $\zeta$ is given by
$
f(s,\eta,\beta) = f_F(s, \beta_\eta) = {(1-\beta_\eta)}/{(1-\beta_\eta s)}.$
One can derive such PGFs for other distributions of $\Friends.$  \textit{Interestingly enough, we find that most  of the analysis does not depend upon the distribution of \  $\Friends$
but only on its expected value.}

 Let $\textbf{s} :=  (s_1, \cdots,s_N )$ and  $\bpgf (\textbf{s}, \eta ):= \sum_{i = 1}^N f(s_i,  \eta, \beta) \rho_i.$
The  post quality factor dependent  PGF,  of the offspring distribution of the overall branching process,  is given by (see equation (\ref{Eqn_offspring_ex})):
\vspace*{-0.2cm}
 \begin{eqnarray}
\label{Eqn_hl_ex}
\boxed{ h_{l,0}(\textbf{s})  =  \theta\left(s_{l+1}\mathbb{1}_{l<N} + \mathbb{1}_{l=N} \right) + (1-\theta)r_l \bpgf (\textbf{s}, \eta ) +(1-\theta)(1-r_l)}.
\end{eqnarray}

{\bf Generator matrix}
The key ingredient required for analysis of any  MTBP  is its generator matrix. We begin with the generator for  MTBP that represents the evolution of unread TLs with CP-post.
We refer to this process briefly as \underline{TL-CTBP},   timeline continuous time branching process. The generator matrix, $A $, is given by $A = (a_{lk})_{N\times N}$, where $ a_{lk}= a_l\left( {\partial h_l({\bf s})}/{\partial s_k} \Big | {_{ {\bf s} = {\bf 1} }} - \mathbb{1}_{\{l=k\}} \right)$  and $a_l$ represents the transition rate of a  type-$l$ particle (see \cite{AthreyaPaper} for details). 
For our case, from previous discussions   $a_l = \lambda + \nu$ for all $l$. 
Further, using equation (\ref{Eqn_hl_ex}),   the  matrix $A$  for our single CP case is given by (with $c := (1-\theta) m\eta$, $ c_l = c\rho_l$)

\begin{eqnarray}
A_{ex}^1 = (\lambda + \nu)
\label{Eqn_genmatrix_ex}
 \left[ \begin{array}{ccccccc}
 c_1 r_1 -1 & c_2r_1 + \theta & \cdots  & c_{N-1}r_1  & c_Nr_1 \\
 c_1 r_2   & c_2r_2 -1 & \cdots & c_{N-1}r_2  & c_Nr_2 \\
  & \vdots    \\
 c_1 r_{N-1}  & c_2r_{N-1}   & \cdots &  c_{N-1}r_{N-1} -1 &  c_{N} r_{N-1} + \theta \\
 c_1 r_N  &  c_2r_N  &  \cdots &  c_{N-1}r_N  & c_{N} r_N -1 \\
 \end{array} 
 \right] .
\end{eqnarray} 

The exclusive types corresponding to CP-2 can be defined in exactly a similar way.

 \subsubsection{Mixed-types:} These are the TLs having the posts of both the CPs (i.e., both  post-${\bf P}$ and post-${\bf Q}$), i.e., the TLs in $\mathcal{G}_1 \cap \mathcal{G}_2$.  Denote by $X_{l,k} (t)$  the number of users with  post-${\bf P}$ on the $l$-th level and  post-${\bf Q}$ on the $k$-th level of their TLs at time $t$. We classify these  TLs as $(l,k)$ type TLs.

We consider the analysis with initial TLs having the post-${\bf P}$ and post-${\bf Q}$ on the top levels, i.e., we begin with either (1,2) or (2,1) type TLs.  It is not difficult to start with other types of TLs, but the expressions become complicated, and we would like to explain the results in a simplified manner.  Now, with a shift transition, a (1,2) type (a (2,1) type) gets converted to a $(2,3)$  type (a (3,2) type respectively), which further gets converted to $(3,4)$ (to (4,3)) type with another shift, and so on. Thus, we have $2(N-1)$ mixed-type TLs, which  at time $t$ are  given by,  $ \textbf{X}_{mx} (t)  = \left( \textbf{X}_{mx1} (t),  \,  \textbf{X}_{mx2}(t) \right)$  with
 \begin{eqnarray*}
\textbf{X}_{mx1} (t) & := & \{ X_{1,2}(t),  X_{2,3}(t), \cdots, X_{N-1,N} (t) \} ,\\
\textbf{X}_{mx2}(t) & := & \{ X_{2,1} (t), X_{3,2}(t), \cdots,  X_{N,N-1} (t) \}.
\end{eqnarray*}
 
And the group  $\mathcal{G}$, as before, has all the other TLs without the post of either CP.\\
\\
\textbf{Transitions:} Recall that a \textit{shift transition} occurs when a user of $\mathcal{G}$ writes to a user/TL of $\mathcal{G}_1/\mathcal{G}_2.$ In this event, the exclusive-type TLs are changed in as described in (\ref{Eqn_offspring_ex}).  
While for the mixed-types, the position of each post slides down by one level. For example an $(l,k)$ type TL  (with  $l = k+1$ or $k-1$) gets converted  to $(l+1,k+1)$ type
when $l, k < N$; and  $(N-1, N)$ and  $(N, N-1)$ type TLs  get converted to exclusive-types   $(N,0)$   and $(0, N)$  respectively.  
  
In the \textit{share transition},  the exclusive-type TLs propagate as in the case of single CP. Whereas a mixed-type TL, say $(l,l+1)$,  undergoes the  following changes when subjected to the share transition
 \begin{enumerate}[label=\roman{*}), ref=(\roman{*})]
 \item The user first views the post-${\bf P}$ with probability (w.p.) $r_l$ and shares the same with some of its friends (as in single CP case).
 \item  The post-${\bf Q}$ is below post-${\bf P}$ and recall that the posts are of similar nature. {\it The interest of the users to read the second post of similar nature  would be lesser.} We assume that  the    user views the second post w.p.  $\delta$.
 \item When the user views/reads both the posts, it can share  post-${\bf P}$ alone with some of its friends, post-${\bf Q}$ alone with some others, and both the posts with some more.  Otherwise, only post-${\bf P}$ is shared. And the TLs of exclusive-types are produced when a user shares post-${\bf P}$ or post-${\bf Q}$ only.  While mixed-type TLs are produced when it shares both the posts.  
\item  When only one CP's post is shared,  e.g.  post-${\bf P}$, it can produce type $(i,0)$ w.p. ${\bar \rho}_i, \ i = 1,2\cdots, N-1$ and that $\sum_{i=1}^{N-1}\bar{\rho_i} =1.$   It can not produce $(N, 0)$ type as the user has already discarded one post, that of  CP2.  Recall that a TL of type $i$ is produced when $(i-1)$ more posts are shared with it after the CP's post.  When both the posts are shared with the same friend, the mixed-type   $(i+1, i)$ and $(i, i+1)$  (with  $i < N$) are  produced  w.p. $p{\bar \rho}_i$ and $(1-p){\bar \rho}_i$ respectively.  With high probability, post-${\bf P}$ is shared first followed by sharing of post-${\bf Q}$,  as we started with $(l, l+1)$ TL.  Hence, the order of the posts in the recipient TLs would be reversed with high probability, and  $p$ would, in general, be larger than $(1-p)$. 
 \end{enumerate} 
We have similar  transitions with   $(l+1,l)$ type TLs also.  
  
\subsubsection{ PGF and the overall generator matrix}  

The PGF for the two CPs case can be obtained using the above modeling details as before. The random number of friends with whom both the posts are shared, is now parametrized by $\eta_1 \eta_2$. Whereas the random number of friends with whom exclusive post-${\bf P}$ or post-${\bf Q}$ is shared, is parametrized by $\eta_1(1-\eta_2)$ or $\eta_2 (1-\eta_1)$ respectively.
For example, if the number of friends, $\Friends$, is Poisson with parameter
$\beta$, then random (sampled) number of friends with whom both the posts are shared
is  Poisson with parameter  $\beta \eta_1 \eta_2$.

Denote by  $h_{l, l+1} ({\bf s})$  the PGF   for   $(l,l+1)$ type with the following notations:  $${\bf s} := \{ \textbf{s}^1_{ex}, \textbf{s}^2_{ex},\textbf{s}_{mx1}, \textbf{s}_{mx2}\} 
, \ \ \textbf{s}_{mx1} = \{s_{l,l+1}\}, \ \ \textbf{s}_{mx2} = \{s_{l+1,l}\}, \ \mbox{and \ }  \bpgftwo \left(\textbf{s}, \eta \right) := \sum_{i =1}^{N-1}f\left( s_{i}, \eta,\beta \right)\bar{\rho_i}.$$
  We obtain $h_{l, l+1} ({\bf s})$ by conditioning on the events of the first transition,  
{\small{
\begin{eqnarray*}
h_{l, l+1} ({\bf s})& = & \theta  \left ( s_{l+1,l+2} \mathbb{1}_{l<N-1} + s_{N,0}\mathbb{1}_{l=N-1} \right)  +  (1-\theta)(1-r_l)  + (1-\theta)  r_l(1-\delta)  \bpgftwo \left( \textbf{s}^1_{ex}, \eta_1 \right)   
\\ 
& & \hspace{-1.7cm} + \ (1-\theta)  r_l 
   \delta \Bigg( \Big ( (1-p)\bpgftwo \left(\textbf{s}_{mx1}, \eta_1\eta_2 \right) + \  p\bpgftwo \left(\textbf{s}_{mx2}, \eta_1 \eta_2 \right) \hspace{-0.5mm} \Big )
 \bpgftwo \left( \textbf{s}^1_{ex}, \eta_1(1-\eta_2) \right)\bpgftwo \left( \textbf{s}^2_{ex}, \eta_2 (1-\eta_1) \right  )  \Bigg )  \\
 \mbox{ and } \\
h_{l+1, l} ({\bf s})& = & \theta  \left ( s_{l+2,l+1} \mathbb{1}_{l<N-1} + s_{0,N}\mathbb{1}_{l=N-1} \right)  +  (1-\theta)(1-r_l)  +(1-\theta)  r_l   (1-\delta)  \bpgftwo \left( \textbf{s}^2_{ex}, \eta_2 \right)  
\nonumber
\\ 
& & \hspace{-1.4cm} + (1-\theta)  r_l  
   \delta \Bigg( \Big ( p\bpgftwo \left(\textbf{s}_{mx1}, \eta_1\eta_2 \right) + \  (1-p)\bpgftwo \left(\textbf{s}_{mx2}, \eta_1 \eta_2 \right) \nonumber \hspace{-0.5mm} \Big )
    \bpgftwo \left( \textbf{s}^1_{ex}, \eta_1(1-\eta_2) \right)\bpgftwo \left( \textbf{s}^2_{ex}, \eta_2 (1-\eta_1) \right  )  \Bigg ). \nonumber
\end{eqnarray*}}}
\\
\noindent And the PGF for exclusive-types is  as in the single CP case, e.g.,   $ h_{l, 0} ({\bf s}) = h_{l,0} ({\bf s}^1_{ex}).$ 
 The generator matrix  $\mathbb{A}$ has the following structure:
\begin{eqnarray}
 \mathbb{A} = \left[\begin{array}{ccc}
 A_{mx} & A^1_{mx,ex} & A^2_{mx,ex} \\
 \textbf{0} & A^1_{ex} & \textbf{0} \\
  \textbf{0} &  \textbf{0}& A^2_{ex}
 \end{array} 
 \right],
\label{Eqn_Matrix_Mix_Generator}
\end{eqnarray}
where: a) matrices  $A^j_{ex}$ for $j = 1, 2$  are  as in (\ref{Eqn_genmatrix_ex});
  b) the matrix $A_{mx}$  corresponds to transitions within the mixed-types and is
given by the following when types are arranged in the following order  $(1,2),(2,1),(2,3),(3,2),\cdots, (N-1,N),(N, N-1)$,
 \begin{eqnarray*}
A_{mx} & = & (\lambda + \nu)\left[  \begin{array}{cccccc}
 z'_1r_1 -1 &  z_1 r_1 &  \theta +z'_2 r_1 &  \dots  & z'_{N-1}r_1 & z_{N-1}r_1 \\
  z_1 r_1 &   z'_1r_1 -1  & z_2 r_1  &   \dots  & z_{N-1}r_1  & z'_{N-1}r_1 \\
  z'_1r_2 &  z_1 r_2 &  z'_2r_2- 1&   \dots  & z'_{N-1}r_2  & z_{N-1} r_2 \\
  z_1 r_2  &  z'_1r_2 & z_2 r_2  & \dots  & z_{N-1}r_2  & z'_{N-1}r_2 
  \\
  \vdots & \vdots & \vdots   & \ddots     & \vdots   \\
  z'_1 r_{N-2}  &  z_1 r_{N-2}  & z'_2r_{N-2} &   \dots  & \theta + z'_{N-1} r_{N-2} & z_{N-1} r_{N-2}
 \\
  z_1 r_{N-2}  &  z'_1r_{N-2}  & z_2 r_{N-2} & \dots  & z_{N-1}r_{N-2}  & \hspace{-2mm}  \theta + z'_{N-1}r_{N-2}   \\
z'_1 r_{N-1}  &  z_1 r_{N-1}  & z'_1 &   \cdots  & z'_{N-1} r_{N-1}-1 & z_{N-1} r_{N-1} \\
  z_1 r_{N-1}  &  z'_1r_{N-1} & z_2 r_{N-1}  &  \dots  & z_{N-1}r_{N-1}  & z'_{N-1} r_{N-1} -1
 \end{array}   
 \right],
\end{eqnarray*}
with 
$c_{mx} := \delta(1-\theta)\eta_1\eta_2 m $, $z'_i := (1-p)c_{mx}\bar{\rho_i}$ and $ z_i := p c_{mx}\bar{\rho_i}$  for all $i$;
  and  c) the matrix $A^j_{mx,ex}$ for $j=1,2$ represents the transitions from mixed-types to exclusive-types (exclusive CP types) and 
   \begin{eqnarray*}
A^1_{mx,ex} & = & \bordermatrix{
~ & 1,0 & 2,0 & \cdots & N-1,0 & N,0 \cr
1,2& c_{mx,1} r_1 \bar{\rho}_1 &  c_{mx,1} r_1 \bar{\rho}_2 & \cdots & c_{mx,1} r_1 
\bar{\rho}_{N-1}  & 0  \cr
2,1&c'_{mx,1} r_1 \bar{\rho}_1 &  c'_{mx,1} r_1 \bar{\rho}_2 & \cdots & c'_{mx,1} r_1 \bar{\rho}_{N-1}  & 0  \cr
2,3 & c_{mx,1} r_2 \bar{\rho}_1 &  c_{mx,1} r_2 \bar{\rho}_2 & \cdots & c_{mx,1} r_2 \bar{\rho}_{N-1}  & 0  \cr
3,2 &  c'_{mx,1} r_2 \bar{\rho}_1 &  c'_{mx,1} r_2 \bar{\rho}_2 & \cdots & c'_{mx,1} r_2 \bar{\rho}_{N-1}  & 0  \cr
 \vdots  &\vdots & \vdots & \ddots   & \vdots     & \vdots     \cr
 N-2, N-1 &  c_{mx,1} r_{N-2} \bar{\rho}_1 &  c_{mx,1} r_{N-2}  \bar{\rho}_2 & \cdots & c_{mx,1} r_{N-2}  \bar{\rho}_{N-1}  & 0  \cr
 N-1, N-2 & c'_{mx,1} r_{N-2}  \bar{\rho}_1 &  c'_{mx,1} r_{N-2}  \bar{\rho}_2 & \cdots & c'_{mx,1} r_{N-1}  \bar{\rho}_{N-1}  & 0  \cr
 N-1, N &  c_{mx,1} r_{N-1} \bar{\rho}_1 &  c_{mx,1} r_{N-1}  \bar{\rho}_2 & \cdots & c_{mx,1} r_{N-1}  \bar{\rho}_{N-1}  & \theta\cr
 N, N-1 & c'_{mx,1} r_{N-1}  \bar{\rho}_1 &  c'_{mx,1} r_{N-1}  \bar{\rho}_2 & \cdots & c'_{mx,1} r_{N-1}  \bar{\rho}_{N-1}  & 0 \cr}  (\lambda + \nu)  
\end{eqnarray*}
where
 \begin{eqnarray*}
c_{mx,j} &:=& (1-\theta)m\eta_j \left[1-\delta  + \delta (1-\eta_{-j} )\right] =  (1-\theta)m\eta_j(1-\delta \eta_{-j})
\\
c'_{mx,j} &:=& (1-\theta)m \delta \eta_j (1-\eta_{-j} ) \\
 \mbox{ and } \ \ -j & :=& 1 {\mathbb 1}_{\{j = 2\}} + 2  {\mathbb 1}_{\{j = 1\}}  \mbox{ is the usual game theoretic notation.}
 \end{eqnarray*}  
 Because of  ${\bf 0}$ sub-matrices of  (\ref{Eqn_Matrix_Mix_Generator}), the matrix  $\mathbb {A} $ is not positive regular, and hence, the underlying MTBP is decomposable (e.g., \cite{Haut,Kesten}).

\subsection{Analysis of the mixed-type population}    
From the structure of generator matrix $\mathbb {A} $   given by (\ref{Eqn_Matrix_Mix_Generator}), it is clear that the  subgroup of types corresponding to mixed populations, 
$ \big \{ (l, k) : l \ge 1, k \ge 1  \mbox{ and } l = k+1  \mbox{ or } k = l+1 \big \}, $
  survive on their own.
A mixed-type can be produced only by another mixed-type. Note that  the mixed-types can produce exclusive CP types, but not the other way round (see the matrix in (\ref{Eqn_Matrix_Mix_Generator})).  Thus, the extinction/virality analysis of the mixed population can be  obtained independently.  
   To begin with, we have the following result.  
\begin{theorem}
 \label{Thm_mixed}i) If \ \ $0 < \theta, p < 1$,  matrix $e^{A_{mx }t}$ for any $t > 0$, is positive regular.  
   \\
ii) Let $\alpha_{mx} $ be the largest eigenvalue of the generator matrix $A_{mx}$. Then  
  $$\alpha_{mx} \in \big ( (c_{mx}{\bf r}.{\bm {\bar \rho}}  -1),  \ (c_{mx}{\bf r}.{\bm {\bar \rho}} -1 + \theta) \big)(\lambda+\nu)$$
   where the reading probability vector
   ${\bf r}$ is redefined as $ (r_1 \cdots r_{N-1})$. 
When $r_l = d_1d_2^l \forall \ l$:

 $$\alpha_{mx}  \to(c_{mx}{\bf r}.{\bm {\bar \rho}} -1 + \theta d_2)(\lambda+\nu), \  c_{mx} := \delta(1-\theta)\eta_1\eta_2 m \mbox{ as } N \to \infty.$$
iii)  Further, the left eigenvector ${\bf u}_{mx} = (u_{mx,1}, \cdots, u_{mx, {2N-2}})$ corresponding to $\alpha_{mx}$ satisfies   for any $2\le l\leq N-1$:

 \begin{eqnarray*}
  u_{mx,2l-1} =  \sum_{i = 0}^{l-1} \frac{\bar{\rho}_{l-i}}{\bar{\rho}_1} \left(\frac{\theta}{\sigma_{mx}}\right)^{i} u_{mx,1} ;  u_{mx,2l}  =  \sum_{i = 0}^{l-1} \frac{\bar{\rho}_{l-i}}{\bar{\rho}_1} \left(\frac{\theta}{\sigma_{mx}}\right)^{i} u_{mx,2}; \ \ \sigma_{mx} = \frac{\alpha_{mx} }{\lambda + \nu} + 1.
 \end{eqnarray*}
 And the right eigenvector ${\bf v}_{mx} = (v_{mx,1}, \cdots, u_{mx, {2N-2}})$ corresponding to $\alpha_{mx}$ satisfies   for any $1 \le l\leq N-2$:
\begin{eqnarray*}
 v_{mx,2l-1} = \sum_{i = 0}^{N-1-l} \frac{r_{l+i}}{r_{N-1}} \left(\frac{\theta}{\sigma_{mx}}\right)^{i} v_{mx,2N-3}; \ \ v_{mx,2l-2} = \sum_{i = 0}^{N-1-l} \frac{r_{l+i}}{r_{N-1}} \left(\frac{\theta}{\sigma_{mx}}\right)^{i} v_{mx,2N-2}.
\end{eqnarray*}
 
iv) The process  $\{  \textbf{v}_{mx}. \textbf{X}_{mx}(t)e^{-\alpha_{mx} t} \}$  
       is a  non-negative martingale, where $\textbf{v}_{mx}$ is right eigenvector and
\begin{eqnarray*}
\lim_{t\to\infty}\textbf{X}_{mx}(t, \omega)e^{-\alpha_{mx} t} = W_{mx}(\omega)\textbf{u}_{mx}  \text{ {\normalsize for   almost all}   } \omega. 
\end{eqnarray*}
\textbf{Proof:} The proof is given in Appendix A. \eop
\end{theorem} 
From part (ii) of the above Theorem, the mixed TLs get viral when $c_{mx}{\bf r}.{\bm {\bar \rho}} > 1$, and the rate of explosion on viral paths equals $\alpha_{mx}$. When  $\alpha_{mx}  < 0$, the mixed population gets extinct surely, i.e.,
 $$
P \Big ( {\bf X}_{mx} (t) = 0 \mbox{ for some } t > 0 \Big | {\bf X}_{mx} (0) = {\bf e}_{l,k} \Big ) = 1 \mbox{ for all $(l,k)$. }
$$
However, before the extinction, the mixed-type TLs generate exclusive-type TLs which then evolve on their own. Further, these exclusive-type TLs can  get viral if  $\alpha_j >0$ ($j=1$ or 2).  And this is possible because $\alpha_{mx}$ is less than $\alpha_j$ for $j =1,2$. For example, consider the case with
 $\rho_l  = {\bar \rho}_l = \mathbb{1}_{\{ l = 1\}}$ and as $N \to \infty $,  $\alpha_j \to (1-\theta) m \eta_j -1+ \theta d_2$ while 
 $\alpha_{mx} \to  c_{mx}-1+ \theta d_2 = (1-\theta) \delta m \eta_1 \eta_2 -1+ \theta d_2$; and clearly $\alpha_j > \alpha_{mx}.$

We will now deviate to derive some results in a special type of decomposable branching process; these results would be used  later for analyzing the propagation characteristics of   competing posts.

\section{Type-changing   decomposable branching process}
\label{Sec_Decompose_BP}
In our example,   `mx' (mixed)  class  particles produce the offspring of all the classes whereas an 'ex' (exclusive) class particle produces offspring of its class only. This allows us to split the generator matrix into two sub-matrices as below, which would facilitate    independent/separate  analysis of each 'ex' class:
\begin{eqnarray}
\left[
\begin{array}{cc}
{\cal A}_{mx} & A^1_{mx,ex} \\
\textbf{0} & A^1_{ex}
\end{array}
\right] \hspace*{0.4cm} \mbox{and} \hspace*{0.4cm} \left[
\begin{array}{cc}
{\cal A}_{mx} & A^2_{mx,ex} \\
\textbf{0} & A^2_{ex}
\end{array}
\right].
\end{eqnarray}
In particular, we are interested in deriving the time evolution of  `expected net progeny'  (a measure  like  total progeny, which would be defined soon), which would represent the expected number of shares in our social network context.   

 
\label{sec_type_change}
The above processes are slightly different from the usual decomposable  branching process;    the difference lies in the  events of the transition/reproduction epoch.  Any parent at the transition epoch either produces a random number of offspring (as is usually considered in branching processes)  \textit{or} its type gets changed;  one of the two events takes place. And then it dies.  We consider a decomposable branching process consists of two irreducible classes, namely mixed (${\cal M}_x$) class and exclusive (${\cal E}_x$) class. Particles of ${\cal M}_x$ class produce particles of ${\cal M}_x$ class as well as that of ${\cal E}_x$ class. While particles of ${\cal E}_x$ produces particles of ${\cal E}_x$ class only. In one of the two events at a transition epoch, a particle of $l \in {\cal M}_x$ type wakes up after exponentially distributed time with parameter $\nu$, i.e., $exp(\nu)$ and produces a random number of offspring. 
 Whereas in the other event, the type of the particle gets changed after $exp(\lambda)$ time. 
 It is easy to see that probability of the former event is $1-\theta$ and that of the latter  is $\theta$ where $\theta = \lambda / (\lambda+ \nu).$ We refer to this process briefly as \textit{type-changing  decomposable branching process} (TC-DBP).

Let $m_{l,k}$ represent the expected number of offspring of type $k$ produced by a parent of type $l$, where $l$ and $k$ can be of ${\cal E}_x$ class or ${\cal M}_x$ class. And $a_{l,k}$ is the probability that a $l \in {\cal M}_x$ particle gets converted to a $ k \in {\cal M}_x$ particle.   In a similar way, type change transitions are allowed within ${\cal E}_x$  class, however, there are no type changes possible from one class to another. 
The generator matrix of such a TC-DBP has the following structure:
\begin{eqnarray}
\label{Eqn_Gen_TC-DBP}
 \left[
\begin{array}{cc}
A_{mx} & A_{mx,ex} \\
\textbf{0} & A_{ex}
\end{array}
\right],
\end{eqnarray} 
where $A_{mx}$ represents all the transitions between types belonging to class ${\cal M}_x$,   $A_{ex}$ represents all the transitions between types belonging to class ${\cal E}_x$, 
while  $A_{mx,ex}$ represents the transition between ${\cal M}_x$ and ${\cal E}_x$ (offspring of  class ${\cal E}_x$ produced  by class ${\cal M}_x$). 
The matrix $A_{mx}$ includes type change as well as real transitions as below:

\vspace{-4mm}
{\small \begin{eqnarray}
 {A}_{mx} := \left [
 \begin{array}{llll}
 \theta a_{1,1} + (1-\theta)  m_{1,1}-1    & \theta a_{1,2} + (1-\theta)   m_{1,2}   & \cdots   & \theta a_{1,M} + (1-\theta)   m_{1,M}      \\
  \theta a_{2,1} + (1-\theta)  m_{2,1}    & \theta a_{2,2} + (1-\theta)  m_{2,2} -1  & \cdots   & \theta a_{2,M} + (1-\theta)  m_{2,M}      \\
 & \vdots   \\
   \theta a_{M,1} + (1-\theta)   m_{M,1}    & \theta a_{M,2} + (1-\theta)  m_{M,2}   & \cdots   & \theta a_{M,M} + (1-\theta)   m_{M,M}   -1   \\
 \end{array}
 \right]. 
 \label{Eqn_TC_DBP}
\end{eqnarray}}  
The matrix $ {A}_{ex}  =   (( \theta a_{l, k} + (1-\theta) m^e_{l,k}  - 1_{l=k})) $ has exactly similar structure, with the only difference being that now $l, k \in {\cal E}_x$, but $a_{l,k}$ are the same (i.e., type change transitions are   the same). 
There are no type changes from one class to another, hence    $ {A}_{mx, ex}  =   ((   m_{l,k} )) $, with $l \in {\cal M}_x$ and $k \in {\cal E}_x$.
Our focus is to investigate the evolution of the number of shares of ${\cal E}_x$ class particles when started with a particle of ${\cal M}_x$. 

Note that this kind of branching processes can model various real-world applications,  including our social network example.  

\subsection{ Analysis: time evolution of the expected  net progeny}

We are now ready to study the time evolution of the expected net progeny  in TC-DBP.   Prior to that, we define the relevant terms. 
 
{\bf Two different notions of `total' progeny:} We emphasize that there are two different notions for the total progeny in   TC-DBP, as opposed to the standard one in the branching processes. 
   One may view the type-changing as the production of one offspring of a different type, and thereby adding one to the total progeny (for each type-change). This phenomenon is the usual way the total progeny is counted in standard BPs.  Alternatively, one may not view type-change as an offspring, which can lead to a different (new) notion of total progeny that counts only the new offspring.    
   For instance, as we already discussed,   in a social network  one needs only the count of total shares (the number of distinct users shared with the post of interest).   
{\it We refer to the progeny that does not count the type-changes as `net progeny', while the one that counts all the transitions as the usual `total progeny.'}  
   
   In this section, we derive the time evolution of the expected net progeny.  {\it Using the results on net progeny one can also obtain corresponding results for total progeny\footnote{To the best of our knowledge there are no results on total progeny of decomposable   branching processes.} in the standard decomposable branching
   processes (which  is of independent importance).} 
It is clear that the total progeny of  a decomposable branching process is obtained by substituting $\theta = 0$ in the expression for the expected net progeny of appropriate TC-DBP.  We obtain these results by using slight tweaks of the existing methods to analyze such branching processes; {\it in particular we study the functions representing the time evolution of  expected net progeny as  fixed points in some appropriate Banach  spaces  and derive the required analysis by obtaining the approximate fixed point solutions  (for applications like our social network). }

The  net (total) progeny at a time instance, say $t$,  represents the total accumulated population (i.e., including the dying particles) of all types till $t$, without (with) considering the type-changes. We study the evolution of the  net progeny for both exclusive class and mixed class particles.

 When one starts only with   'exclusive' parents, i.e.,  parents from class ${\cal E}_x$ and consider expected net progeny of exclusive class,  then the resulting branching process is well known irreducible multi-type continuous time branching process (e.g., viral branching process representing the propagation of single post in Par1-\cite{Ranbir2}). We derived the net progeny (a.k.a., number of shares)  in Part-I using standard tools, but here we would require the same without considering the count of initial population.  That is provided in Appendix B and the solution in general case   is given by (when the matrix is invertible):
\begin{eqnarray} 
\boxed{
 {\bf y}^e (t)   =  \left (  e^{ A_{ex} t }  - I  \right )  \lambda_\nu  (1-\theta ) A_{ex}^{-1}    \left [  
\begin{array}{lll}
  \sum_{k \in {\cal E}_x  }  m^e_{1,k}  \\
    \sum_{k \in {\cal E}_x  }  m^e_{2,k}     \\
    \vdots \\
      \sum_{k \in {\cal E}_x  }  m^e_{N,k}  
\end{array}
\right ]  .}   \nonumber
\end{eqnarray}
For the special case of OSNs, we have the following simplification (see Appendix B):
\begin{eqnarray}
\label{Eqn_ye_special_case_in_text}
y^e_l (t)  \approx ( e^{\alpha_e t}  - 1)   h_l^e \mbox{ with }
  h_l^e =  \frac{  \lambda_\nu  (1-\theta)   r_l   m \eta_1    }{ \alpha_e }   \mbox{ for all } l \in {\cal E}_x,
\end{eqnarray}and $  
  \alpha_e \approx  m \eta_1 \sum_{l}  r_l    \rho_l  + \lambda_\nu \theta \Delta_r  -1, $ when $r_{l+1} / r_l = \Delta_r$ and ${\bar \rho}_l = \rho_l$ for all $l$. The above result is true even if $ A_{ex}$ is not invertible;  further when  $\alpha_e > 0$, i.e.,  when the exclusive can  get  viral,   then $\alpha_e$ is the Perron root of $A_{ex}$ (see Appendix B).

We now focus on investigating the evolution of the expected net progeny of  exclusive class when the process starts with a mixed class particle. We obtain this by first deriving appropriate fixed point equations.

\subsubsection{Derivation of an appropriate fixed point (FP) equation}
 Denote by $ y_{l} (t) $ the  expected net progeny   of ${\cal E}_x$ class till time $t$ when the process is initiated with one type-$l$ particle of  ${\cal M}_x$ class, and $ \textbf{y}(t) = \{y_{l} (t)\}_l $ represents the `net progeny' vector till time $t$ and ${\bf y} = {\bf y} (\cdot)$ represents the vector of  net progeny time evolution.
We will show that ${\bf y} $ satisfies a fixed point equation in an appropriate functional space, i.e.,  ${\bf y} = {\bf z}$ where  $\textbf{z}: = \{z_l  (\cdot ) \}_{l}$ represent finite number of waveforms on time interval $[0, \infty)$ and satisfies an appropriate fixed point equation,   $z_{l} (t)  =  G_l (\textbf{z}) (t)$  (for all $l$,  $t$). We arrive at the fixed point function ${\bf G} = \{G_l\}$ by conditioning on the events related to the first transition epoch. Let the random variable $\tau$ represent the time instance of the first transition epoch, which  is exponentially distributed  with parameter $\lambda+\nu$.   Conditioning on the first transition events, we observe that the net progeny  ${\bf y} (\cdot) $ should satisfy the following fixed point equation (for all $l$ and $t$):
\begin{eqnarray*}
y_{l} (t) = G_l ({\bf y}  ) (t)   &:=&  \theta \int_0^t    \sum_{k \in {\cal M}_x   } a_{l, k}  y_k (t-\tau)   (\lambda+\nu) e^{- (\lambda+\nu) \tau } d\tau   \\
&& + \ \  (1-\theta ) \int_0^t    \sum_{k \in {\cal M}_x  }  m_{l,k} (1+ y_k (t-\tau)  )  (\lambda+\nu) e^{- (\lambda+\nu) \tau } d\tau   \\
&&
+ \ \ (1-\theta)  \int_{0}^t   \sum_{k  \in  {\cal E}_x }  m_{l,k}  ( 1 +  y^e_k (t-\tau) )   \   (\lambda+\nu)e^{- (\lambda+\nu) \tau } d\tau. 
\end{eqnarray*}
The above is due to the following reasons:
\begin{itemize}
\item The type-$l$ undergoes a shift transition w.p. $\theta$, its type gets changed to type-$k$ of the  same class  (i.e., $ {\cal M}_x   $).
\item The type-$l$ undergoes a share transition w.p. $1-\theta$,  it  produces $m_{l,k}$ offspring belonging to either class  $ {\cal M}_x  $ or $ {\cal E}_x$. As per example, it produces particles of $ {\cal M}_x  $ when  both post-${\bf P}$ and post-${\bf Q}$ are shared, whereas particles of $ {\cal E}_x  $ are produced when only one of the posts is shared.
\end{itemize}

{\it 
In Appendix C, we showed that  ${\bf G}( \cdot)$ is a contraction mapping and hence has unique fixed point solution.}
In view of the uniqueness, it suffices to obtain any solution of ${\bf G}$, which is considered immediately.

 \subsubsection{Solution of the fixed point equation}

The solution of the fixed point equation is derived in Appendix C and it equals (when the matrices are invertible)
\begin{eqnarray}
 \label{Eqn_net_progeny_mixed}
{\bf y} (t) &=&    \left ( e^{A_{mx} t}   -I  \right ) {\bf c_v}_0 +  e^{A_{mx} t} \int_0^t  e^{- A_{mx} s }  A_{mx, ex}      e^{ A_{ex} s }     A_{ex}^{-1}  {\bf c_v}_1  ds \\
{\bf c_v}_0 &=&  A_{mx}^{-1} \lambda_\nu  (1-\theta )  \left(    \left [  
\begin{array}{lll}
  \sum_{k \in {\cal M}_x  }  m_{1,k}  \\
    \sum_{k \in {\cal M}_x  }  m_{2,k}     \\
    \vdots \\
      \sum_{k \in {\cal M}_x  }  m_{N,k}  
\end{array}
\right  ] +   A_{mx, ex}  \left (   {\bf 1}    -     A_{ex}^{-1}  {\bf c_v}_1    \right )  \right ) \nonumber \\ 
{\bf c_v}_1 &=&   \lambda_\nu (1-\theta)  \left [  
\begin{array}{lll}
  \sum_{k \in {\cal E}_x  }  m^e_{1,k}  \\
    \sum_{k \in {\cal E}_x  }  m^e_{2,k}     \\
    \vdots \\
      \sum_{k \in {\cal E}_x  }  m^e_{N,k}  
\end{array}
\right ] .
\end{eqnarray} 
We then showed that the net progeny has the following simplified    form for our OSN example (details in Appendix C,  and approximation is good as $N \to \infty$):
\begin{eqnarray}
 y_{l} (t) \approx  g_l + h_l e^{\alpha_e t}+ o_le^{{\bar \alpha} t} ,  \label{Shares2}
\end{eqnarray}for some appropriate coefficients $\{h_l, o_l, g_l\}$. 
As in exclusive case, this approximation is true even the matrices are not invertible, we would only require that some  eigen values are positive (e.g., Case I in Appendix C).
We discuss more details of this representation in the coming sections.

Thus, we notice that for decomposable branching processes, the growth rates of the  current population as well as the total shares (\ref{Shares2}) are  influenced by two distinct exponential functions.  
Also,  both of them (like in other variants of the branching process) are influenced by the same growth patterns. 

We now return to our OSN example.

\section{CP-wise performance measures}   
Mixed-type TLs keep producing the exclusive-types  as well as their own type TLs till driven to extinction, i.e., when none of the TLs contain both the CP posts. Once the mixed-types get extinct,  the leftover exclusive-types do not influence each other (matrix (\ref{Eqn_Matrix_Mix_Generator})), and hence, they evolve on their own.  
  Nevertheless,  their survival/growth depends upon the effects created by mixed-types before death.  When the mixed population gets viral,  total CP population (sum of exclusive-type and mixed-type TLs having that particular CP's post) is clearly influenced by mixed-types. The mixed population thus gives an impetus to the propagation of exclusive CP-1 and CP-2 posts with different degrees, and consequently, introduces competition between these posts for relative visibility. In other words, the more the number of  exclusive-type TLs generated by the mixed-types, the better it is for the corresponding CP (owning the said exclusive-type TLs).  
To summarize, the evolution of the population corresponding to a particular CP  depends upon the competition regardless of whether the source of the competition (i.e., mixed-types) dies out or not.

Recall that the underlying MTBP is decomposable. This MTBP behaves significantly different from the MTBP in the single CP scenario. Here it may happen that the population corresponding to a particular CP gets extinct with probability one while the other can get viral with positive probability. Further, they can have different growth rates in the event of virality. 
\subsection{CP-wise extinction probabilities}
We say \CPj is extinct when all the mixed-type and  exclusive CP-$j$ type TLs get extinct.  
By Lemma 1 of Part-I \cite{Ranbir2}, the sub-matrix $A_{ex}^j$ of (\ref{Eqn_Matrix_Mix_Generator}) is irreducible.  Thus, all exclusive-type TLs of one CP survive/die together when the process starts with exclusive-type TL of the same CP. And the same is the case for mixed population when started with a mixed type TL, as matrix $A_{mx}$ is irreducible.  With   ${\bf e}_{l,k}$ as the unit vector with one only at $(l,k)$ position where $l, k =l+1 $ or $l-1$, we define the extinction probability of \CPj as below:
 $$ q_{l,k}^j := P\Big(\textbf{X}^j_{ex}(t) = \textbf{0},\ \textbf{X}_{mx} (t)  = \textbf{0} \mbox{ for some  } t > 0 \Big | {\bf X} (0) = {\bf e}_{l,k} \Big).
 $$
  $$\mbox{Let\ } {\bf q}^j := \left \{ {\bf q}^j_{ex}, {\bf q}^j_{mx1}, {\bf q}^j_{mx2}  \right \} \ \mbox{with\ \ }{\bf q}^j_{ex}  := \{q^j_{l, 0}\}_l,  {\bf q}^j_{mx1} := \{q^j_{l, l+1}\}_l \ \mbox{and\ } {\bf q}^j_{mx2} := \{q^j_{ l+1, l}\}_l.$$
By  conditioning again on the events of first transition, we obtain ${\bf q}_{mx1}^1$ via fixed point (FP) equations:
{\normalsize{
 \begin{eqnarray*}
q_{l,l+1}^1  & =& 
\theta \Big ( q_{l+1,l+2}^1  \mathbb{1}_{\{l < N-1\}} + \mathbb{1}_{\{l < N-1 \}} q_{N, 0}^1 \Big)  +  (1-\theta)(1-r_l)  
\\
&+ &
  (1-\theta)  r_l\Bigg [ (1-\delta) \bpgftwo\left(\textbf{q}^1_{ex},\eta_1 \right)  
\bpgftwo\left(\textbf{q}^1_{ex},\eta_1(1-\eta_2) \right)  
 \\
 &  & \hspace{-1cm} +\ \ \delta\ \Bigg ( p\bpgftwo\left(\textbf{q}^1_{mx2},\eta_{1 2} \right)   +  (1-p)  \bpgftwo\left(\textbf{q}^1_{mx1},\eta_{12} \right) \Bigg )    \Bigg ] ,  \mbox{ with }
  \bpgftwo \left(\textbf{s}', \eta \right) := \sum_{i =1}^{N-1}f\left( s_{i}', \eta,\beta \right)\bar{\rho_i}
 \nonumber
\end{eqnarray*}}}
and $\eta_{12} := \eta_1 \eta_2.$
One can write the fixed point (FP) equations for ${\bf q}_{mx2}^1$, ${\bf q}_{mx1}^2$ and ${\bf q}_{mx2}^2$ in a similar way. And the expression for ${\bf q}_{ex}^j$ starting from exclusive CP types is same as that in the  single CP scenario (\cite{Ranbir2}).   We further have the following:
\begin{lemma}
\label{Lemma_mixed_unique_extn}
When ${\bf q}^j_{ex} < {\bf 1} = (1, \cdots, 1)$,  we have unique solution   in the interior of \ $[0, 1]^{2N-2}$,
i.e., ${\bf q}^j_{mx1}  < {\bf 1}$,  $ {\bf q}_{mx2}^j < {\bf 1}$.
  When ${\bf q}^j_{ex} = {\bf 1}$,  we have that $({\bf q}^j_{mx1}, {\bf q}_{mx2}^j) = {\bf 1}$ is the unique solution,
under extra assumption that $\rho_N = 0$ and  ${\bar \rho}_i = \rho_i$ for all $i < N$.  
\end{lemma}
\textbf{Proof:} The proof is given in Appendix A. \eop

Thus when the exclusive types get extinct with probability one when started with exclusive types,  then  they get extinct with probability one even upon starting with  mixed types. When the exclusive can survive upon starting with their own types, then the exclusives can survive with positive probability even when started with mixed types.

 \subsection{Evolution of  exclusive-type NU-TLs in  presence of mixed-types }
 The evolution of exclusive CP population when started with exclusive-type TLs is same as in single CP scenario, and we have the following result  for sufficiently large $t$
 $$ 
 \textbf{X}_{ex}^i(t)^T \textbf{v}^{i} \approx W_ie^{\alpha_i t}, \ \mbox{ and also \ \ } \textbf{X}_{ex}^i(t)  \approx W_ie^{\alpha_i t} \textbf{u}^{i} ; \ i = 1,2; \ \mbox{ where}
 $$
  \begin{itemize}
  \item[-]  $W_i$ is a non negative random variable, with $P(W_i = 0) = $ extinction probability (\cite{AthreyaPaper})
  \item[-]  $\textbf{u}^{i},\textbf{v}^{i}$ are the normalized left and right eigenvectors corresponding to the eigenvalue $\alpha_i$ (Perron root) of the matrix $A_{ii}$ for $i = 1,2.$
  \end{itemize}
We derive the evolution of exclusive-types when started with a mixed-type TL.
Consider any $i$, by    \cite[Theorem 3]{arx} we have the following result, when $\alpha_{mx} \ne \alpha_i$. 
  \begin{lemma} \label{Lemma_DecomBP}
  Let  $( \textbf{X}_{ex}^i(t), \textbf{X}_{mx}  (t ))$ be the number of unread TLs having post of CP-$i$ at various levels  at time $t$ as before and let $\mathcal{F}^i_t := \sigma \big\{\textbf{X}_{mx}(t'), \textbf{X}^i_{ex}(t'); \ t' \le t \big\}$ be the natural sigma algebra.  When this process starts with a mixed-type TL, then the stochastic process:
{\normalfont \begin{eqnarray}
 \left\{\textbf{X}^i_{ex}(t) \cdot  \textbf{v}^{i}  e^{-\alpha_{i} t} \mbox{ + }  e^{-\alpha_{i} t} \textbf{X}_{mx}(t)  \big(\alpha_i I - A_{mx} \big)^{-1} A^i_{mx,ex} \cdot \textbf{v}^{i};\ \ \mathcal{F}^i_t; \  t \ge 0 \right\}
 \end{eqnarray}}
is a martingale, where the vectors {\normalfont $\textbf{v}^{i}, \textbf{u}^{i}$} are as defined above.  \eop
 \end{lemma}
\textbf{Remarks:} 1) Thus even in the presence of mixed population, the population corresponding to exclusive types eventually evolves with a growth rate given by $\alpha_i$;  when multiplied with $e^{-\alpha_i t}$ the otherwise exploding process (exploding with time) converges to a limit.\\  2)  By Martingale property,
the   weighted sum of expected values of  the individual components  (all are  column vectors and $\cdot$ is the dot product),
{\normalfont \begin{eqnarray*}
 E\big[\textbf{X}^i_{ex}(t)     \big] \cdot  \textbf{v}^{i}   +  E [\textbf{X}_{mx}(t) ] \big(\alpha_i I - A_{mx} \big)^{-1} A^i_{mx,ex}  \cdot  \textbf{v}^{i}  \hspace{-20mm} \\
 & = & e^{\alpha_{1} t} \left ( \textbf{X}_{mx}(0)  \big(\alpha_1 I - A_{mx} \big)^{-1} A^1_{mx,ex} \cdot  \textbf{v}^{i}  \right ).
\end{eqnarray*}}
\ignore{
This is due to the following. From Lemma \ref{Lemma_DecomBP}  (without loss of generality consider $i=1$) and since the expected value of martingales are constant with respect to time, 

\vspace{-5mm}
{\fontsize{9.5}{6}{
 \begin{eqnarray*}
 E\Big[\textbf{X}^1_{exm}(t).\textbf{v}^{1}  e^{-\alpha_{1} t} +  e^{-\alpha_{1} t} \textbf{X}_{mx}(t)  \big(\alpha_1 I - A_{mx} \big)^{-1} A^1_{mx,ex}\textbf{v}^{1} \Big] & = &
\textbf{X}_{mx}(0)  \big(\alpha_1 I - A_{mx} \big)^{-1} A^1_{mx,ex}\textbf{v}^{1}
\\
 E\Big[\textbf{X}^1_{exm}(t) \Big]\textbf{v}^{1}.\textbf{u}^{1}  +  E \Big[\textbf{X}_{mx}(t)  \big(\alpha_1 I - A_{mx} \big)^{-1} A^1_{mx,ex} \Big]\textbf{v}^{1}. \textbf{u}^{1}& = & e^{\alpha_{1} t}\textbf{X}_{mx}(0)  \big(\alpha_1 I - A_{mx} \big)^{-1} A^1_{mx,ex} \textbf{v}^{1}.\textbf{u}^{1} 
\\
 E\Big[\textbf{X}^1_{exm}(t) \Big]  +  E \Big[\textbf{X}_{mx}(t)   \Big] \big(\alpha_1 I - A_{mx} \big)^{-1} A^1_{mx,ex} & = & e^{\alpha_{1} t}\textbf{X}_{mx}(0)  \big(\alpha_1 I - A_{mx} \big)^{-1} A^1_{mx,ex}
 \end{eqnarray*}}} 
 Now using $  E \Big[\textbf{X}_{mx}(t)   \Big]\big(\alpha_1 I - A_{mx} \big)^{-1} A^1_{mx,ex} \approx  \sum_{l} v_l^{mx}e^{\alpha_{mx} t}\textbf{u}_{mx}\big(\alpha_1 I - A_{mx} \big)^{-1} A^1_{mx,ex}  $
 \begin{eqnarray*}
   E\Big[\textbf{X}^1_{exm}(t) \Big]  + \sum_{l} v_l^{mx}e^{\alpha_{mx} t}\textbf{u}_{mx}\big(\alpha_1 I - A_{mx} \big)^{-1} A^1_{mx,ex} & \approx & e^{\alpha_{1} t}\textbf{X}_{mx}(0)  \big(\alpha_1 I - A_{mx} \big)^{-1} A^1_{mx,ex}.  \hspace{4mm} 
 \end{eqnarray*}
}
From Theorem \ref{Thm_mixed}.{\it iv} (under appropriate second moment conditions), we have for all large $t$:
{\normalfont \begin{eqnarray*}
 E\big[\textbf{X}^i_{ex}(t)     \big] \cdot  \textbf{v}^{i}   & \approx&   c_{ons1}   e^{\alpha_{mx}  t} + c_{ons2} e^{\alpha_{1} t}\\
 c_{ons1} & =&  -  \left (  E [W_{mx} ] \textbf{u}_{mx}   \big(\alpha_i I - A_{mx} \big)^{-1} A^i_{mx,ex}  \cdot  \textbf{v}^{i} \right )   \mbox{ and } \\
c_{ons2}  &=&  \left ( \textbf{X}_{mx}(0)  \big(\alpha_1 I - A_{mx} \big)^{-1} A^1_{mx,ex} \cdot  \textbf{v}^{i}  \right ).
\end{eqnarray*}}
Thus the growth of the expected number  of TLs with a CPi-post ($ E\big[\textbf{X}^i_{ex}(t) \big]$), when the process starts with mixed-type TLs,    is governed  by    two exponential curves.   
\ignore{
\begin{eqnarray*}
 E\left[\textbf{X}^1_{exm}(l,t) \right]   & = & \bar{h}_{l} e^{\alpha_{1} t} - \bar{o}_l e^{\alpha_{mx} t}; \ \mbox{where}
\end{eqnarray*}
$\bar{h}_{l} $  and $\bar{o}_{l} $ are the $l$-th component of the vectors $ \textbf{X}_{mx}(0)  \big(\alpha_1 I - A_{mx} \big)^{-1} A^1_{mx,ex}$ and $  \sum_{l} v_l^{mx}\textbf{u}_{mx}\big(\alpha_1 I - A_{mx} \big)^{-1} A^1_{mx,ex} $ respectively.}
 \subsection{ Evolution of the expected number of shares} 
\label{sec_twocp_shares} 
In the single CP scenario, we derived the evolution of the expected number of shares, which serves as an important performance measure for the spread of the content. We now derive the time evolution of the expected number of shares in the two-CP scenario; which, in this case, is instrumental in obtaining performance measures such as relative visibility. Let $Y^j_{l,k}(t)  $ be the total number of shares of \CPj post  till time $t$ and  $y^j_{l,k}(t) $ represents its expected value, when started with one TL of the type $(l,k)$ with $k = l + 1$ or $l-1 $. Note that these shares include mixed-type shares as well as exclusive-type shares belonging to CP-$j.$ 
We present the evolution for the number of shares to CPs' posts in non-viral and viral scenarios.
 
 \subsubsection {Number of shares in non-viral scenario:}
With $m < 1$, any post (CP1-Post or CP2-Post)  gets extinct with probability one. Even with $m > 1$,  the CPj-post can  get extinct with probability one,  depending on other parameters.  We refer to this as non-viral scenario.
We obtain the total expected shares (before extinction) when the process starts with a TL of mixed/exclusive type. Recall that  $\{y_{l,k }^j\} $ is expected number of shares of \CPj post when started with $(l,k)$ type TL

$$y_{l,k}^j = E[\lim_{t \to \infty }Y^j (t) | {\bf X}(0)  = {\bf e}_{l,k}]; \ \ j = 1,2.$$
These $\{ y_{l,k}^j \}$  can  be obtained by  solving appropriate FP equations (below) as in the single CP scenario (\cite{Ranbir2}). 
Without loss of generality we consider shares of CP-1.
Let ${\bf y}_{mx1}^j :=  \{y^j_{l, l+1} \} $,   ${\bf y}_{mx2}^j :=  \{y^j_{ l+1, l} \} $ and ${\bf y}^j_{mx}  := {\bf y}^j_{mx1}+{\bf y}^j_{mx2}$,   by appropriate conditioning
\begin{eqnarray*}
y_{l,l+1} & = & \theta \Big ( \mathbb{1}_{\{l < N-1\}} y_{l+1,l+2} + \mathbb{1}_{\{l = N-1\}} y_{N,0}  \Big )  +(1-\theta)r_l (1-\delta) m\eta_1 (1 + \textbf{y}_{ex1}.\bm{\bar{\rho} })  \nonumber \\
&& + \ (1-\theta)r_l\delta  m\eta_1\Bigg [(1-\eta_2)(1+ \textbf{y}_{ex1}.\bm{\bar{\rho}})  \nonumber  
+\eta_2(1+ p\textbf{y}_{mx1}.\bm{\bar{\rho}} +  (1-p)\textbf{y}_{mx2}.\bm{\bar{\rho}})\Bigg ] \nonumber 
\end{eqnarray*}
where ${\bf  y}_{ex1} = \{ y_{1,0}^1, y_{2, 0}^1, \cdots, y_{N-1,0}\} $. And again for any $l < N$, 
{\footnotesize{
\begin{eqnarray}
y_{l+1,l} & = & \mathbb{1}_{\{l < N-1\}}\theta y_{l+2,l+1} 
 + \ (1-\theta)r_l\delta  m\eta_1\Bigg [(1-\eta_2)(1+ \textbf{y}_{ex1}.\bm{\bar{\rho}})   
 \   \eta_2 \big (1+ p\textbf{y}_{mx1}.\bm{\bar{\rho}} +  (1-p)\textbf{y}_{mx2}.\bm{\bar{\rho}} \big ) \Bigg ]. \label{Eqn_y_ll+1}
\end{eqnarray}}} 
For the special case with   $ \bar{\rho}_l = \ \bar{\rho}^l/\sum_{i=1}^{N-1} \bar{\rho}^i,\  r_l = d_1 d_2^l$
 and as $N \to \infty$, we have  (with $-j := 2 \mathbb{1}_{\{j=1\}} + 1 \mathbb{1}_{\{j=2\}} $)
\begin{eqnarray*}
\textbf{y}^j_{mx}.\bm{\bar{\rho}} 
& \to &
 \frac{\Big ( 2 c_j\delta \left[ (1 + \textbf{y}^j_{ex}.\bm{\bar{\rho} })(1-\eta_{-j}) + \eta_{-j}\right]  +  c_j (1-\delta) (1 + \textbf{y}^j_{ex}.\bm{\bar{\rho} }) \Big )O^*_{mx}}{1-c_{mx}  O^*_{mx}};  
\\
 O_{m_x}^* &= &  \frac{d_1  d_2  (1-\bar{\rho})}{\left(1-d_2 \bar{\rho}\right)\left(1-\theta d_2\right)} \ \mbox{and \ } \textbf{y}^1_{ex} = \big\{ y_{1,0},\cdots, y_{N-1,0}\big\}; \ \ \textbf{y}^2_{ex} = \big\{ y_{0,1},\cdots, y_{0,N-1}\big\}.
\end{eqnarray*} 
Here $ \textbf{y}^j_{exj}$  is similar to that in Part-I \cite{Ranbir2}, and  $\{ y_{l,k}^j \}$ with $k = l+1$ or $l-1$ 
can be computed  uniquely using $\{\textbf{y}^j_{mx} \idot \bm{\bar{\rho}} \}$ and equation (\ref{Eqn_y_ll+1}) and backward induction.  
The derivation of these limits and expressions are  provided in Appendix D.

\subsubsection{Number of shares in viral scenario}
In viral scenarios,   we have two sub-cases: 1) one when the mixed survives with non zero probability, and 2)  when the mixed population gets extinct, but exclusive types may survive.  In both these cases, the expected number of shares explode with time.   
This analysis can be obtained using net progeny of section \ref{sec_type_change}, which is provided by (\ref{Eqn_net_progeny_mixed}). 
We explain the results for the case with $p=0$ for ease of explanation, one can easily extend the results to general case using the results of Appendix B and C.

Consider without loss of generality the number of shares for CP1-post.  When one starts with $(l, l+1)$-type (and with $p=0$)   one should consider the following modelling details for using expression (\ref{Eqn_net_progeny_mixed}):
\begin{eqnarray*}
a_{l, k} := a_{(l,l+1), \ (k, k+1)} &=&  a_{(l,0), (k,0)} \ = \  \mathbb{1}_{ \left\{k =   l+1\right\} } {\mathbb 1}_{l < N}, \\ 
&& \hspace{20mm}\mbox{ i.e., either for both } l,k \in {\cal M}_x \mbox{ or  for both } l, k \in {\cal E}_x,\\
 m_{l,k} :=  m_{(l,l+1),  (k, k+1)} & =&   m \eta_1 \eta_2 \delta  {\bar \rho}_k  r_l,  \ \ \  \  \  \   \mbox{i.e.,  for  } l,k \in {\cal M}_x ,   \\
m_{l,k} :=   m_{(l,l+1),  (k,0)} &=&  m \eta_1  r_l \bar{\rho}_k
(1-\delta + \delta (1-\eta_{2})),  \ \   \mbox{ i.e.,  for } \ l \in {\cal M}_x, \ k \in{\cal E}_x \mbox{ and }  \\   
m^e_{l,k}  &:=&   m_{(l, 0),  (k, 0)}  \ =\  m \eta_1 r_l \rho_k,    \ \      \ \   \   \  \  \ \  \mbox{i.e.,  for }  \ l \in {\cal E}_x, \ k \in{\cal E}_x.
\end{eqnarray*}
 When one starts with $( l+1, l)$-type (and with $p=0$)   one should consider the following modelling details:
\begin{eqnarray*}
a_{l, k} := a_{(l+1, l), \ (k+1, k)} &=&  a_{(l,0), (k,0)} \ = \  \mathbb{1}_{ \left\{k =   l+1\right\} }1_{l < N}  ; \  \  \ \ \\
&& \hspace{20mm}\mbox{ i.e., either for both } l,k \in {\cal M}_x \mbox{ or  for both } l, k \in {\cal E}_x,\\
 m_{l,k} := m_{(l+1, l),  (k+1, k)} &=&   m \eta_1 \eta_2 \delta  {\bar \rho}_k  r_l,      \ \ \  \  \  \   \mbox{i.e.,  for  } l,k \in {\cal M}_x ,    \\
 m_{l,k} :=  m_{(l+1, l),  (k,0)}  &=&  m \eta_1 r_{l} \bar{\rho}_k
(1-\delta + \delta (1-\eta_{2}))   \   \mbox{ for } \ l \in {\cal M}_x, \ k \in{\cal E}_x  \mbox{ and }  \\   
    m^e_{l,k} & := & m_{(l, 0),  (k, 0)}  \ =\  m \eta_1 r_l \rho_k \ \   \   \  \  \ \ \ \ \   \mbox{ for }  \ l \in {\cal E}_x, \ k \in{\cal E}_x.
\end{eqnarray*}
We derived  much simplified expressions  for the same in Appendix C (case I and  II),  when $
{\bar \rho}_l = \rho_l$ for all $l$,  $p =0$  and when $r_{l+1}/r_l = \Delta_r$ for all $l < N$ and we reproduce the same here:\\
{\bf When  one starts with $(l, l+1)$ type particle}, irrespective of whether $\alpha_{mx} >0 $ or not, the expected number of shares/net progeny   evolves exactly as when one started with an exclusive particle at $l$ level. It is not difficult to see this equality with $p=0$, if one closely introspects the two evolutions (when one counts all shares that belonging to CP1). This case is studied as Case I in Appendix C and the final result is the following:
 \begin{eqnarray}
 y_{l,l+1}^1 (t) & \approx & 
  \frac{  \lambda_\nu  (1-\theta)   m \eta_1 r_l   }{ \alpha_1 } \left (  e^{\alpha_1 t} - 1  \right )  \mbox{ for all } l  , \mbox{ with } \label{EqnShare_mixed_non_viral}   \\
  \alpha_1 & \approx& \lambda_\nu   \left  ( \theta  \Delta_r - 1    + (1-\theta)  m \eta_1  \sum_k \rho_k  r_k  \right ).  \nonumber 
 \end{eqnarray}
 One can have another approximate solution for this sub-case (see Appendix C, Case 2 for more details).\\
{\bf When one starts with $(l+1,l)$ type particle}, 
 from Case II we have: 
\begin{eqnarray}
 y_{l+1, l}^1 (t)  & \approx &  h_{l} (e^{\alpha_1 t} -1 )+ o_{l} (e^{\alpha_{mx} t} - 1) \label{EqnShare2CPViral}  \mbox{ with }  \\
 %
o_l  &=&  r_l    \frac{ (1-\theta) \lambda_\nu     }{  (  {\alpha}_{mx} -\alpha_1)   {\alpha}_{mx} }  \bigg (     \lambda_\nu \Delta_r \theta + (1-\theta)   \lambda_\nu  m^2 \eta_1^2 \delta (1-\eta_2)+ \eta_2  \delta  )     \sum_{k   } {\bar \rho}_k  r_k  \nonumber \\
&& \hspace{40mm} - (\lambda_\nu +\alpha_1)  m \eta_1    \delta
 \bigg  ),     \nonumber  \\
  \\
 h_l  &=&  \frac{  \lambda_\nu  (1-\theta) m \eta_1\delta  r_l   -   o_l   {  \alpha }_{mx} }  { \alpha_1 }, \ \ \
  \nonumber   \\
{  \alpha }_{mx}  &\approx &  \lambda_\nu   \left ( \theta \Delta_r - 1   + (1-\theta) m\eta_1 \delta \eta_2 \sum_k {\bar \rho}_k r_k \right )    \mbox{ and }  \nonumber  \\
\alpha_1 & \approx& \lambda_\nu   \left  ( \theta  \Delta_r - 1    + (1-\theta)  m \eta_1  \sum_k \rho_k  r_k  \right ).  \nonumber 
\end{eqnarray}

One can write similar expressions for CP2-post and for arbitrary $p$.
As mentioned before, $\alpha_{mx} \le \alpha_j $ for $j=1,2$. One can argue that $ y_{l,l+1}^j (t) $ grows with rate $\alpha_j$ in the long run, i.e., the growth rate of the expected number of shares to \CPj  when mixed-type TLs get viral is carried by dominating rate ($\alpha_j$) which is same as the growth rate of  shares in the single CP case.

\subsubsection*{Numerical Evidence}  
\begin{figure}
\begin{minipage}{7.6cm}
\includegraphics[scale=0.25]{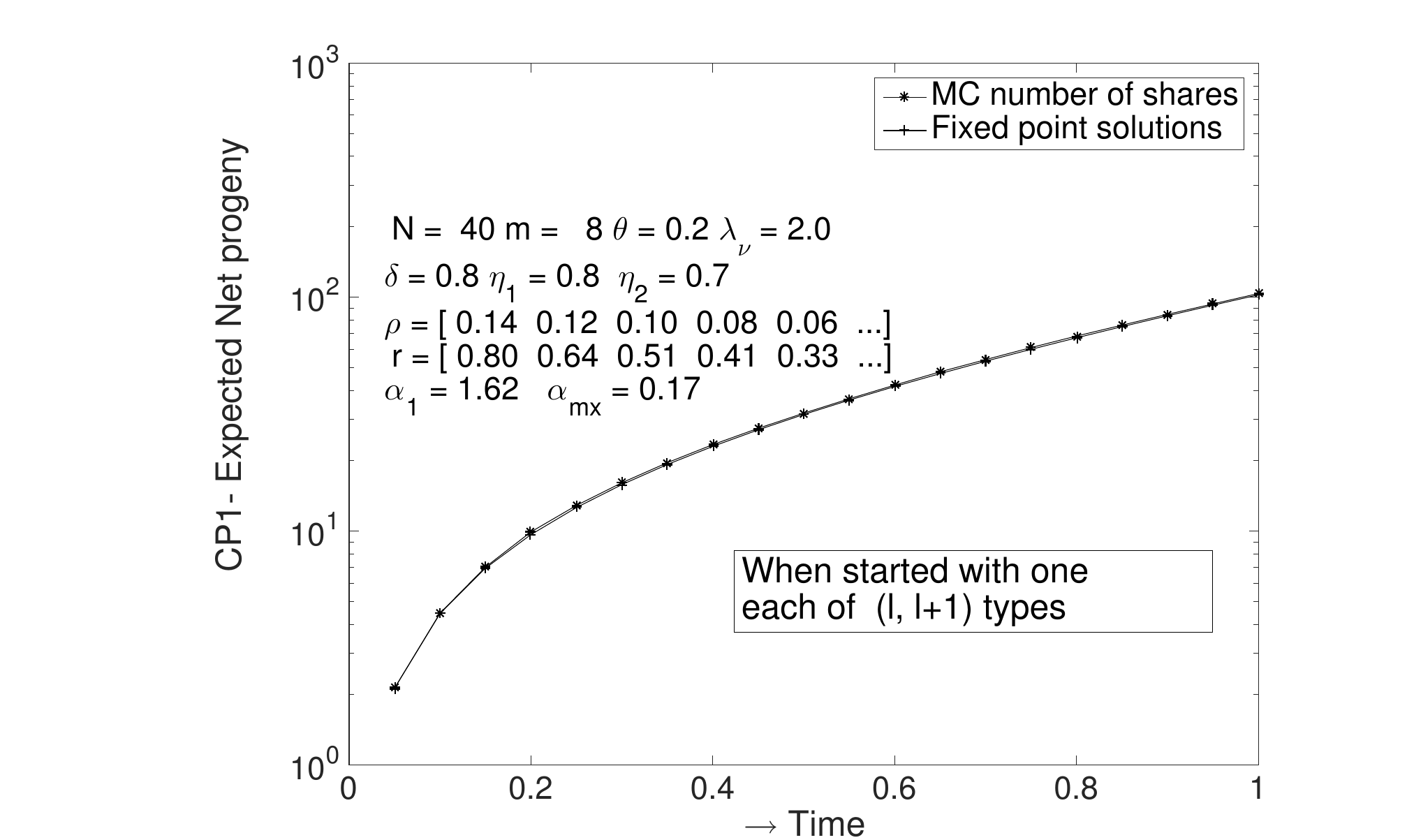}
\end{minipage} 
\begin{minipage}{6.6cm}
\includegraphics[scale=0.25]{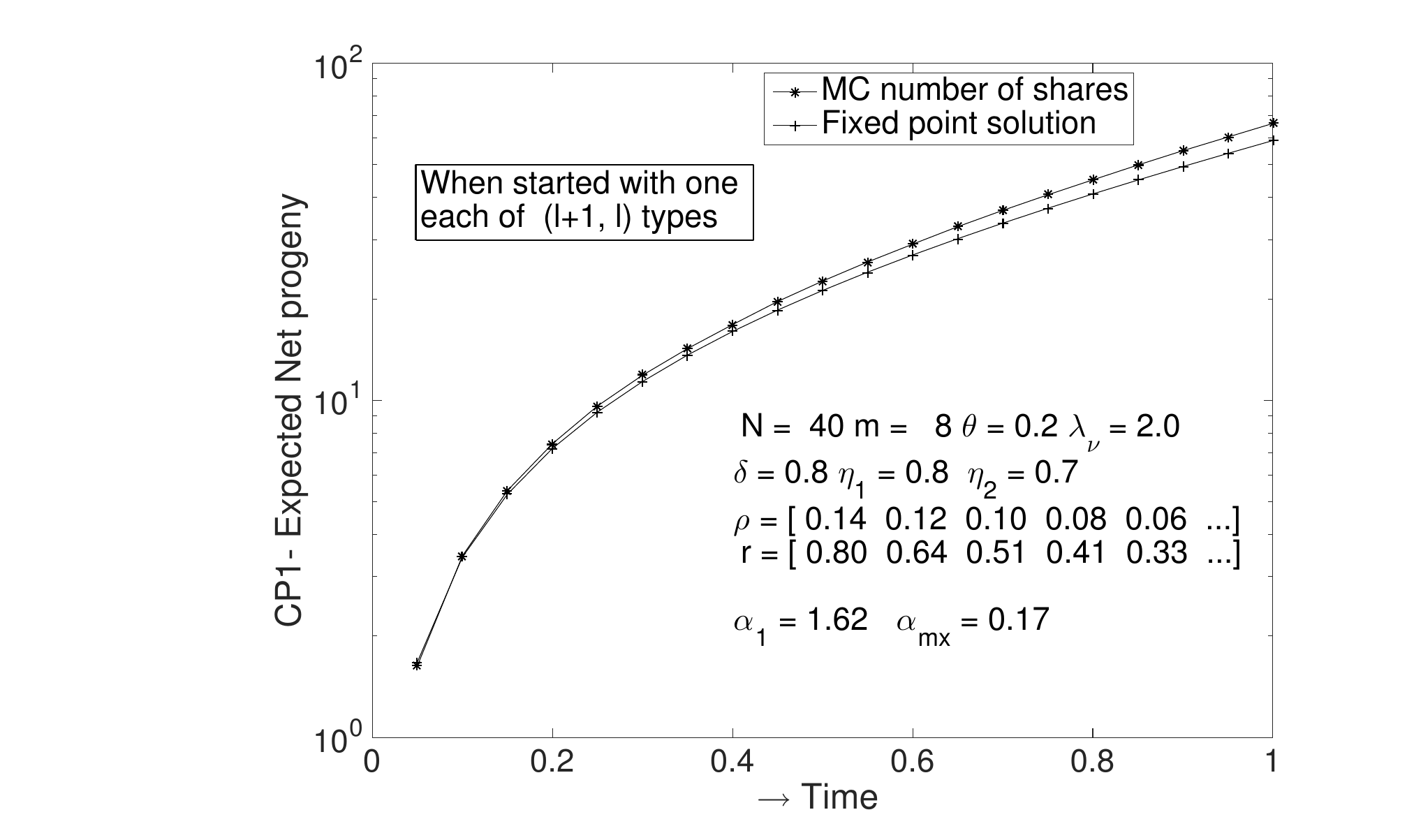}
\end{minipage}
\caption{Expected Net progeny evolution:  theory and MC estimates, when mixed  types also get viral \label{Fig_net_prog}}
\end{figure}

\begin{figure}
\begin{minipage}{7.6cm}
\includegraphics[scale=0.25]{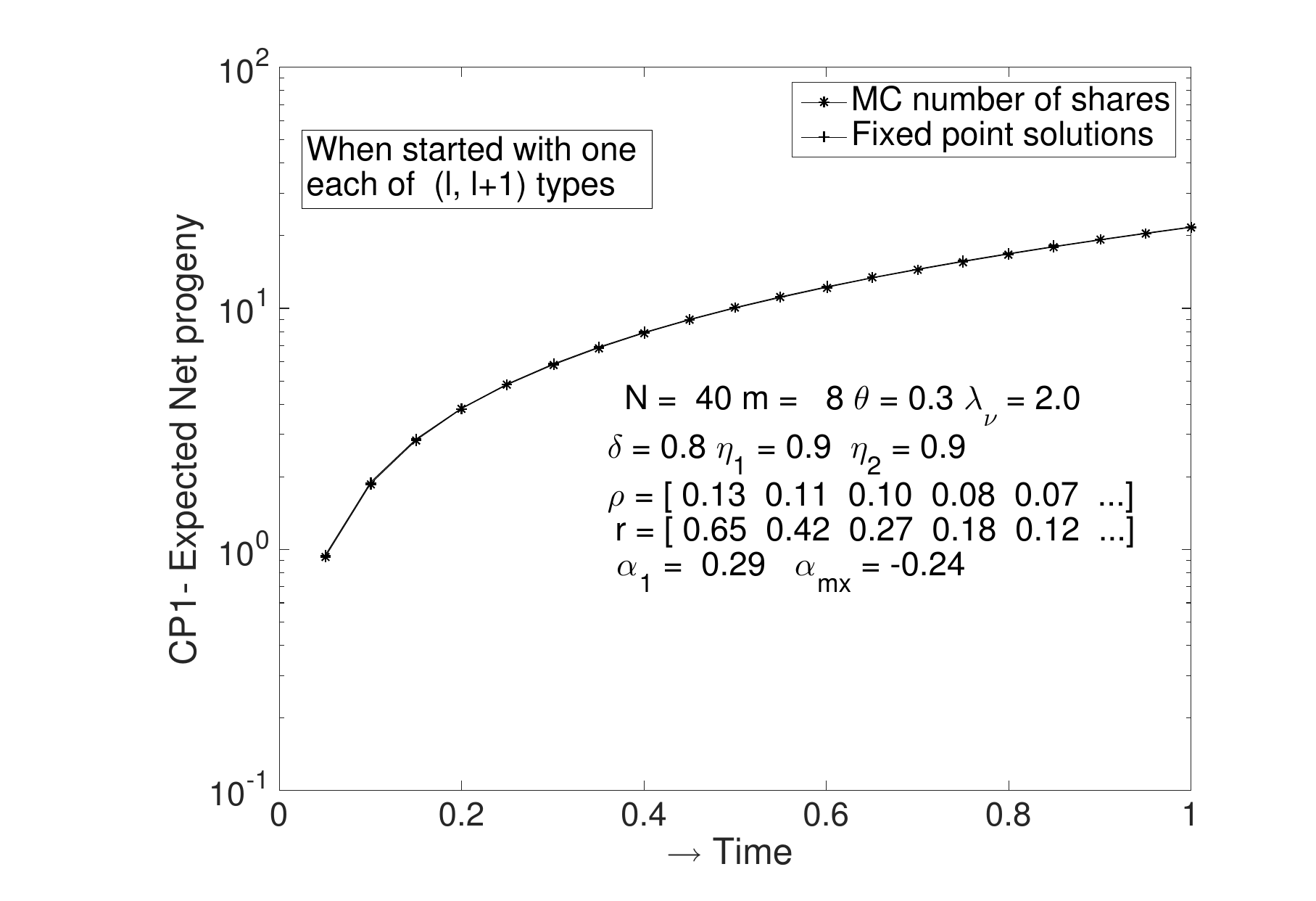}
\end{minipage} 
\begin{minipage}{6.6cm}
\includegraphics[scale=0.25]{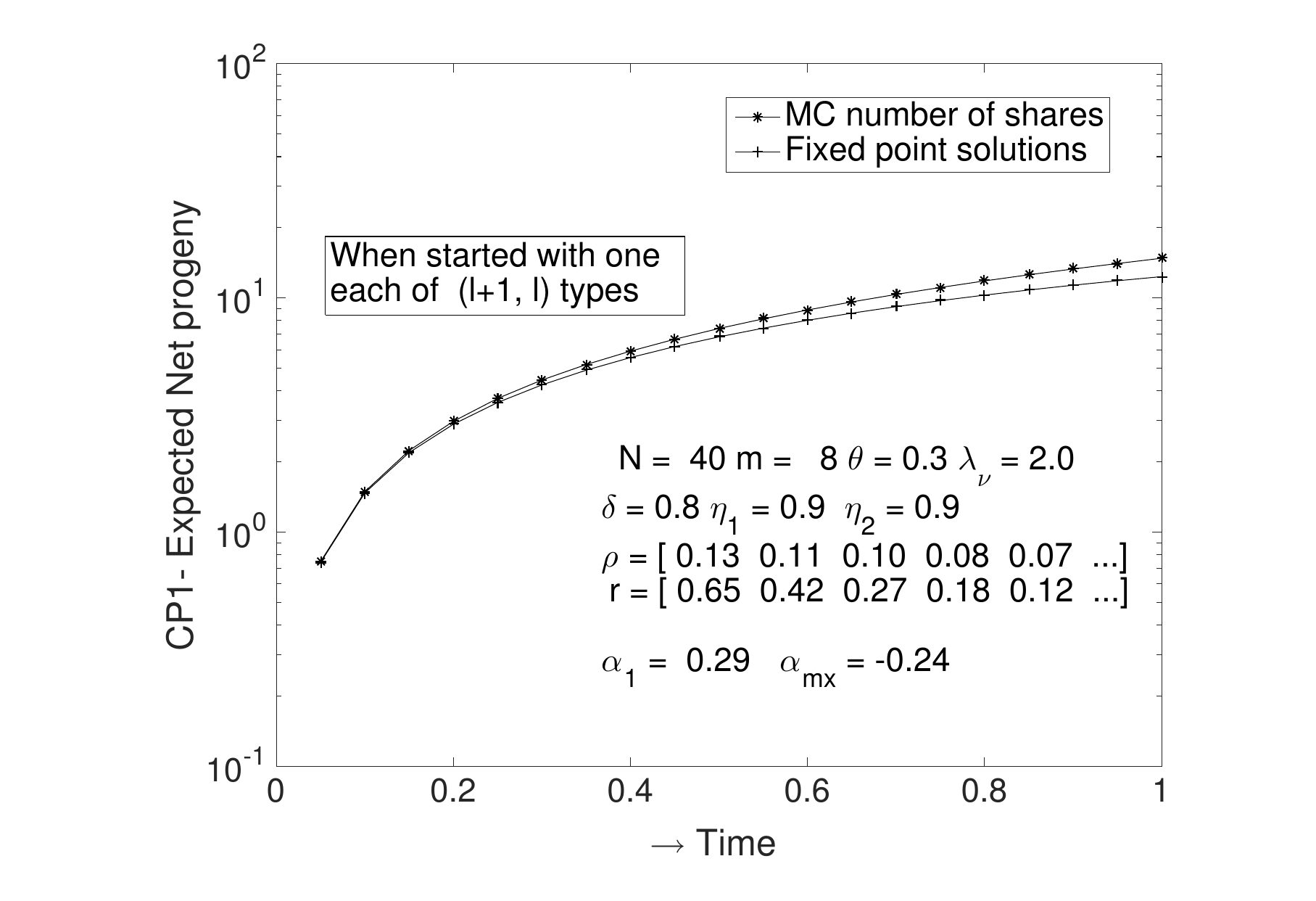}
\end{minipage}
\caption{Expected Net progeny evolution:  theory and MC estimates, when mixed  types  get extinct with probability one  \label{Fig_net_prog_2}}
\end{figure}

In Figures \ref{Fig_net_prog} and \ref{Fig_net_prog_2}, we consider Social network examples with details as given in the figure.  Here we plotted the Monte Carlo (MC) based estimates of expected net progeny (or the expected number of shares) at different time points and also the theoretical estimates  obtained using approximate fixed point solutions given by (\ref{EqnShare2CPViral}) and (\ref{EqnShare_mixed_non_viral}).  In these figures we are starting with either one each of $(l, l+1)$  (for all $l < N$) or one each of $(l+1, l)$ type, i.e.,  with a total of  $N-1$ particles. We observe that the theoretical approximate solutions well approximate the net progeny trajectories estimated by MC based  simulations.  Based on these examples and many more such examples  (not mentioned in this paper), we found that: a) the theoretical approximate solutions well approximate the MC based estimates when one starts with $\{(l, l+1)\}$ types for almost all sets of the parameters;  b) there is a good match between the two, for the sub-cases when one starts with $\{(l+1, l)\}$ types only for small $\theta$; the approximation is very good at small $\theta$ (for most of the cases with  $\theta < 0.2$).

Thus we derived various relevant performance measures related to content propagation in the presence of competing posts. One can use these performance measures for any relevant optimization or game theoretic problems. In the following section we use some of these measures for problems related to online auctions and viral marketing.

\section{Viral marketing and real time bidding} \label{Ch_oxn}
The performance measures obtained in the previous sections can be useful in many advertisement/campaign related objectives such as brand awareness, search engine optimization, maximizing the number of clicks to a post/advertisement (ad), etc. 
In this section, we will study online auctioning for advertisements in viral marketing using the performance measures as obtained in the previous sections.

The publishers of OSNs sell the advertisement inventory/space to various content providers (CPs) via auction mechanism commonly known as \textit{real time bidding} (\cite{BidEst}). For example, Facebook auctions billions of advertisement space inventory every day, and the advertisements  (ads) of the winners are served. Real-time bidding enables the CPs to automatically submit their bids in real time, and the advertisement of the highest worth (based on bid amount and its performance) is thus served. By virtue of auctioning, a natural competition occurs among the CPs for winning auctions. Further, a content provider (CP) has to win the auction to get sufficient number of seed (initial) timelines. 
The virality/sharing of the post further depends upon the quality of the advertisement/post (recall the post quality factor $\eta$). On summarizing, the CP has to invest in two aspects: a) the bid amount to win the auction, and b) the amount spent to the design of the post ($\eta$).
Recall that designing of a post could include providing authentic information about your services/products, or providing quality content,  or giving offers, etc.  Inappropriately tailored post can make users lose interest in the post, and thereby reducing the virality chances.

Content providers (CPs) typically have wide-ranging objectives while advertising on OSNs. For example, a CP may be in interested in enhancing the brand awareness of its products. Brand awareness plays a central role in users' decision making for a purchase. Such an objective is achieved if the brand promotional post gets viral. 
Recall, we say a post gets viral if it spreads on a massive scale via its sharing among the users.  Given that a post gets viral, a CP may be interested in knowing how fast the post spreads, i.e., the rate of virality. Other objectives, a CP may be interested in, include:  maximizing the number of clicks on its post, improving its reputation, increasing its presence in the marketplace, etc. 

In previous sections, we derived some of these performance measures. For example, we obtained the time evolution of the number of shares and NU-TLs which characterize the rate of virality. We also obtained the expression for the probability of virality. On the other hand, in non-viral (sure extinction) scenarios, we computed the expected number of total shares before extinction. We provide the explicit expressions for some of the performance measures as a function of controllable parameters while others are represented as the solutions of appropriate  FP (fixed point) equations. One can use these measures to study a relevant optimization problem taking auctions into account. In particular, and without loss of generality, we take the expected shares/NU-TLs as an indicative of the performance of CP's posts.

\subsection{Games with auction}
We now discuss this problem in the context of direct competition between the two CPs. We formulate a game, as before, considering online auctions additionally.

As mentioned before, a natural competition is induced between the CPs due to the propagation of the competing posts through the same OSN. As before, we also have competition due to winning auctions. We study this competition by formulating an appropriate game theoretic framework. We begin with the description of the utility functions of the CPs in the game. 
 As before, we take the utility of the CP as the ``number of TLs having its post at time $t$, i.e., $\textbf{X}(t)$" as one among those choices. Further, each CP incurs twofold costs as before: 1) winning the auction, and 2) cost for post quality. The net utility is thus obtained by subtracting this total cost from the revenue generated from $\textbf{X}(t)$. As there are two CPs, we assume that two sequential auctions are conducted and each CP participates only in one auction. The auctions so conducted can result in the following outcomes: 1)  each CP wins its own auction, 2) CP-1 wins while CP-2 loses the auction or vice versa, and 3) both CPs lose their respective auctions.
 We disregard the third outcome as we do not have the post of interest.
  The second outcomes gives rise to the propagation of the post corresponding to one specific CP only, whose analysis is carried out in a single CP scenario. Whereas the first outcome leads to the propagation of both the CPs' posts, i.e., production of the mixed population in addition to the exclusive CP types populations. 
 With this, we now describe the net utility, say $\mathbf{C}_i\left(x_1,x_2, \eta_1, \eta_2 \right)$ for $i = 1,2$, derived by the CP-$i$ as below:
 \begin{eqnarray*}
\mathbf{C}_1\left(x_1,x_2, \eta_1, \eta_2 \right)  & = &
  \Big(  \log E\Big(\sum {X^1_{ex}}(t) \Big) -  \kappa_2(x_1 + \kappa_1 \eta_1) \Big)P( \mathbf{B} < x_1 \eta_1) P( \mathbf{B}> x_2 \eta_2)  
 \\[0.15cm] 
& & \hspace*{-3cm} \  +\ \  0 \times P(\mathbf{B} > x_1 \eta_1) + \Big(\log E\Big(\sum {X^1_{exm}}(t) \Big) -  \kappa_2(x_1 + \kappa_1 \eta_1) \Big)P( \mathbf{B} < x_1 \eta_1) P(\mathbf{B} < x_2 \eta_2) 
 \\[0.2cm]
 \mathbf{C}_2\left(x_1,x_2, \eta_1, \eta_2 \right)  & = &
   \Big(  \log E\Big(\sum {X^2_{ex}}(t) \Big) -   \kappa_2(x_2 + \kappa_1 \eta_2) \Big)P(\mathbf{B} < x_2 \eta_2) P(\mathbf{B} > x_1 \eta_1) 
 \\[0.15cm] 
& + &  \Big(\log E\Big(\sum {X^2_{exm}}(t) \Big) -  \kappa_2(x_2 + \kappa_1 \eta_2) \Big)P(\mathbf{B} < x_1 \eta_1) P(\mathbf{B} < x_2 \eta_2).
 \end{eqnarray*}
We study the game theoretic problem in the budget constraint framework $ \mathit{ B}_i(x_i, \eta_i) := x_i + \kappa_1 \eta_i \le \bar{B}$ for CP-$i$ where $i = 1,2.$  We obtain the well-known solution concept to this game, \textit{Nash Equilibrium} (NE), using \textit{best response} method.

\begin{proposition} \label{Prop_2cp}
The best response of any CP say CP-1, $(x_1^*, \ \eta_1^*)$, against any strategy  $(x_2, \eta_2)$ of the other CP  \ satisfies \  $x_1^* + \kappa_1 \eta_1^* = \bar{B}.$
\\
\textbf{Proof:} Observe that the best response of CP-1 is computed by solving the optimization problem which is, essentially, similar to the problem $O2$. Now appealing to Proposition 1 of Part-I \cite{Ranbir2}, it follows. \eop
\end{proposition}

By the above proposition, each CP now has only one variable to choose and the other one $x_i$ is obtained through the equality constraint, i.e., $x_i =  \bar{B} -\kappa_1\eta_i; \ i = 1,2$. And this suffices to study the game with one controllable variable only.
  NE is an equilibrium strategy say $(\eta_1^*, \eta_2^*)$  deviating unilaterally from which neither of the CPs would  benefit.  Existence  and   uniqueness would be a topic of future research.
We numerically compute the  NE using  gradient and  best-response dynamics based algorithm. We vary the $m$ and study the NE in two scenarios: 1) when $m$ is directly proportional to $\mu_b$, and 2)  when $m$ is directly proportional to $\lambda$ in the figures below.

We see that in both the Figure \ref{Fig_2CPM0p1_w1p3} and Figure \ref{Fig_6To102CP_Lam0p1m_w1p2}, the more influential CP (i.e., CP-1) chooses $\eta_1^* = 1/w_1$ (maximum value of $\eta$) under the Nash strategy.
\vspace*{-1.5cm}
 \begin{figure}[H]
\begin{minipage}{9cm}
\begin{center}
\includegraphics[width = 7.5cm, height = 10cm]{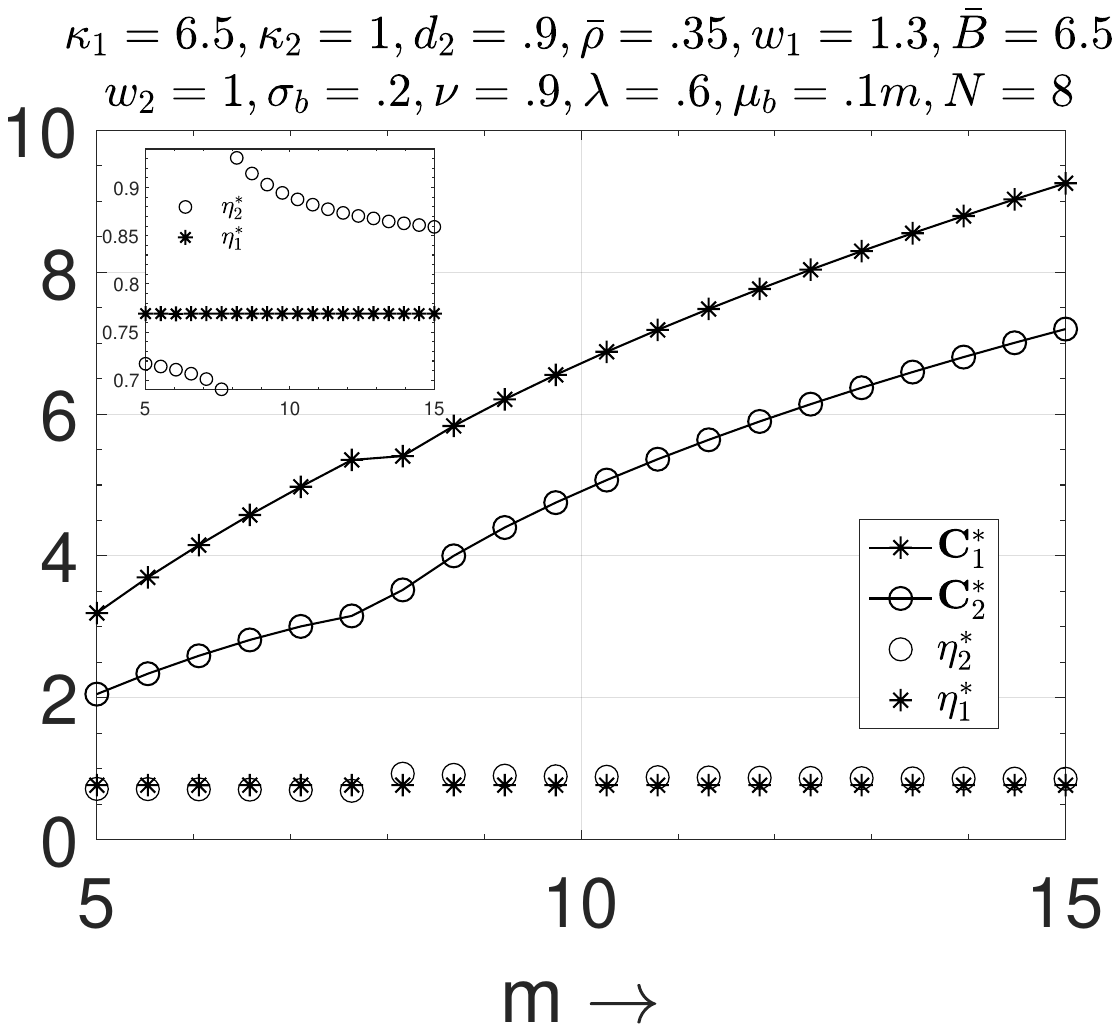}
\vspace*{-2.3cm}
\caption{When $m \propto \mu_b$ \label{Fig_2CPM0p1_w1p3}}
\end{center}
\end{minipage}
\hspace*{-1.8cm}
\begin{minipage}{9cm}
\begin{center}
\includegraphics[width = 7.5cm, height = 10cm]{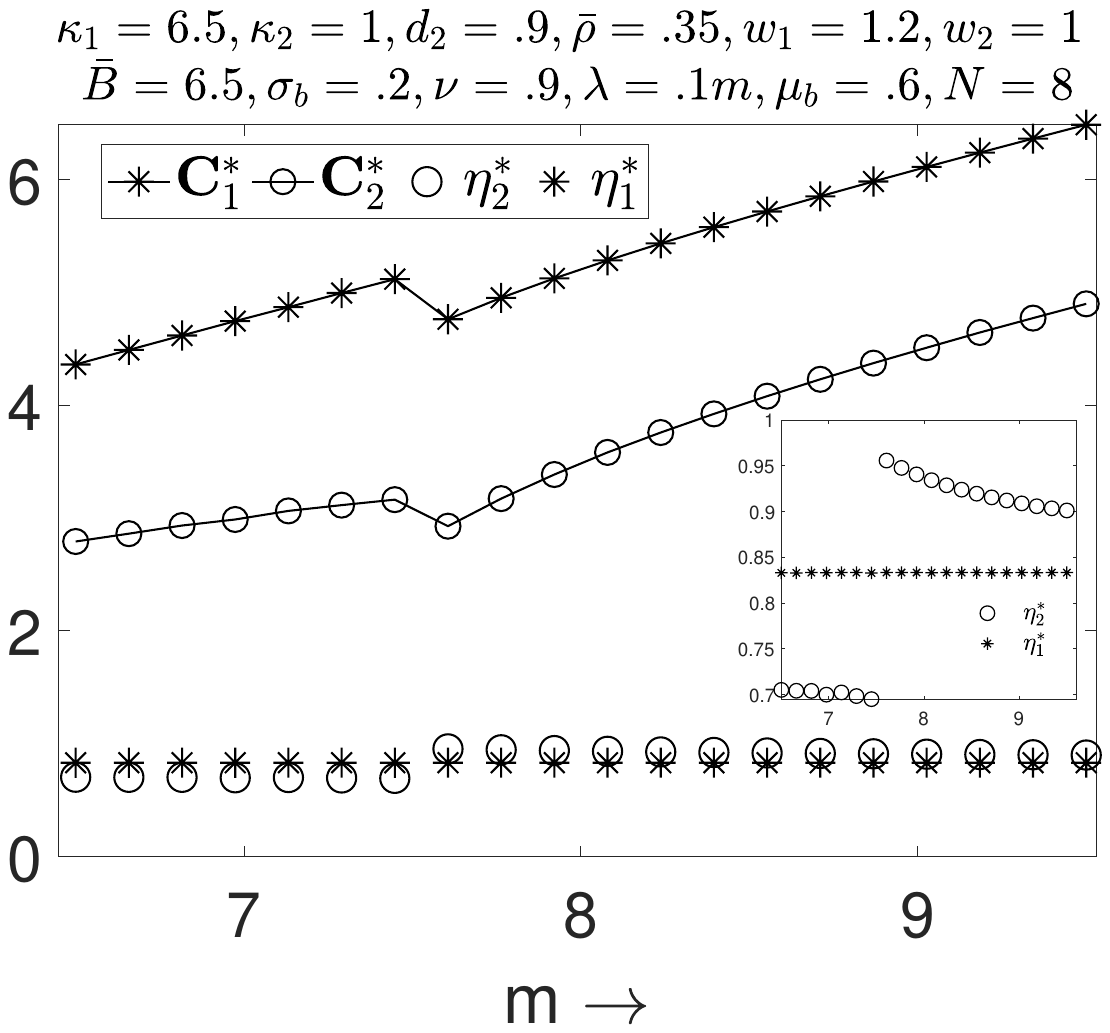} 
\vspace*{-2.5cm}\caption{When $m \propto \lambda$ \label{Fig_6To102CP_Lam0p1m_w1p2}}  
\end{center}
\end{minipage}
\end{figure}

 Whereas the Nash strategy of less influential one (CP-2) decreases initially, and then there is a jump discontinuity followed by a gradual decrease again. This is because the weaker CP has to invest more money in winning the auction initially ($m$ is small) as it gives higher returns compared to that in $\eta$, and consequently the allocation on $\eta$ shrinks (in both figures). Beyond a threshold $m$, investing in $\eta$ earns more revenue. And the gradual decrease is due to the following: in Figure \ref{Fig_2CPM0p1_w1p3}, winning auction gets difficult as $m$ is directly proportional to $\mu_b$ and hence CP needs to invest more in auction;  whereas in Figure \ref{Fig_6To102CP_Lam0p1m_w1p2}, the growth rate $\alpha$ decreases as $\lambda$ increases investing in $x$ yields higher returns compared to that in $\eta.$ 
 
\begin{figure}[h]
\begin{center}
\includegraphics[width = 10.5cm, height = 12cm]{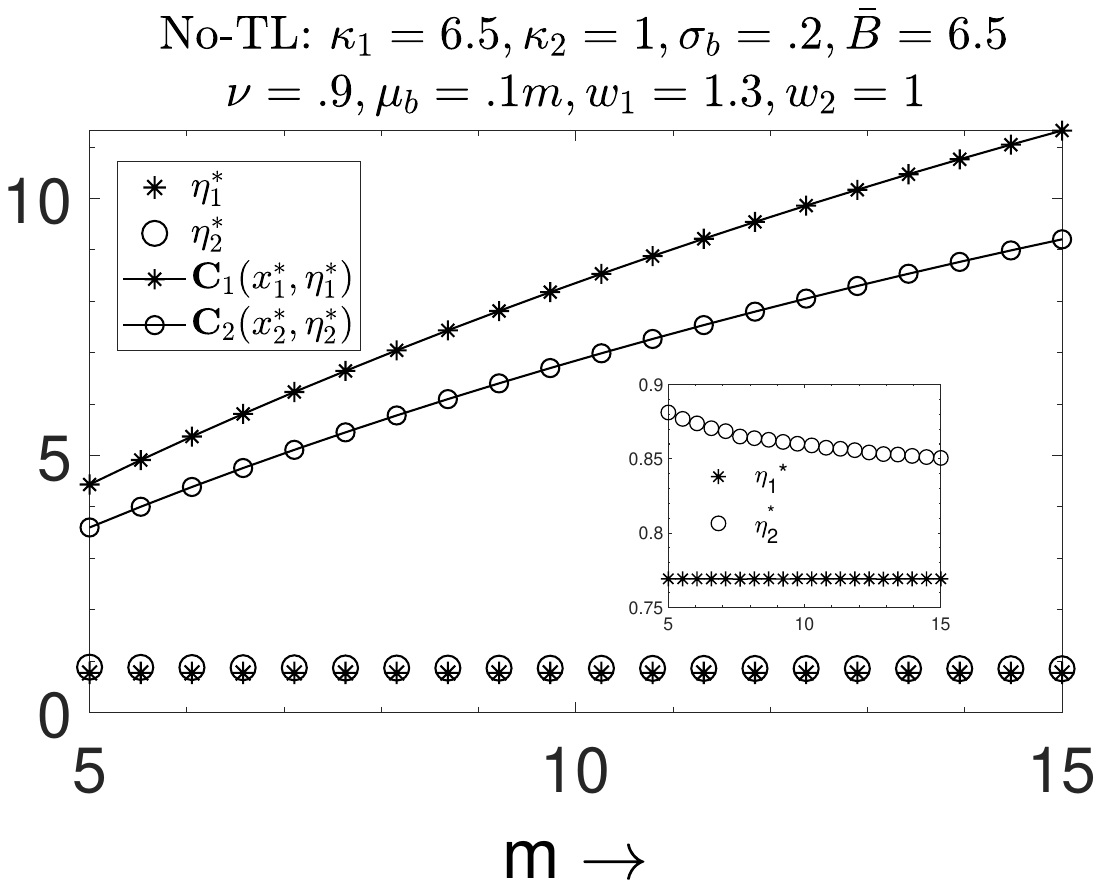} 
\vspace*{-2.5cm}\caption{No-TL case  \label{NTL2}}  
\end{center}
\end{figure}

In No-TL case,  we see in Figure \ref{NTL2} that the Nash strategies  show monotonous behavior as the network activity increase ($m$). This pattern is quite different from that seen in Figure \ref{Fig_2CPM0p1_w1p3} in the Nash strategy of CP-2, $\eta^*_2$. 

In all, we see that the performance measures and the optimizers  are drastically different in the study with and without considering the TL structure. Further, as the network activity increases, one anticipates that a good quality content can easily get viral. However, this is not true because of the shifting effect. Recall that as $m$ increases, contents get pushed down rapidly by the arrival of new posts, and hence the content gets missed more often before a user visits its TL.  These important aspects are missed when the TL structure is ignored.

\section*{Conclusions}

We modeled the propagation of competing posts by multi-type branching process. As the underlying branching process is decomposable, the entire analysis is   different from that in the Part-I. We found that the dichotomy no longer holds, i.e., one of the competing posts may get viral while the other gets extinct;  different types of populations can have different growth rates. We obtained various performance measures, using the previous results, specific to individual CP such as CP-wise extinction probabilities, the expected number of shares to a CP's post in the viral and non-viral scenario, etc. We conjectured using partial theoretical arguments that the expected number of shares corresponding to one CP  grow exponentially fast with time (if viral) in the presence of the competing posts. We verified the same numerically.  We found that the virality chances of a post are greatly influenced by the competing post propagation. We then formulated a non-cooperative game between competing CPs using the CP-wise performance measures and studied the relevant Nash equilibria. Again, we found that the study without considering the TL structure cannot capture accurately the competition induced due to the propagation of competing posts.
 
More importantly, we also observe that without TL effects, one cannot capture  some interesting paradigm shifts/phase transitions in certain behavioral patterns. For example, as the network becomes more active, one anticipates  that it is more beneficial to engage in the network. The studies which do not include the effects of TL often leads to this erroneous conclusion; and argue that the virality chances increase monotonically as the mean number of friends increases ($m$). We demonstrated that the virality chances does not increase monotonically with the number of friends. After a certain value of $m$, it actually decreases for some intermittently active networks (medium $m$ values).   To be more specific, for some range of parameters,  less active networks are preferable to more active networks.

\section*{Appendix A}
 \textbf{Proof of Theorem \ref{Thm_mixed} :}  The proof of \textit{part-i} is as follows. The generator matrix $A_{mx}$ is 
 \begin{eqnarray}
\label{GenMatrixTwoCP}
 \left[ \begin{array}{cccccc}
 z'_1r_1 -1 &  z_1 r_1 &  \theta +z'_2 r_1 &  \dots  & z'_{N-1}r_1 & z_{N-1}r_1 \\
  z_1 r_1 &   z'_1r_1 -1  & z_2 r_1  &   \dots  & z_{N-1}r_1  & z'_{N-1}r_1 \\
   z'_1r_2 &  z_1 r_2 &  z'_2r_2- 1&   \dots  & z'_{N-1}r_2  & z_{N-1} r_2 \\
  z_1 r_2  &  z'_1r_2 & z_2 r_2  & \dots  & z_{N-1}r_2  & z'_{N-1}r_2 \\
  \vdots & \vdots & \vdots   & \ddots     & \vdots  & \vdots  \\
   z'_1 r_{N-2}  &  z_1 r_{N-2}  & z'_2r_{N-2} &   \dots  & \theta + z'_{N-1} r_{N-2} & z_{N-1} r_{N-2}  \\
  z_1 r_{N-2}  &  z'_1r_{N-2}  & z_2 r_{N-2} & \dots  & z_{N-1}r_{N-2}  & \theta + z'_{N-1}r_{N-2}\\
z'_1 r_{N-1}  &  z_1 r_{N-1}  & z'_1 &   \cdots  & z'_{N-1} r_{N-1}-1 & z_{N-1} r_{N-1} \\
  z_1 r_{N-1}  &  z'_1r_{N-1} & z_2 r_{N-1}  &  \dots  & z_{N-1}r_{N-1}  & z'_{N-1} r_{N-1} -1
 \end{array} 
 \right]. 
\end{eqnarray}
First, we prove that  $e^{A_{mx}}$ is positive regular for any $0< \theta, p < 1$. 
As in the case of Lemma 1 of Part-I \cite{Ranbir2}, we prove the result for a special  case with  $z_l = z'_l = 0 \ \forall l > 1$ and $z_1 > 0$, $z'>0$. The  result again follows for the general case because all the terms involved are non-negative.  For this special case, the matrix $A_{mx}+I$  has the following  form with all $z_l r_k$ or $z_l'r_k$ terms being strictly positive because 
 $0 < p < 1$:
\begin{eqnarray*}
 A_{mx}+I =  \left[ \begin{array}{ccccccccc}
z_1' r_1 &  z_1 r_1 &  \theta & 0 & \dots  & 0 &0 & 0 &0 \\
  z_1 r_1 &   z_1' r_1  & 0  & \theta &  \dots  & 0  & 0 & 0  & 0  \\
 z_1' r_2 &  z_1 r_2 &  0 & 0 &  \dots  & 0 & 0 & 0  & 0  \\
  z_1 r_2  &  z_1' r_2 & 0 & 0 & \dots  & 0 & 0 & 0 & 0  \\
  \vdots & \vdots & \vdots   & \vdots & \ddots     & \vdots  & \vdots & \vdots  & \vdots  \\
  z_1' r_{N-2}  &  z_1 r_{N-2}  & 0 & 0&  \dots & 0  & 0  & \theta  & 0  \\
  z_1 r_{N-2}  &   z_1' r_{N-2}  & 0& 0& \dots & 0  & 0  & 0  & \theta \\
 z_1' r_{N-1}   &  z_1 r_{N-1}  & 0  & 0&  \cdots & 0  & 0   & 0 & 0\\
  z_1 r_{N-1}  &   z_1' r_{N-1}  & 0  & 0& \dots  & 0  & 0 & 0  & 0
 \end{array} 
 \right];  
\end{eqnarray*}
Positive regularity of a matrix is determined by the existence of all positive terms in some power of the given matrix. Since the matrix $A_{mx}+I$ has only non-negative entries, it is sufficient to check zero, non-zero structure (the location of zero and non zero terms in the given matrix and not the exact values) of the resulting powers of the matrices 
 $(A_{mx}+I)^n$.
The matrix $(A_{mx}+I)$  is exactly similar in zero non-zero structure  as  the second power $ A_1^2$ given in Part-I \cite{Ranbir2}. 
Thus, positive regularity follows in exactly similar lines.

 {\bf  Proof of parts (ii)-(iii): } We follow exactly the same procedure as in the proof of parts (ii)-(iii) of Lemma 1 of Part-I \cite{Ranbir2}. We mention only the differences with respect to that proof. 

Let $\textbf{u}_{mx}  = \{ u_{mx,1}, u_{mx,2}, \cdots,  u_{mx,2N-3},  u_{mx,2N-2}\}$ be the left eigenvector of $A_{mx}$, corresponding to largest eigenvalue $\alpha_{mx}$, both of  which  exist because of positive regularity given by part (i). 
 
On solving $\textbf{u}_{mx} A_{mx} = \alpha_{mx} \textbf{u}_{mx}$ as before, we have the following system of equations with $\sigma_{mx} = \alpha_{mx}/(\lambda + \nu) + 1$:
 \begin{flalign}
 z_1 \textbf{r}.\textbf{u}_{mx,e} +  z'_1 \textbf{r}.\textbf{u}_{mx,o} & = & \sigma_{mx}u_{mx,1}, \hspace{0.2cm} z_l \textbf{r}.\textbf{u}_{mx,e} + z'_l \textbf{r}.\textbf{u}_{mx,o} + \theta u_{mx,2l-3} \ = \ \sigma_{mx} u_{mx,2l-1}; \hspace{0.1cm} \forall l \geq 2 \nonumber \\
  z'_1 \textbf{r}.\textbf{u}_{mx,e} +   z_1 \textbf{r}.\textbf{u}_{mx,o}  & = &  \sigma_{mx} u_{mx,2}, \hspace{0.2cm} z'_l \textbf{r}.\textbf{u}_{mx,e} + z_l\textbf{r}.\textbf{u}_{mx,o} +  \theta u_{mx,2l-2}  \ = \ \sigma_{mx} u_{mx,2l}; \hspace{0.1cm} \forall l \geq 2. 
  \label{Eqn_zz'} 
 \end{flalign}
 where $\textbf{r}.\textbf{u}_{mx,o}: = \sum_{i = 1}^{N-1} r_{i}u_{mx,2i-1} $,   $ \textbf{r}.\textbf{u}_{mx,e} := \sum_{i = 1}^{N-1} r_{i} u_{mx,2i} $ and  $u_{mx,-1}(u_{mx,-2}) :=0.$ 
 Now, we write the expression of $u_{mx,2l-1}$ and $u_{mx,2l}$ in terms of $u_{mx,1}$ and $u_{mx,2}$ receptively as done before in the single CP case.  After simplifying equations (\ref{Eqn_zz'}) , we have the following for any $2 \le l \leq N-1$
 \begin{eqnarray*}
 u_{mx,2l-1}  =  \sum_{i = 0}^{l-1} \frac{\bar{\rho}_{l-i}}{\bar{\rho}_1} \left(\frac{\theta}{\sigma_{mx}}\right)^{i} u_{mx,1}; u_{mx,2l}  =  \sum_{i = 0}^{l-1} \frac{\bar{\rho}_{l-i}}{\bar{\rho}_1} \left(\frac{\theta}{\sigma_{mx}}\right)^{i} u_{mx,2}. 
 \end{eqnarray*}
 Following the same procedure, we obtain the relation among various components of   right  eigenvector  $\textbf{v}_{mx}$. Thus, we have  ($\ \forall l = 1,\cdots N-2$)
\begin{eqnarray*}
\label{Eqn_vl_mx_pattern}
v_{mx,2l-1} = \sum_{i = 0}^{N-1-l} \frac{r_{l+i}}{r_{N-1}} \left(\frac{\theta}{\sigma_{mx}}\right)^{i} v_{mx,2N-3}; \ \ v_{mx,2l-2} = \sum_{i = 0}^{N-1-l} \frac{r_{l+i}}{r_{N-1}} \left(\frac{\theta}{\sigma_{mx}}\right)^{i} v_{mx,2N-2}. 
\end{eqnarray*}

Recall $c_{mx}  =\delta(1-\theta)m\eta_1\eta_2$, the above equations can be rewritten as
 \begin{eqnarray*}
 c_{mx}\bar{\rho}_1 \left(p \textbf{r}.\textbf{u}_{mx,e} +  (1-p) \textbf{r}.\textbf{u}_{mx,o}\right) \ = \  \sigma_{mx} u_{mx,1}, \hspace{0.2cm}  c_{mx}\bar{\rho}_1 \left((1-p)\textbf{r}.\textbf{u}_{mx,e} +  p \textbf{r}.\textbf{u}_{mx,o}\right) = \sigma_{mx} u_{mx,2}
 \end{eqnarray*} 
 and \vspace*{-0.2cm}
 \begin{eqnarray*}
  c_{mx}\bar{\rho}_l \left(p \textbf{r}.\textbf{u}_{mx,e} +  (1-p) \textbf{r}.\textbf{u}_{mx,o}\right)+ \theta u_{mx,2l-3} & = & \sigma_{mx}  u_{mx,2l-1} \hspace{0.1cm} \forall l \geq 2
 \\
  c_{mx}\bar{\rho}_l \left((1-p) \textbf{r}.\textbf{u}_{mx,e} +  p \textbf{r}.\textbf{u}_{mx,o}\right)+ \theta u_{mx,2l-2} & = & \sigma_{mx}  u_{mx,2l-2} \hspace{0.1cm} \forall l \geq 2. 
 \end{eqnarray*}
On multiplying with $r_i$ and adding all even and odd term equations separately, we obtain 
\begin{eqnarray*}
 c_{mx}\textbf{r.}\bm{\bar{\rho}} \left(p \textbf{r}.\textbf{u}_{mx,e} +  (1-p) \textbf{r}.\textbf{u}_{mx,o}\right)  + \theta \sum_{i = 1}^{N-1} r_{i} u_{mx,2i-3} = \sigma_{mx} \textbf{r}.\textbf{u}_{mx,o} \ \ \text{; for odd terms,}
 \\
 c_{mx}\textbf{r.}\bm{\bar{\rho}} \left((1-p)\textbf{r}.\textbf{u}_{mx,e} +  p \textbf{r}.\textbf{u}_{mx,o}\right)  + \theta \sum_{i = 1}^{N-1} r_{i} u_{mx,2i-2} = \sigma_{mx} \textbf{r}.\textbf{u}_{mx,e}   \ \ \text{; for even terms.} 
\end{eqnarray*}
On adding the above equations, we obtain the following linear equation
{\footnotesize{
\begin{eqnarray}
P(\sigma_{mx}) =   c_{mx}\textbf{r.}\bm{\bar{\rho}}  \left(\textbf{r.u}_{mx,o} +  \textbf{r.u}_{mx, e}\right) + \theta \sum_{i = 1}^{N-1} r_i \left(u_{mx,2i-3} + u_{mx,2i-2}\right)- \sigma_{mx} \left(\textbf{r.u}_{mx,o} +  \textbf{r.u}_{mx,e}\right), \label{TwoCPEqaution}
\end{eqnarray}}}
and $\sigma_{mx} = (\alpha_mx  + \lambda + \nu)/(\lambda + \nu)$ would be the only zero of it. 
Now $P\left( c_{mx}\textbf{r.}\bm{\bar{\rho}} \right) > 0$ and $P\left( c_{mx}\textbf{r.}\bm{\bar{\rho}} + \theta \right) < 0$ (again using monotonicity of the reading probabilities). And due to similar reasons as in the singe CP case, the largest eigenvalue lies in the  interval  $\alpha_{mx} \in \left(c_{mx}\textbf{r.}\bm{\bar{\rho}}-1 , \ \ c_{mx}\textbf{r.}\bm{\bar{\rho}}  + \theta -1\right)(\lambda+\nu).$
\\
Let us assume $r_l = d_1d_2^l$, once  $u_{mx,2N-2} + u_{mx,2N-3}$ are both bounded for any $N$. 
In what follows, the  root of equation (\ref{TwoCPEqaution}) for this special  case 
\begin{eqnarray*}
\sigma_{mx} = c_{mx}\textbf{r.}\bm{\bar{\rho}} + \theta d_2 \dfrac{\sum_{i = 1}^{N-2} r_i \left( u_{mx,2i-1} + u_{mx,2i} \right)}{\left(\textbf{r.u}_{mx,o} +  \textbf{r.u}_{mx,e}\right)  } 
 c_{mx}\textbf{r.}\bm{\bar{\rho}} + 
  \theta d_2 \left(1-\frac{r_{N-1} \left( u_{mx,2N-2} + u_{mx,2N-3}\right)}{\textbf{r.u}_{mx,o} +  \textbf{r.u}_{e}} \right),
  \end{eqnarray*}
  converges to the following because $r_N = d_1 d_2^N \to 0$,
\begin{eqnarray*} 
\sigma_{mx}   \to c_{mx}\textbf{r.}\bm{\bar{\rho}}  + \theta d_2 \ \ \text{as \ } N \to \infty .
\end{eqnarray*}  

Thus, as the number of TL levels increases the largest eigenvalue, $\alpha_{mx}$ of matrix $A_{mx}$ converges to $\left(c_{mx}\textbf{r.}\bm{\bar{\rho}} + \theta d_2 -1 \right)(\lambda + \nu) $.

\textbf{Part-iv:} Observe that  $\{  \textbf{X}_{mx}(t) \}$ evolves according to non-decomposable BP when the process starts with a mixed-type TL, as in the single CP case. Now appealing to Theorem 1 of Part-I \cite{Ranbir2}, the proof follows.
\eop
\\

\textbf{Proof of Lemma \ref{Lemma_mixed_unique_extn}:}
\textit{Existence:}  We consider ${\bf q}_{ex}^1$ a constant vector, as explained below. We have a continuous mapping from $[0,1]^{2N-2}$  into $[0,1]^{2N-2}$ i.e. over compact set. By \textit{Brouwer's fixed point theorem} there exists a solution to the given system of equations. 

\textit{Uniqueness:} The exclusive CP1 types evolve on their own, and by Lemma 2 of Part-I \cite{Ranbir2}, we have a unique solution to the relevant fixed point equations in unit cube $[0,1]^N$ which provide the extinction probabilities of CP1 population when started with one of its exclusive-types. That is, we have a unique  
$\textbf{q}^1_{ex}=\{q^l_{l,0} \}_l$ which represents the extinction probabilities for any given set of system parameters.  We treat them as constants while studying the fixed point equations of the other equations that provide the extinction probabilities when started with a mixed population, $(\textbf{q}_{mx1}\textbf{q}_{mx2})$.  One can rewrite the fixed point equations corresponding to this set of the extinction probabilities   as below for any $l < N$, after suitable simplification:  
{\footnotesize{
\begin{eqnarray}
q_{l,l+1}^1 & = &  K_{1l} (pg_{mx1} +(1-p) g_{mx2} )  \bpgftwo \left(\textbf{q}^1_{ex},\eta_1(1-\eta_2) \right) +  K_{2l}(1-\delta) \bpgftwo\left(\textbf{q}^1_{ex},\eta_1 \right)+ 
K_{3l}  + \theta^{N-l}q^1_{N,0}   \nonumber
\\
q_{l+1,l}^1  & = & K_{1l} \left( (1-p)g_{mx1} + pg_{mx2} \right)  \bpgftwo \left(\textbf{q}^1_{ex},\eta_1(1-\eta_2) \right)  + K_{2l}(1-\delta) +K_{3l} +
 \theta^{N-l}   \label{Eqn_simple_mix_qs}
\end{eqnarray}}}
where
\begin{eqnarray}
K_{2l}  & = &  (1-\theta) \sum_{i=0}^{N-l-1}\theta^i r_{l+i}; \hspace{0.2cm} 
K_{3l}  = (1-\theta) \sum_{i=0}^{N-l-1}\theta^i \left(1-r_{l+i} \right) \nonumber
 \\
 K_{1l} & = & K_{2l} \delta  \nonumber
 \\
 g_{mx1}  & = & \bpgftwo\left(\textbf{q}_{mx1},\eta_1 \eta_2 \right); \hspace{0.3cm} g_{mx2} =  \bpgftwo\left(\textbf{q}_{mx2},\eta_1\eta_2 \right)\nonumber.
 \end{eqnarray}
Consider the following weighted sum over $l$, of terms $f (q_{l,l+1}^1,  \    \eta_1\eta_2,  \   \beta)$ and $f(q_{l+1,l}^1,  \ \  \eta_1\eta_2,  \  \beta)$
$$
\sum_{l=1}^{N-1}  {\bar \rho}_l f (q_{l,l+1}^1) \, \ \ \mbox{ and }  \,\sum_{l=1}^{N-1}  {\bar \rho}_l   f(q_{l+1,l}^1)   ,
$$  and note that these precisely equal $g_{mx1}$ and $g_{mx2}$ respectively. 
Thus using  the right hand side (RHS) of equation (\ref{Eqn_simple_mix_qs}), 
we have the following two dimensional   equation,  $\Psi = (\Psi_1, \Psi_2)$, whose fixed point provides $(g_{mx1}, g_{mx2})$: 
{\footnotesize{
\begin{eqnarray*}
\Psi_1 (g_{1}, g_2 ) =   \sum_{i=1}^{N-1}f \left( \theta^{N-i}q^1_{N,0} +  K_{1i} \big (pg_{1} +(1-p) g_{2} \big )  \bpgftwo \left(\textbf{q}^1_{ex},\eta_1(1-\eta_2) \right) 
+ K_{2i} (1-\delta) \bpgftwo\left(\textbf{q}^1_{ex},\eta_1 \right) +K_{3i} , \eta_1\eta_2,  \beta \right) {\bar \rho}_i \nonumber 
\\
\Psi_2 (g_{1}, g_2 ) =  \sum_{i=1}^{N-1} f \left( \theta^{N-i} +  K_{1i}\big ( (1-p)g_{1} + pg_{2} \big) \bpgftwo \left(\textbf{q}^1_{ex},\eta_1(1-\eta_2) \right)  + K_{2i} (1-\delta)   + K_{3i}, \ \eta_1\eta_2,\ \beta  \right){\bar \rho}_i. \nonumber 
\end{eqnarray*}}}
 It is easy verify for any $l$ that 
  \[
K_{1l}+ K_{2l} (1-\delta)   + K_{3l} \ =\  K_{2l}    + K_{3l}  \   = \   (1- \theta) \sum_{i=0}^{N-l-1} \theta^i   =  \Big(1-\theta^{N-l}\Big);
\] 
 and hence that 
 $$
 \theta^{N-l} + K_{1l} +  K_{2l} (1-\delta)   + K_{3l}  = 1.
 $$
 Thus for any $\textbf{q}^1_{ex} \le {\bf 1}$ we\footnote{Here $\le$ represents the usual partial order between two Euclidean vectors, i.e., ${\bf a} < {\bf b}$
 if and only if $a_i < b_i$ for all $i$ and ${\bf a} \le {\bf b}$  if $a_i \le b_i$ for all $i$ .}   have:
  
\begin{eqnarray}
\label{eqn_Ks_one}
\theta^{N-l} +  & \hspace{-15mm} K_{1l} \bpgftwo \left(\textbf{q}^1_{ex},\eta_1(1-\eta_2) \right)  + K_{2l} (1-\delta)     + K_{3l}   \le 1 \\
\theta^{N-l}q^1_{N,0} +&  K_{1l}   \bpgftwo \left(\textbf{q}^1_{ex},\eta_1(1-\eta_2) \right) + K_{2l} (1-\delta) \bpgftwo\left(\textbf{q}^1_{ex},\eta_1 \right) +K_{3l}  \le 1\nonumber .
\end{eqnarray}

 {\bf Case 1  When $\textbf{q}^1_{ex} < {\bf 1}$:}
 When $\textbf{q}^1_{ex} < 1$,  $\bpgftwo\left(\textbf{q}^1_{ex},\eta_1 \right) < 1$ as well as  $ \bpgftwo \left(\textbf{q}^1_{ex},\eta_1(1-\eta_2) \right)  < 1$ and so we have strict inequality 
 in (\ref{eqn_Ks_one})  and thus 
  $\Psi_j (1,1) < 1$ for each $j.$ 
Consider  $j=1$ without loss of generality.
Thus, $\Psi_1 (1,g_2) < 1$ for any $g_2  \le 1$. 
Consider the one-variable 
function $g \to \Psi_1 (g, g_2)$, represented by
$$
\Psi_1^{g_2} (g) := \Psi_1 (g, g_2),
$$ for any fixed $g_2$, which is clearly a continuous and monotone function.   Let $id(g) := g$ represent the identity function.
 From the definition of $\Psi$, clearly  $\Psi_1^{g_2} (0) > 0$  for any $g_2$.   Hence $\Psi_1^{g_2} (0) - id(0) > 0$ while $\Psi_1^{g_2} (1) - id(1) < 0$.
Thus, by intermediate value theorem as applied to the (continuous) function $\Psi_1^{g_2} (\centerdot) - id(\centerdot)$, 
 there exists at least one point at which it crosses the 45-degree line, the straight line  through origin (0,0)  and (1,1). 
 Note that the intersection points of this 45-degree line and a function are precisely the fixed points of that function.

 It is easy to verify that the derivative of the function $\Psi_1^{g_2}$ (partial derivative of $\Psi_1$ with respect to the second variable)  is positive.  Thus $
\Psi_1^{g_2} (\centerdot) $ 
 for any fixed $g_2$ is continuous increasing strict convex function. 
 If $\Psi_1^{g_2} (\centerdot)$ function were to cross 45-degree line more than once  before reaching $\Psi_1^{g_2} (1) < 1$ at 1, then it would  have to cross the 45-degree line three times (recall  $\Psi_1^{g_2} (0)>0$).  However, this is not possible because 
 any strict convex real function crosses any  straight line at maximum twice.   Thus, there exists exactly one point   in interval $[0,1]$
at which $\Psi_1^{g_2} (\centerdot)$ crosses 45-degree line, which would be its  unique fixed point.  

Thus, for any   $g_2$ there exists a unique fixed point of
the mapping $\Psi_1^{g_2} (\centerdot)$ in the interval $[0,1]$ and call the unique fixed point  as ${\bf g}^*(g_2)$.  It is easy to verify that this fixed point is minimizer of the following objective function parametrized by $g_2$:
$$
\min_{g \in [0,1] } \Phi (g, g_2) \mbox{ with }  \Phi (g, g_2)  := \left (  \Psi_1 (g, g_2) - g \right )^2.
$$ The function $ \Phi$ is jointly continuous, convex in $(g, g_2)$ and the domain of optimization is same for all $g_2$. Further, for each $g_2$ by previous arguments there exists unique optimizer in $[0,1]$. Thus, by \cite[Maximum Theorem for Convex Functions]{FirstOpt}, the fixed point function ${\bf g}^*(.) $ is continuous, 
and convex function.

We now obtain the overall (two dimensional) fixed point via the solution of the following one-dimensional fixed point equation. 
$$
{\bm \Gamma}(g) := \Psi_2 ({\bf g}^*(g) ,  g).
$$ 
$$
f \left( \theta^{N-i} +  K_{1l}\big ( (1-p)g_{1} + pg_{2} \big) \bpgftwo \left(\textbf{q}^1_{ex},\eta_1(1-\eta_2) \right)  + K_{2l} (1-\delta)   + K_{3l}  , \ \ \ \eta_1\eta_2,  \ \ \  \beta  \right)
$$
Let $K_{4l} : = \theta^{N-l} + K_{2l} (1-\delta)   + K_{3l}$ and $K_{5l} :=  K_{1l} \bpgftwo \left(\textbf{q}^1_{ex},\eta_1(1-\eta_2) \right) $. With these 
definitions:
$$
{\bm \Gamma}(g) = \sum_{l=1}^{N-1} f \left( K_{4l}  +  K_{5l}  \big ( (1-p) {\bf g}^*(g) + pg \big) , \ \eta_1\eta_2,\beta  \right){\bar \rho}_l 
$$ 

Consider any $0\le \gamma, g, g' \le 1 $ and  by convexity of $ {\bf g}^* $ and  monotonicity of $\Psi_2$ we have 
{\footnotesize{
\begin{eqnarray*}
{\bm \Gamma}(\gamma g + (1-\gamma) g')  & = & \sum_{l=1}^{N-1} f \left( \ K_{4l}  +  K_{5l} \left( (1-p){\bf g}^*(\gamma g + (1-\gamma) g') + p \Big [\gamma g + (1-\gamma) g'\Big ] \right) ,\eta_1\eta_2,\beta  \right){\bar \rho}_l
 \\
  & \le & \sum_{l=1}^{N-1} f \left(  K_{4l}  +  K_{5l} \left( (1-p) \Big [\gamma {\bf g}^*(g) + (1-\gamma) {\bf g}^*(g')\Big ] + p \Big [\gamma g + (1-\gamma) g' \Big ] \right)  ,\eta_1\eta_2,\beta  \right){\bar \rho}_l 
  \\
  &  &\hspace*{-4cm} = \ \sum_{l=1}^{N-1} f \left(  K_{4l}  +  K_{5l} \left(\gamma \Big [ (1-p) {\bf g}^*(g) + p  g\Big ] + (1-\gamma) \Big [ (1-p) {\bf g}^*(g')+ p   g' \Big ]\right),  \eta_1\eta_2,\beta  \right){\bar \rho}_l; \ \ \ \mbox{using convexity of } f 
   \\
    &  & \hspace*{-4cm}  \le \ \ 
    \sum_{l=1}^{N-1}\Bigg ( \gamma f \left(  K_{4l}  +  K_{5l}  \Big [  (1-p) {\bf g}^*(g) + p  g \Big ] , \ \eta_1\eta_2,\beta  \right) +  (1-\gamma)  f \left( K_{4l}  +  K_{5l}  \Big [ (1-p) {\bf g}^*(g')+ p   g' \Big ]  , \ \eta_1\eta_2,\beta  \right) \Bigg )
       \\
      &=&  
      \gamma {\bm \Gamma}(  g)  + (1-\gamma) {\bm \Gamma}(  g').
\end{eqnarray*}}}
This shows that   $ {\bm \Gamma}$ is convex, further we have ${\bm \Gamma} (1) < 1$ and ${\bm \Gamma} (0) > 0$. 
Note here that ${\bf g}^*(0) > 0$ because $\Psi_1(0,0) >0$.
Thus, using similar arguments as before we establish the existence of unique fixed point $g_2^*$ for function  ${\bm \Gamma}$.  Therefore, $({\bf g}^*(g_2^*), g_2^*)$ represents the unique fixed point, in unit cube $[0, 1]^2$, of the two dimensional function $\Psi$.  This establishes the existence and uniqueness of extinction probabilities $(g_{12}, g_{21})$.

The uniqueness of other extinction probabilities is now direct from equation (\ref{Eqn_simple_mix_qs}).

 {\bf Case 2  When $\textbf{q}^1_{ex} = 1$:} 
Consider that we start with one of the following three TLs: one exclusive  type $(l, 0)$, one  mixed  type $(l, l+1)$ or one mixed-type $(l+1, l)$.  Consider  the scenario in which the  CP-1 population gets extinct at the first transition  epoch itself,  when started with one $(l, 0)$ type. This can happen  if one  of the following two  events occur: a)  the TL does not view post-${\bf P}$ (w.p. $r_l$);  or b) the TL shares  post-${\bf P}$  to none (0) of its friends. In either of the two events the CP-1 population gets extinct even when started with mixed TLs  $(l, l+1)$ or $(l+1, l)$. Thus, the event of extinction at first transition epoch starting with one $(l, 0)$ TL implies extinction at first transition epoch when started with either one  $(l, l+1)$ TL or one $(l+1, l)$ TL.  Say the number of shares at first transition epoch were non-zero and say they equal $x_i$ of type $(i,0)$ for each $i$  when started with one $(l,0)$ type TL. This proof is given under extra assumption that $\rho_N = 0$ and that ${\bar \rho}_i = \rho_i$.  {\it We assume the following is the scenario under assumption.}  When we start with mixed-type  $(l, l+1)$ (or $(l+1, 1)$ type respectively),    post-${\bf P}$ is shared with $\sum_i x_i$ number of Friends as   when started with  exclusive $(l,0)$ type.  Out of these, some are now converted to mixed TLs because the parent TL also shares CP-2 post.  And a converted  type $(i, 0)$ offspring becomes $(i, i+1)$ offspring w.p. $p$  (w.p. $(1-p)$  respectively) and 
 $(i+1, i)$ w.p. $(1-p)$  (w.p. $(1-p)$ respectively). 
When started with mixed-type $(l+1, l)$ it is possible that some out of $\sum x_i$ shares of CP1 post are discarded (w.p. $\delta$) because the TL would have viewed the post-${\bf Q}$ first and would be discouraged to view post-${\bf P}$.  
  Thus,  in either case, with or without extinction at  first transition epoch,  the resulting events  are inclined towards survival with bigger probability  when started with one exclusive $(l,0)$ type than when started with either of the mixed-type TLs.
Basically, the aforementioned 
  arguments can be applied recursively to arrive at  this conclusion, and hence   the probabilities of extinctions satisfy the following inequalities:
$$
q_{l, 0}^1  \le q_{l, l+1}^1  \,\,  \mbox{ and }  \,\, q_{l, 0}^1  \le q_{ l+1, l}^1   \mbox{ for any }   l < N-1.
$$
Further,   
with $\textbf{q}^1_{ex} = 1$, it easy to verify that  $\Psi_i(1, 1)  =1$ for $i=1$ as well as $2$.  
Thus, we have unique extinction probabilities, ${\bf q}_{mx1} = {\bf 1}$ and  ${\bf q}_{mx2} = {\bf 1}$.  \eop
\\
\\

 \section*{Appendix B: Evolution of  Exclusive net progeny   }

The evolution of the size of the population of exclusive/mixed class when initiated with its own class particle(s) is obtained using the well-known theory of non-decomposable BPs, in particular   Lemma 3 of  Part-I (\cite{Ranbir2})  provides the time evolution of expected net progeny; 
it is easy to observe that the `number of shares' is  the net progeny. 
  According to \cite[Lemma 3]{Ranbir2} the  net progeny  of an ${\cal E}_x$ class  till time $t$ when initiated with type-$l$ particle of  its own class, i.e., with $l \in {\cal E}_x$,  is represented by 
 $y^e_l(t)$  and is provided in the vector form.  
 
 In \cite{Ranbir2}  the study is about the net progeny when one started with one $l$-particle, for any  $l \in {\cal E}_x$ and  after setting $y_l^e(0) = 1$. We require a small change to facilitate study   of net progeny when  started in mixed class; we require that  $y^e_l(0) = 0$, i.e., net progeny at time 0 is set to 0; in  other words, the initial particle is not counted as an offspring, but its offsprings, offsprings of offsprings so on to  form the progeny.  We consider this study, using similar fixed point equations as in section \ref{sec_type_change}. 
 
%

Consider only exclusive types, evolving on their own;  when you start in ${\cal E}_x$  class, the particles produce offsprings of only  ${\cal E}_x$  class.
Conditioning on the events of first transition,  $y_l (t)$ (for any $l \in {\cal E}_x$) satisfies the following fixed point equation 
\begin{eqnarray}
y^e_{l} (t)  &=&  \theta \int_0^t    \sum_{k \in {\cal E}_x   } a_{l, k}  y^e_k (t-\tau)   (\lambda+\nu) e^{- (\lambda+\nu) \tau } d\tau  \nonumber  \\
&& 
+  (1-\theta)  \int_{0}^t   \sum_{k  \in  {\cal E}_x }  m^e_{l,k}  ( 1 +  y^e_k (t-\tau) )   \   (\lambda+\nu)e^{- (\lambda+\nu) \tau } d\tau  .\label{Eqn_FP}
\end{eqnarray} 
By change of variable from $t-\tau = s$ we obtain:
\begin{eqnarray}
y^e_{l} (t)  &=&  e^{- (\lambda+\nu) t }  \theta \int_0^t    \sum_{k \in {\cal E}_x   } a_{l, k}  y^e_k (s)   (\lambda+\nu) e^{ (\lambda+\nu) s  } ds  \nonumber  \\
&& 
+ e^{- (\lambda+\nu) t }  (1-\theta)  \int_{0}^t   \sum_{k  \in  {\cal E}_x }  m^e_{l,k}  ( 1 +  y^e_k (s) )   \   (\lambda+\nu) e^{ (\lambda+\nu) s } d s  . \label{Eqn_FP_2}
\end{eqnarray} 
Differentiating we obtain (with $\lambda_\nu := \lambda + \nu$):
\begin{eqnarray*}
\boxed{
\frac{d y_l^e (t) }{dt }  = - \lambda_\nu y^e_l (t)  + \lambda_\nu     \sum_{k \in {\cal E}_x   } \left ( \theta  a_{l, k}  + (1- \theta )  m^e_{l,k}  \right )   y^e_k (t) + \lambda_\nu (1-\theta)   \sum_{k  \in  {\cal E}_x }  m^e_{l,k} .}
\end{eqnarray*}
In vector form,  \vspace{-6mm}
\begin{eqnarray*}
\frac{d {\bf y} ^e (t) }{dt } =  A_{ex}  {\bf y} (t) +\lambda_\nu   (1-\theta )    \left [  
\begin{array}{lll}
  \sum_{k \in {\cal E}_x  }  m^e_{1,k}  \\
    \sum_{k \in {\cal E}_x  }  m^e_{2,k}     \\
    \vdots \\
      \sum_{k \in {\cal E}_x  }  m^e_{N,k}  
\end{array}
\right ] .
\end{eqnarray*}
Thus  the solution is given by:
\begin{eqnarray} 
\boxed{
 {\bf y}^e (t)   =  \left (  e^{ A_{ex} t }  - I  \right )  \lambda_\nu  (1-\theta ) A_{ex}^{-1}    \left [  
\begin{array}{lll}
  \sum_{k \in {\cal E}_x  }  m^e_{1,k}  \\
    \sum_{k \in {\cal E}_x  }  m^e_{2,k}     \\
    \vdots \\
      \sum_{k \in {\cal E}_x  }  m^e_{N,k}  
\end{array}
\right ]  .}  \label{Eqn_ye_single}
\end{eqnarray}
\ignore{
For the special case, we have:
\begin{eqnarray}
\label{Eqn_ye_special_case}
y^e_l (t)  = ( e^{\alpha_e t}  - 1)   h_l^e \mbox{ with }
  h_l^e =  \frac{  \lambda_\nu  (1-\theta)   r_l   m \eta_1    }{ \alpha_e }   \mbox{ and }  
  \alpha_e =  m \eta_1 \sum_{l}  r_l    \rho_l  -1 .
\end{eqnarray}  \eop} 


We claim  that the solution of (\ref{Eqn_FP}) has the following approximate (approximation good as $N \to \infty$)  structure for the special case of social networks:
\begin{eqnarray}
\label{Eqn_anticipated_solution}
y_l(t) =  g^e_l + h^e_l e^{ \alpha_e t }    \ \ \ \mbox{ for } \ \ l \in {\cal E}_x \mbox{ and for all } t. 
\end{eqnarray}

Directly substituting the above representation of $\{y_l(t)\}$ in both sides of the fixed point equation, we have the following ($\lambda_\nu := \lambda+ \nu$):
\begin{eqnarray*}
g^e_l + h^e_l e^{ \alpha_e t }   &=& y_l (t) 
= \theta \int_0^t    \sum_{k \in {\cal E}_x   } a_{l, k}    \left (  g^e_k + h^e_k e^{ \alpha_e (t - \tau) }   \right )  \lambda_\nu e^{-\lambda_\nu \tau } d\tau    \\
&& 
\ \ +  (1-\theta)  \int_{0}^t   \sum_{k  \in {\cal E}_x }  m^e_{l,k}  \left ( 1 +    \left (  g^e_k + h^e_k e^{  \alpha_e (t - \tau) }   \right )   \right  )   \   \lambda_\nu e^{-\lambda_\nu \tau } d\tau  
\end{eqnarray*} 
\begin{eqnarray*}
& = &   \left [\theta \sum_{k \in {\cal E}_x   } a_{l, k}   g^e_k     + (1-\theta) \sum_{k \in  {\cal E}_x }  m^e_{l,k}  (1+g^e_k) 
\right ]  (1  - e^{-\lambda_\nu t} )
 \\
&& + \ \  e^{\alpha_e t} \left [\theta \sum_{k \in {\cal E}_x   } a_{l, k} h^e_k   + (1-\theta) \sum_{k \in   {\cal E}_x }  m^e_{l,k}  h^e_k \right ]  \left (  1  -  e^{- (\lambda_\nu +  \alpha_e) t }\right )  \frac{ \lambda_\nu}{\lambda_\nu + \alpha_e } 
.
\end{eqnarray*}
Thus, if one can find the solution (for coefficients $\{g_l^e, h_l^e \}$ and $\alpha_e$)   of the following  equations,  it is easy to verify that the waveforms given by (\ref{Eqn_anticipated_solution})  are a   fixed point solution of ${\bf G} $:
\begin{eqnarray}
g^e_l &=& 
\theta \sum_{k \in {\cal E}_x   } a_{l, k}   g^e_k      + (1-\theta) \sum_{k  \in {\cal E}_x }  m^e_{l,k}  (1+g^e_k) 
 ,   \label{Eqn_gle}  
 \\
h^e_l &=&   \left (\theta \sum_{k \in {\cal E}_x   } a_{l, k} h^e_k    + (1-\theta) \sum_{k \in  {\cal E}_x }  m^e_{l,k}  h^e_k 
\right ) \frac{ \lambda_\nu}{ \lambda_\nu  + \alpha_e  }  .    \label{Eqn_hlee}   
\end{eqnarray}
and that
\begin{eqnarray}
&&- \theta \sum_{k \in {\cal E}_x   } a_{l, k}   g^e_k     - (1-\theta) \sum_{k  \in {\cal E}_x }  m^e_{l,k}  (1+g^e_k) 
\nonumber \\ 
&&   -\left [\theta \sum_{k \in {\cal E}_x   } a_{l, k} h^e_k   + (1-\theta) \sum_{k  \in {\cal E}_x }  m^e_{l,k}  h^e_k 
\right ] \frac{ \lambda_\nu}{\alpha_e + \lambda_\nu   }  
   = 0 \mbox{ for all } l. \label{Eqn_zero_FP}
\end{eqnarray} 
The last one is required   as we must  have the coefficients of $ e^{- \lambda_\nu  t }$ term zero for all $t$.  
The above implies $g_l + h^e_l = 0$. 
Further we also require
 $y_l(0) =0$, i.e., the user with post of interest has not shared,  for which again we need  $g_l + h^e_l = 0$. 
Now using  (\ref{Eqn_gle}) and  (\ref{Eqn_hlee}), we have for any $l \in \cal{E}_x$:
\begin{eqnarray}
g^e_l + \frac{\alpha_e + \lambda_\nu }{\lambda_\nu}h^e_l & = & \theta \sum_{k \in {\cal E}_x   } a_{l, k} (g^e_k +h^e_k)   +    (1-\theta) \sum_{k  \in {\cal E}_x }  m^e_{l,k}  (1+g^e_k +h^e_k)  \nonumber
\\[0.2cm]
 & = &    (1-\theta) \sum_{k  \in {\cal E}_x }  m^e_{l,k}, \ \ \mbox{and hence using   $-g_l^e = h_l^e$}   \nonumber
\\
h_l^e & = &  \left (   (1-\theta) \sum_{k  \in  {\cal E}_x   }  m^e_{l,k}  \right ) \frac{ \lambda_\nu}{    \alpha_e  } \mbox{ for each } l.   \label{Eqn_hle}
\end{eqnarray} 
 Thus, 
\begin{eqnarray*}
\boxed{  y^e_l(t )  \  =  \   \frac{  \lambda_\nu  (1-\theta) \sum_{k  \in  {\cal E}_x   }  m^e_{l,k}   }{ \alpha_e } \left (  e^{\alpha_e t} - 1  \right )  },
\end{eqnarray*}if the $\{ h_l^e \}$ given by (\ref{Eqn_hle}) satisfies  (\ref{Eqn_zero_FP}) with an appropriate $\alpha_e.$
 Note that  \underline{\it $\alpha_1 > 0$ is a required} condition for virality;  otherwise (by positiveness of $\{m_{l,k}^e\}$)
the coefficients $\{h_l^e\}$ are negative and then the  solution $y_l^e (t)$  settles to a  positive limit (which in our OSN context represents the eventual expected number of shares before extinction) and does not explode (i.e., no virality) . 

\noindent{\bf Derivation of $\alpha_e$:} From (\ref{Eqn_zero_FP}) we have:
\begin{eqnarray}
   (1-\theta) \sum_{k  \in {\cal E}_x }  m^e_{l,k}   =  \left [\theta \sum_{k \in {\cal E}_x   } a_{l, k} h^e_k   + (1-\theta) \sum_{k  \in {\cal E}_x }  m^e_{l,k}  h^e_k 
\right ] \frac{ \alpha_e }{\alpha_e + \lambda_\nu   } ,  \label{Eqn_zero_FP_Case1}
\end{eqnarray} and this equation is to be used to derive $\alpha_e$. We consider this for the example of our Social network. 
\\{
{\bf Social network example:} Here ${\cal E}_x = \{1, 2, \cdots, N\}$,  $m^e_{l, k } = r_l c^e_k$  with $c^e_k := m \eta_1 \rho_k $  and  $a_{l,k} = 1_{k = l+1} 1_{l < N} $. 
Assume $r_{l+1}/r_l = \Delta_r$ a constant  (independent of $l$).
For this  case  we require for each $l$:
$$(1-\theta) r_l \sum_k c^e_k  (\alpha_e + \lambda_\nu ) = \theta r_{l+1} (1-\theta)  \lambda_\nu \sum_k c^e_k  + \lambda_\nu  (1-\theta)^2  r_l  \sum_k c^e_k r_k \sum_{k'} c^e_{k'}   $$
Or we need
\begin{eqnarray}
\label{Eqn_need_alpha_single}
    (\alpha_e + \lambda_\nu ) = \theta r_{l+1}/r_l    \lambda_\nu \  + (1-\theta) \lambda_\nu   \sum_k c^e_k r_k  \mbox{ for all }  l \in {\cal E}_x.   \end{eqnarray}
Thus   $\alpha_e$  equals:
\begin{eqnarray}
\alpha_e \approx \lambda_\nu   \left  ( \theta  \Delta_r - 1    + (1-\theta)    m \eta_1 \sum_k \rho_k r_k  \right ).
\label{Eqn_alppha_e}
\end{eqnarray}
  This is only approximation, which is accurate as $N \to \infty$, because the above equations are not satisfied for  $l = N$. For this case, $a_{N, k}=0$ for all $k$ in social network example, thus  one needs to satisfy the following equation 
$$ r_N  (\alpha_e + \lambda_\nu ) =    r_N  (1-\theta)  m \eta_1 \sum_k \rho_k r_k  ,  $$
which is approximately true because $r_N \approx 0$ (recall $r_N \to 0$ as $N \to \infty$). 
From equations  (\ref{Eqn_zero_FP_Case1}) and  (\ref{Eqn_hle}) it is clear that $\alpha_e$ and $\{h_l^e\}$ are the eigen value and eigen vector of matrix $A_{ex}$; further $\alpha_e$ is the Perron root (largest eigen value)  as 
  $\{h_l^e\}$  is the vector of all positive elements. In all, 
$$
y_l^e (t) \approx  \frac{  \lambda_\nu  (1-\theta)   m \eta_1 r_l   }{ \alpha_e } \left (  e^{\alpha_e t} - 1  \right )  \mbox{ for all } l ,
$$  
as $\sum_l c_k^e = m \eta_1.$
 \eop
}

\section*{Appendix C:  Expected Net Progeny,  with mixed population}

\subsection*{Uniqueness of solution of ${\bf G}$} Define the following norm on the space of waveforms $ \textbf{y}(\cdot ) = \{y_{l} (\cdot )\}_l $ 
$$
||  \textbf{y} ||_\phi := \sum_l ||  y_l  ||_\phi  \mbox{ with }   \   ||  y_l  ||_\phi \ :=\  \int_0^\infty  |  y_l (t) |  \phi e^{- \phi t} dt.
$$
and observe that 
\begin{eqnarray*}
\int_0^\infty  \left |   \int_0^t    \sum_{k \in {\cal M}_x   } a_{l, k} \left (  y_k (t-\tau)  - z_k (t-\tau) \right )  (\lambda+\nu) e^{- (\lambda+\nu) \tau } d\tau   \right |  \   \phi e^{- \phi t} dt \hspace{-90mm}\\
&\le &   \int_0^\infty  \int_\tau ^\infty    \sum_{k \in {\cal M}_x   } a_{l, k}     \left |  y_k (t-\tau)  - z_k (t-\tau)   \right |  \phi e^{- \phi t} dt  \       (\lambda+\nu) e^{- (\lambda+\nu) \tau } d\tau   \\
&\le &   \int_0^\infty  \int_0^\infty    \sum_{k \in {\cal M}_x   } a_{l, k}     \left |  y_k (s)  - z_k (s)   \right |  \phi e^{- \phi (s+\tau) } ds  \       (\lambda+\nu) e^{- (\lambda+\nu) \tau } d\tau   \\
&\le &   \sum_{k \in {\cal M}_x   } a_{l, k}    \int_0^\infty \phi e^{- \phi s}    ||  {y}_k - z_k  ||_\phi     (\lambda+\nu) e^{- (\lambda+\nu) \tau } d\tau   \\ 
&=&     \frac{\lambda+\nu }{ \lambda+\nu + \phi}   \sum_{k \in {\cal M}_x   } a_{l, k}   ||  {y}_k - z_k  ||_\phi .
\end{eqnarray*}
Working in a similar way,  we have
\begin{eqnarray*}
|| {\bf G}( {\bf y} )  - {\bf G} ({\bf z}) ||_\phi &  \le &    \theta   \frac{\lambda+\nu }{ \lambda+\nu + \phi} \sum_{k \in {\cal M}_x   }   \sum_l  a_{l, k}   ||  {y}_k - z_k  ||_\phi  \\
&&
+ (1-\theta)   \frac{\lambda+\nu }{ \lambda+\nu + \phi}     \sum_{k \in {\cal M}_x   }   \sum_l  m_{l, k}   ||  {y}_k - z_k  ||_\phi 
  \end{eqnarray*}
If one chooses a suitable $\phi_s$ such that we have the following   for some $ \zeta_s < 1$:
$$
\max_k     \left (  \sum_l  \left ( \theta  a_{l, k}  + (1-\theta) m_{l,k} \right )  \right ) \  \frac{\lambda+\nu }{ \lambda+\nu + \phi_s}    =   \zeta_s  ,    
$$then
\begin{eqnarray*}
|| {\bf G}( {\bf y} )  - {\bf G} ({\bf z}) ||_{\phi_s}
&\le&     \zeta_s    \sum_{k \in {\cal M}_x   }      ||  {y}_k - z_k  ||_{\phi_s}   =  \zeta_s  |  \textbf{y} - {\bf z} ||_{\phi_s} .
\end{eqnarray*} 
Thus ${\bf G}( \cdot)$ is a contraction mapping and hence has unique fixed point solution.  \eop

\subsection*{Solution of the fixed point equation ${\bf G}$ }
  As $y_l (t)$ satisfies the following fixed point equation  for any $l \in {\cal M}_x$,
\begin{eqnarray*}
y_{l} (t)  &=&  \theta \int_0^t    \sum_{k \in {\cal M}_x   } a_{l, k}  y_k (t-\tau)   (\lambda+\nu) e^{- (\lambda+\nu) \tau } d\tau   \\
&& + (1-\theta ) \int_0^t    \sum_{k \in {\cal M}_x  } m_{l,k}   (1+ y_k (t-\tau)  )  (\lambda+\nu) e^{- (\lambda+\nu) \tau } d\tau   \\
&&
+  (1-\theta)  \int_{0}^t   \sum_{k  \in  {\cal E}_x }  m_{l,k}  ( 1 +  y^e_k (t-\tau) )   \   (\lambda+\nu)e^{- (\lambda+\nu) \tau } d\tau  .
\end{eqnarray*}
After change of variable from $t-\tau$ to $s$ and differentiating we get :
\begin{eqnarray*}
\frac{d 
y_{l} (t)}{dt}  &=&  -  \lambda_\nu  y_l (t) + \lambda_\nu  \bigg  (    \theta     \sum_{k \in {\cal M}_x   } a_{l, k}   y_k(t)  + (1-\theta )    \sum_{k \in {\cal M}_x  }  m_{l,k} (1+ y_k (t)  )   \\
&&
+  (1-\theta)     \sum_{k  \in  {\cal E}_x }  m_{l,k}  ( 1 +  y^e_k (t) )  \bigg  )   \mbox{simplifying}  \\
\frac{d 
y_{l} (t)}{dt}  &=&    \lambda_\nu  \bigg  (       \sum_{k \in {\cal M}_x   }  \big  (  \theta a_{l, k} +  (1-\theta )  m_{l,k} - 1_{l =k} \big )  y_k(t)  + (1-\theta )    \sum_{k \in {\cal M}_x  }  m_{l,k}  \\
&&
+  (1-\theta)     \sum_{k  \in  {\cal E}_x }  m_{l,k}  ( 1 +  y^e_k (t) )  \bigg  )   .
\end{eqnarray*}
In other words in vector notation
\begin{eqnarray*}
\frac { d {\bf y} (t)  } {d t }  =  A_{mx}  {\bf y} (t)  + \lambda_\nu (1-\theta )    \left [  
\begin{array}{lll}
  \sum_{k \in {\cal M}_x  }  m_{1,k}  \\
    \sum_{k \in {\cal M}_x  }  m_{2,k}     \\
    \vdots \\
      \sum_{k \in {\cal M}_x  }  m_{N,k}  
\end{array}
\right ] +  \lambda_\nu (1-\theta )   A_{mx, ex}
 \left [  
\begin{array}{lll}
 1 +   y_1^e (t) \\
1 +   y_2^e (t)     \\
    \vdots \\
 1 +   y_N^e (t) 
\end{array}
\right ].
\end{eqnarray*}
Using the standard tools of ODEs and the solution  $\{y^e_l (.)\}$ derived in (\ref{Eqn_ye_single}) we have
\begin{eqnarray*}
{\bf y} (t) &=&    \left ( e^{A_{mx} t}   -I  \right ) {\bf c_v}_0 +  e^{A_{mx} t} \int_0^t  e^{- A_{mx} s }  A_{mx, ex}      e^{ A_{ex} s }     A_{ex}^{-1}  {\bf c_v}_1  ds \\
{\bf c_v}_0 &=&  A_{mx}^{-1} \lambda_\nu  (1-\theta )  \left(    \left [  
\begin{array}{lll}
  \sum_{k \in {\cal M}_x  }  m_{1,k}  \\
    \sum_{k \in {\cal M}_x  }  m_{2,k}     \\
    \vdots \\
      \sum_{k \in {\cal M}_x  }  m_{N,k}  
\end{array}
\right  ] +   A_{mx, ex}  \left (   {\bf 1}    -     A_{ex}^{-1}  {\bf c_v}_1    \right )  \right ) \\
{\bf c_v}_1 &=&   \lambda_\nu (1-\theta)  \left [  
\begin{array}{lll}
  \sum_{k \in {\cal E}_x  }  m^e_{1,k}  \\
    \sum_{k \in {\cal E}_x  }  m^e_{2,k}     \\
    \vdots \\
      \sum_{k \in {\cal E}_x  }  m^e_{N,k}  
\end{array}
\right ] .
\end{eqnarray*} \eop

\subsection*{
 Simplified form for Social network example: } To keep the explanations simple we do not consider\footnote{In actuality, the group ${\cal M}_x$ (in social network example)  has two sub-classes:  the CP1-post is in higher level  than CP2-post in  in one sub-class and it is  vice versa in other.  We consider the former case ($p=0$ and when we start with CP1-post at higher level), while the  
  calculations easily extend to the other case ($p=1$) and a combination of the two cases ($ 0 < p < 1$) also.   } $p$, the probability of swapping the order of the two posts in the recipient TLs. 
We consider the number of shares of CP1-post, which can be obtained by computing the net progeny.   
  Here ${\cal E}_x = \{1^e, 2^e, \cdots, N^e\}$ (the notation $^e$   is only required to differentiate between mixed and exclusive types, but we discard this notation  when things are obvious),  and
${\cal M}_x = \{1, 2, \cdots, N-1\}$.  When one starts with $(l, l+1)$- mixed type of section \ref{sec_Diff_TLs}  (i.e., with CP1-post at level $l$ and CP2-post at $l+1$ level and   since $p=0$,   we have, 
    $m_{l,k} = r_l c^m_k$  (for any $l \in {\cal M}_x$) with $c^m_k := m\eta_1 \delta \eta_2 {\bar \rho}_k$ and $c^m_k := m \eta_1 (1- \delta \eta_2)  {\bar \rho}_k$ respectively for  
 $k \in {\cal M}_x$ and $k \in  {\cal E}_x $.  We also have $m^e_{l,k} = r_l c^e_k$ with $c^e_k := m\eta_1  \rho_k$.     \underline{\it Further we assume}
  \begin{eqnarray}
  \label{Eqn_OSN}
  {\bar \rho}_l = \rho_l \mbox{  for all } l  \mbox{  and that } \Delta_r = r_{l+1} /  r_l \mbox{ when } l < N.   
  \end{eqnarray}
 Now we   have the following:
\begin{eqnarray}
\label{Eqn_SpecialCase}
\sum_{k \in {\cal M}_x}  m_{l,k}  &=&   m  \eta_1 \eta_2 \delta   r_l    \mbox{ for any }   l \in {\cal M}_x,     \nonumber  \\
\sum_{k \in {\cal E}_x}  m_{l,k}  & =&    m  \eta_1 (1- \eta_2 \delta )  r_l   \mbox{ for any }   l \in {\cal M}_x,     \nonumber  \\
\sum_{k \in {\cal E}_x}  m^e_{l,k}  & =&   m  \eta_1   r_l    \mbox{ for any }   l \in {\cal E}_x,  \mbox{ and finally, }
\sum_{k \in {\cal E}_x \cup {\cal M}_x } c_k^m  \  =   m  \eta_1    =  \  \sum_{ k \in {\cal E}_x}  c_k. \hspace{10mm}
\end{eqnarray}Here $r_l$ remains the same for mixed as well as exclusive types. 
 We claim  under this special case that the net progeny has the following approximate (approximation good as $N \to \infty$) simplified structure:
\begin{eqnarray}
\label{Eqn_anticipated_solution}
y_l(t) =  g_l + h_l e^{ \alpha_e t } + o_l e^{ {\bar \alpha} t }  \ \ \ \mbox{ for } \ \ l \in {\cal M}_x \mbox{ and for all } t,
\end{eqnarray}with ${\bar \alpha} = \alpha_{mx}$ the largest eigen value of matrix ${ A}_{mx}$ and appropriate $\{g_l, h_l, o_l\}$ and $\alpha_e$ is the Perron root/largest (positive) eigen value of $A_{ex}$. We prove our claim in the following:

Directly substituting the above representation of $\{y_l(t)\}$ in both sides of the fixed point equation, we have the following ($\lambda_\nu := \lambda+ \nu$):
\begin{eqnarray*}
g_l + h_l e^{ \alpha_e t } + o_l e^{ {\bar \alpha} t }   &=& y_l (t) 
= \theta \int_0^t    \sum_{k \in {\cal M}_x   } a_{l, k}    \left (  g_k + h_k e^{ \alpha_e (t - \tau) } + o_k e^{ {\bar \alpha} (t - \tau)  }  \right )  \lambda_\nu e^{-\lambda_\nu \tau } d\tau   \\
&& + (1-\theta ) \int_0^t    \sum_{k \in {\cal M}_x  }  m_{l,k}   \left (  1+ g_k + h_k e^{ \alpha_e (t - \tau) } + o_k e^{ {\bar \alpha} (t - \tau)  }  \right )    \lambda_\nu e^{-\lambda_\nu \tau } d\tau   \\
&& 
\ \ +  (1-\theta)  \int_{0}^t   \sum_{k  \in {\cal E}_x }  m_{l,k}  \left ( 1 +    \left (  g^e_k + h^e_k e^{  \alpha_e (t - \tau) }   \right )   \right  )   \   \lambda_\nu e^{-\lambda_\nu \tau } d\tau  
\end{eqnarray*} 
\begin{eqnarray*}
& = &   \left [\theta \sum_{k \in {\cal M}_x   } a_{l, k}   g_k  + (1-\theta) \sum_{k \in {\cal M}_x  }  m_{l,k}   (1+ g_k)   + (1-\theta) \sum_{k \in  {\cal E}_x }  m_{l,k}  (1+g^e_k) 
\right ]  (1  - e^{-\lambda_\nu t} )
 \\
&& + \ \  e^{\alpha_e t} \left [\theta \sum_{k \in {\cal M}_x   } a_{l, k} h_k + (1-\theta) \sum_{k \in {\cal M}_x  }  m_{l,k}    h_k   + (1-\theta) \sum_{k \in   {\cal E}_x }  m_{l,k}  h^e_k \right ]  \left (  1  -  e^{- (\lambda_\nu +  \alpha_e) t }\right )  \frac{ \lambda_\nu}{\lambda_\nu + \alpha_e } 
\\
&&   + \ \   e^{ {\bar \alpha} t} \left [\theta \sum_{k \in {\cal M}_x   } a_{l, k} o_k + (1-\theta) \sum_{k \in {\cal M}_x  }  m_{l,k}    o_k    
\right ]  \left (  1  -  e^{- (\lambda_\nu +  {\bar \alpha}) t }\right )  \frac{ \lambda_\nu}{\lambda_\nu + {\bar \alpha} }.
\end{eqnarray*}
Thus, if one can find the solution (for coefficients $\{g_l, h_l, o_l\}$ and ${\bar \alpha}$)  of the following  equations,  it is easy to verify that the waveforms given by (\ref{Eqn_anticipated_solution})  are a   fixed point solution of ${\bf G} $:
\begin{eqnarray}
g_l &=& 
\theta \sum_{k \in {\cal M}_x   } a_{l, k}   g_k  + (1-\theta) \sum_{k \in {\cal M}_x  }  m_{l,k}   (1+ g_k)   + (1-\theta) \sum_{k  \in {\cal E}_x }  m_{l,k}  (1+g^e_k) 
 ,   \label{Eqn_gl}  
 \\
h_l &=&   \left (\theta \sum_{k \in {\cal M}_x   } a_{l, k} h_k + (1-\theta) \sum_{k \in {\cal M}_x  }  m_{l,k}    h_k   + (1-\theta) \sum_{k \in  {\cal E}_x }  m_{l,k}  h^e_k 
\right ) \frac{ \lambda_\nu}{ \lambda_\nu  + \alpha_e  }  \mbox{ and }    \label{Eqn_hll}   
\\
o_l &=&  \left (\theta \sum_{k \in {\cal M}_x   } a_{l, k} o_k + (1-\theta) \sum_{k \in {\cal M}_x  }  m_{l,k}    o_k    
\right )  \frac{ \lambda_\nu}{\lambda_\nu  +{\bar \alpha}  }  \label{Eqn_ol}  ,
\end{eqnarray}and  we also  require
\begin{eqnarray}
g_l + h_l + o_l = 0 \mbox{ for all } l,  \label{Eqn_lamnu_zero}
\end{eqnarray}
 because we must  have the coefficients of $ e^{- \lambda_\nu  t }$ term zero for all $t$.
\ignore{
\begin{eqnarray}
&&- \theta \sum_{k \in {\cal M}_x   } a_{l, k}   g_k  - (1-\theta) \sum_{k \in {\cal M}_x  }  m_{l,k}    (1+ g_k)    - (1-\theta) \sum_{k  \in {\cal E}_x }  m_{l,k}  (1+g^e_k)   \nonumber
 \\ 
&&   -\left [\theta \sum_{k \in {\cal M}_x   } a_{l, k} h_k + (1-\theta) \sum_{k \in {\cal M}_x  }  m_{l,k}    h_k   + (1-\theta) \sum_{k  \in {\cal E}_x }  m_{l,k}  h^e_k  \nonumber 
\right ] \frac{ \lambda_\nu}{\alpha_e + \lambda_\nu   } 
\\
&& -
\left [\theta \sum_{k \in {\cal M}_x   } a_{l, k} o_k + (1-\theta) \sum_{k \in {\cal M}_x  }  m_{l,k}    o_k    
\right ]  \frac{ \lambda_\nu}{{\bar \alpha} +\lambda_\nu  } = 0 \mbox{ for all } l.  \label{Eqn_lamnu_zero}
\end{eqnarray}  }
The above condition is also required as 
 $y_l(0) =0$ (for all $l$), i.e., the user with post of interest has not shared. 
 Further computations are considered in two cases.

{ 
\textbf{Case I, without   ${\bar \alpha}$ term:}  We will first try for a solution  with $o_l = 0$.  If such a solution is possible, then by uniqueness this is the solution and in this case we will not require an additional  positive ${\bar \alpha}$. In our OSN example, such a thing is possible.  
Of course we would require some conditions as shown below. 

Further and more importantly, we are finding here approximate solution  for the example of OSNs; the solution might be unique, however one may have more than one approximation (for the same unique solution); we will see that one sub case (when started with $(l, l+1)$ types) has another approximation in Case 2 provided below.}

With $o_l = 0$, one requires    $g_l = - h_l$ for all $l$. 
 Now using  (\ref{Eqn_gl}) and  (\ref{Eqn_hll}), we have:
\begin{eqnarray*}
g_l + \frac{\alpha_e + \lambda_\nu }{\lambda_\nu}h_l & = & \theta \sum_{k \in {\cal M}_x   } a_{l, k} (g_k +h_k) + (1-\theta) \sum_{k \in {\cal M}_x  }  m_{l,k}   (1+ g_k + h_k)
\\
&& +  \ \  (1-\theta) \sum_{k  \in {\cal E}_x }  m_{l,k}  (1+g^e_k +h^e_k)
\\[0.2cm]
 & = &  (1-\theta) \sum_{k \in {\cal M}_x  }  m_{l,k} +  (1-\theta) \sum_{k  \in {\cal E}_x }  m_{l,k} .
\end{eqnarray*}
By substituting $-g_l = h_l$ and solving we get (for each $l$):
\begin{eqnarray}
h_l& = &
 \frac{  \lambda_\nu  (1-\theta) \sum_{k  \in  {\cal E}_x \cup {\cal M}_x  }  m_{l,k}   }{ \alpha_e } = - g_l. \nonumber
\\
&& \hspace{-2.5cm} \mbox{Thus, \ } y_l(t )\ =  -\frac{  \lambda_\nu  (1-\theta) \sum_{k  \in  {\cal E}_x \cup {\cal M}_x  }  m_{l,k}   }{ \alpha_e } + \frac{  \lambda_\nu  (1-\theta) \sum_{k  \in  {\cal E}_x \cup {\cal M}_x  }  m_{l,k}   }{ \alpha_e } e^{\alpha_e t}.  \label{Eqn_potential_hl}
\end{eqnarray}

Further from (\ref{Eqn_lamnu_zero}) $\alpha_e$ should satisfy (for all $l \in {\cal M}_x$) the following  as well as (\ref{Eqn_need_alpha_single}):
{\small \begin{eqnarray}
  (1-\theta) \sum_{k \in {\cal M}_x \cup   {\cal E}_x  }  m_{l,k}       =   \left [\theta \sum_{k \in {\cal M}_x   } a_{l, k} h_k + (1-\theta) \sum_{k \in {\cal M}_x  }  m_{l,k}    h_k   + (1-\theta) \sum_{k  \in {\cal E}_x }  m_{l,k}  h^e_k   
\right ] \frac{ \alpha_e }{\alpha_e + \lambda_\nu   } .  \label{Eqn_gl_hl_zero}
\end{eqnarray}  }
{Consider now  the example of OSNs (see equations (\ref{Eqn_OSN}) and (\ref{Eqn_SpecialCase}) and the description above it), i.e., for which  $\alpha_e$  has to equivalently satisfy the following for all 
$l$:
\begin{eqnarray}
\alpha_e + \lambda_\nu = \theta \Delta_r \lambda_\nu {\mathbb 1}_{ l  < N} + \lambda_\nu  (1-\theta) \sum_{k \in {\cal M}_x \cup {\cal E}_x} c_k^m  r_k . 
\label{Eqn_alpha_e_Case1}
\end{eqnarray}
So we are done (for all $l < N$) if same $\alpha_e$ satisfies the above as well as (\ref{Eqn_need_alpha_single}), i.e.,  if  the following is true:
$$
m \eta_1 \sum_k \rho_k r_k = \sum_{k \in {\cal E}_x} c_k r_k =  \sum_{k \in {\cal M}_x \cup {\cal E}_x} c_k^m  r_k = m \eta_1 \sum_k {\bar \rho}_k r_k; 
$$and this is true by conditions (\ref{Eqn_OSN}) and (\ref{Eqn_SpecialCase}). 
Once again, the   equation (\ref{Eqn_gl_hl_zero}) is satisfied only approximately for the case when $l = N$. However when $\theta = 0$, it is satisfied exactly and hence the solution is the exact and hence is the unique solution. 
}\\{
{\bf Another sub case:}  When we start  with  $(l+1, l)$ types and  consider net progeny corresponding to CP1-post, then it is easy to verify that we will have 
\begin{eqnarray}
\label{Eqn_Special_sub_Case}
\sum_{k \in {\cal M}_x}  m_{l,k}  &=&   m  \eta_1 \eta_2 \delta   r_l    \mbox{ for any }   l \in {\cal M}_x \ \ \ (\mbox{means when offspring is of   } (l+1,l) \mbox{ type} ),     \nonumber  \\
\sum_{k \in {\cal E}_x}  m_{l,k}  & =&    m  \eta_1 (1- \eta_2  ) \delta  r_l   \mbox{ for any }   l \in {\cal M}_x,    \nonumber  \\
\sum_{k \in {\cal E}_x}  m^e_{l,k}  & =&   m  \eta_1   r_l    \mbox{ for any }   l \in {\cal E}_x, . \hspace{10mm}
\end{eqnarray}
For this sub-case we will have
$$ 
    \  \sum_{ k \in {\cal E}_x}  c^e_k  =    m  \eta_1   \mbox{,  but  }  \sum_{k \in {\cal E}_x \cup {\cal M}_x } c_k^m  \  = m \eta_1 \delta.
$$
The rest of the details can go through, however equation (\ref{Eqn_gl_hl_zero}) is not satisfied. Thus we can't have this kind of a solution. For this sub-case the solution would be given by Case II, provided the value ${\bar \alpha}$ (computed below) is positive. 
}
\ignore{
Thus we need that the potential values of $\{h_l\}$ given by (\ref{Eqn_potential_hl}) satisfy the following:
\begin{eqnarray*}
h_l \alpha_e =   \lambda_\nu \sum_{k \in {\cal M}_x   }   \left ( \theta a_{l,k}  +   (1-\theta)  m_{l,k}  - 1_{l=k}   \right ) h_k   + \lambda_\nu   (1-\theta) \sum_{k  \in {\cal E}_x }  m_{l,k}  h^e_k \mbox{ for all } 
l \in {\cal M}_x.
\end{eqnarray*}
We would consider this only for special case. Recall $\{h_l^e \}$ is eigen vector of $A_{ex}$ matrix corresponding to the largest eigen value $\alpha_e$ for the special case and equals $\lambda_\nu {\bar c}^e {\bf r} /\alpha_e $ and for the special case,  $\{h_l\}$ given by (\ref{Eqn_potential_hl}) also  simplifies to   $\lambda_\nu {\bar c}^e {\bf r} /\alpha_e $ (because ${\bar c}^m + {\bar c}^{me} = {\bar c}^e$ as seen from (\ref{Eqn_SpecialCase})   and from (\ref{Eqn_ye_special_case})
$$
\alpha_e = \lambda_\nu \left (
m \eta_1 \sum_{l}  r_l    \rho_l  -1  \right )
$$

However unfortunately the RHS of the equation becomes  (after cancelling $\lambda_\nu$, $r_l$   $\lambda_\nu {\bar c}^e / \alpha_e $ etc and note $\sum_{k \in {\cal E}_x} m_{l,k} + \sum_{k \in {\cal M}_x} m_{l,k}  =   $) 
$$
\lambda_\nu \left (
m \eta_1 \sum_{l}  r_l    {\bar \rho}_l  -1  \right )
$$
and hence are not equal }

\textbf{Case II, with additional ${\bar \alpha}$ term:}  When Case I is not possible, we should try with $o_l \ne 0$  and appropriate ${\bar \alpha}$, and recall by uniqueness  if such a solution is possible, it would be the solution\footnote{Once again we are only obtaining approximate solutions and one may have more than one approximate solution even when the solution is unique.}.  When $o_l \ne 0$  using   $g_l + h_l + o_l = 0$ and  $g_l^e + h_l^e = 0$ for each $l$  and further using (\ref{Eqn_gl})-(\ref{Eqn_ol}) we get
%

\vspace{-4mm}
{\small 
\begin{eqnarray*}
g_l + \frac{ \lambda_\nu + \alpha_e}{\lambda_\nu} h_l +  \frac{ \lambda_\nu + \bar{\alpha}}{\lambda_\nu} o_l & \ \ = \ &
 \theta \sum_{k \in {\cal M}_x   } a_{l, k} \left(g_k + h_k + o_k  \right) + \sum_{k \in {\cal M}_x} (1-\theta)   m_{l,k}  \left(1+g_k + h_k + o_k \right)  \nonumber
\\
 && \hspace*{.1cm} \ +  \sum_{k \in {\cal E}_x} (1-\theta)   m_{l,k}  \left(1+g^e_k + h^e_k   \right) \ \  = \  \sum_{k  \in  {\cal E}_x \cup {\cal M}_x } (1-\theta)   m_{l,k}, \mbox{ and so, }
 \\
 h_l & =& \left [ - g_l +  (1-\theta) \sum_{k  \in  {\cal E}_x \cup {\cal M}_x }  m_{l,k}  - o_l   \frac{{\bar \alpha} +\lambda_\nu  }{ \lambda_\nu} \right ] \frac{ \lambda_\nu}{  \lambda_\nu + \alpha_e  }  \\
& =& \left [   -g_l -o_l +  (1-\theta) \sum_{k  \in  {\cal E}_x \cup {\cal M}_x }  m_{l,k}  - o_l   \frac{{\bar \alpha}   }{ \lambda_\nu} \right ] \frac{ \lambda_\nu}{  \lambda_\nu + \alpha_e  }. 
\end{eqnarray*} }
Thus

\vspace{-4mm}
{\small
\begin{flalign}
\label{Eqn_hl_gl_in_ol}
h_l =  \frac{  \lambda_\nu  (1-\theta) \sum_{k  \in  {\cal E}_x \cup {\cal M}_x  }  m_{l,k} -   o_l    {\bar \alpha}}  { \alpha_e }, \ \mbox{ and } \ 
g_l = -h_l - o_l = -  \frac{  \lambda_\nu  (1-\theta) \sum_{k  \in  {\cal E}_x \cup {\cal M}_x  }  m_{l,k} -   o_l  (   {\bar \alpha}  -  \alpha_e  ) }  { \alpha_e }. 
\end{flalign}}

Now summing equations (\ref{Eqn_gl})-(\ref{Eqn_ol}),  and  by using $h_l+ g_l + o_l = 0$  and $h_l^e = -g_l^e $ in equation (\ref{Eqn_hll}),  we have:

\vspace{-4mm}
{\small\begin{eqnarray*}
 0 = g_l + h_l + o_l  &=& 
 \left [\theta \sum_{k \in {\cal M}_x   } a_{l, k}   g_k  + (1-\theta) \sum_{k \in {\cal M}_x  }  m_{l,k}   (1+ g_k)   + (1-\theta) \sum_{k \in  {\cal E}_x }  m_{l,k}  (1+g^e_k) 
\right ]  
 \\
 && +    \left [\theta \sum_{k \in {\cal M}_x   } a_{l, k} h_k + (1-\theta) \sum_{k \in {\cal M}_x  }  m_{l,k}    h_k   + (1-\theta) \sum_{k  \in {\cal E}_x }  m_{l,k}  h^e_k 
\right ]  \left ( 1 -  \frac{  \alpha_e }{ \lambda_\nu  + \alpha_e  } \right )  
  \\
&& +   \left [\theta \sum_{k \in {\cal M}_x   } a_{l, k} o_k + (1-\theta) \sum_{k \in {\cal M}_x  }  m_{l,k}    o_k    
\right ] \left ( 1-   \frac{ { \bar \alpha}  }{\lambda_\nu  +{\bar \alpha}  }  \right )
 \\
\\ 
 &=&   - \left [\theta \sum_{k \in {\cal M}_x   } a_{l, k} h_k + (1-\theta) \sum_{k \in {\cal M}_x  }  m_{l,k}    h_k   + (1-\theta) \sum_{k  \in {\cal E}_x }  m_{l,k}  h^e_k 
\right ] \frac{ \alpha_e  }{ \lambda_\nu  + \alpha_e  }   \\
 &&
-  \left [\theta \sum_{k \in {\cal M}_x   } a_{l, k} o_k + (1-\theta) \sum_{k \in {\cal M}_x  }  m_{l,k}    o_k    
\right ]  \frac{  {\bar \alpha} }{\lambda_\nu  +{\bar \alpha}  }  
+    (1-\theta) \sum_{k \in {\cal M}_x \cup {\cal E}_x  }  m_{l,k}.      
\end{eqnarray*}}
Substituting  the value of $h_k$ as given in equation (\ref{Eqn_hl_gl_in_ol}), we get:

\vspace{-4mm}
{\small{
\begin{eqnarray*}
&& \hspace{-1cm} - \sum_{k \in {\cal M}_x   } \left(\theta  a_{l, k}  +  (1-\theta)  m_{l,k}\right) \left(  \frac{  \lambda_\nu  (1-\theta) \sum_{k'  \in  {\cal E}_x \cup {\cal M}_x  }  m_{k,k'} -   o_k    {\bar \alpha}}  { \alpha_e } \right)    \frac{ \alpha_e  }{ \lambda_\nu  + \alpha_e  } - (1-\theta) \sum_{k  \in   {\cal E}_x }  m_{l,k}  h^e_k \frac{ \alpha_e  }{ \lambda_\nu  + \alpha_e  }
\\
 && \hspace*{.51cm}  -  \left [\theta \sum_{k \in {\cal M}_x   } a_{l, k} o_k + (1-\theta) \sum_{k \in {\cal M}_x  }  m_{l,k}    o_k    
\right ]  \frac{  {\bar \alpha} }{\lambda_\nu  +{\bar \alpha}  } 
+   (1-\theta) \sum_{k \in {\cal M}_x \cup {\cal E}_x  }  m_{l,k}   =0.
\end{eqnarray*}}}
Simplifying,
{\footnotesize{
\begin{eqnarray*}
\left [\theta \sum_{k \in {\cal M}_x   } a_{l, k} o_k + (1-\theta) \sum_{k \in {\cal M}_x  }  m_{l,k}    o_k    
\right ]   \left ( \frac{  {\bar \alpha} }{\lambda_\nu  + \alpha_e  }  -  \frac{  {\bar \alpha} }{\lambda_\nu  +{\bar \alpha}  }  \right )  
&=& (1-\theta) \sum_{k  \in  {\cal E}_x }  m_{l,k}  h^e_k 
  \frac{ \alpha_e  }{ \lambda_\nu  + \alpha_e  }   -   (1-\theta) \sum_{k \in {\cal M}_x \cup {\cal E}_x  }  m_{l,k}
\\
   && \hspace*{-3cm}+ \ \ \left [ \sum_{k \in {\cal M}_x   } \left (  \theta  a_{l, k} + (1-\theta)   m_{l,k}  \right )  \sum_{k'  \in  {\cal E}_x \cup {\cal M}_x }  (1-\theta)   m_{k, k'}      
\right ] \frac{ \lambda_\nu  }{ \lambda_\nu  + \alpha_e  }.
\end{eqnarray*}}
Finally using equation (\ref{Eqn_ol}),
\begin{eqnarray*}
  o_l \frac{\lambda_\nu  +{\bar \alpha} }{\lambda_\nu  }\left ( \frac{  {\bar \alpha} }{\lambda_\nu  + \alpha_e  }  -  \frac{  {\bar \alpha} }{\lambda_\nu  +{\bar \alpha}  }  \right )  
&=& ( 1-\theta) \sum_{k  \in   {\cal E}_x }  m_{l,k}  h^e_k 
  \frac{ \alpha_e  }{ \lambda_\nu  + \alpha_e  }  -   (1-\theta) \sum_{k \in {\cal M}_x \cup {\cal E}_x  }  m_{l,k} 
  \\
 && \hspace*{.01cm} + \ \ (1-\theta) \left [ \sum_{k \in {\cal M}_x   } \left (  \theta  a_{l, k} + (1-\theta)   m_{l,k}  \right )  \sum_{k'  \in  {\cal E}_x \cup {\cal M}_x }   m_{k, k'}     
\right ] \frac{ \lambda_\nu  }{ \lambda_\nu  + \alpha_e  } \ \ \ \mbox{}.
\end{eqnarray*}}
Thus for any $l \in {\cal M}_x$,
\begin{eqnarray*}o_l  &=&  
    \frac{ (1-\theta) \lambda_\nu \alpha_e     \sum_{k  \in {\cal E}_x }  m_{l,k}  h^e_k   }{  ( {\bar \alpha}-\alpha_e) {\bar \alpha} }
     + \frac{  (1-\theta) \lambda_\nu  } { ( {\bar \alpha}-\alpha_e) {\bar \alpha}} \Bigg (  \lambda_\nu       \sum_{k \in {\cal M}_x   } \left (  \theta  a_{l, k} + (1-\theta)   m_{l,k}  \right )  \sum_{k'  \in  {\cal E}_x \cup {\cal M}_x }   m_{k, k'}   \nonumber
     \\
     & & \hspace{1cm} - \ \  ({ \lambda_\nu  + \alpha_e  }) \sum_{k \in {\cal M}_x \cup {\cal E}_x  }  m_{l,k} 
   \Bigg ) .\hspace{10mm} \label{Eqn_expr_ol}\end{eqnarray*}
Similarly, one can obtain the closed form expressions for $h_l$ and $g_l$ by substituting the value of $o_l$ in equation (\ref{Eqn_hl_gl_in_ol}). 
$$
\mbox{Recall and we will use, }  \  \  h_l^e  =   \left (   (1-\theta) \sum_{k  \in  {\cal E}_x   }  m^e_{l,k}  \right ) \frac{ \lambda_\nu}{    \alpha_e  } \approx    \left (   (1-\theta)  m \eta_1 r_l  \right ) \frac{ \lambda_\nu}{    \alpha_e  }  \mbox{ for each } l \in {\cal E}_x,   
$$for Social network example.
{ From (\ref{Eqn_expr_ol}), for Social network example (\ref{Eqn_SpecialCase}), $o_l$ simplifies as below
\begin{eqnarray*}
o_l &=&  r_l    \frac{ (1-\theta) \lambda_\nu     }{  ( {\bar \alpha}-\alpha_e) {\bar \alpha} }  \left (     \lambda_\nu \Delta_r \theta \mathbb{1}_{ l < N} + (1-\theta)   \lambda_\nu  m^2 \eta_1^2    \sum_{k   } {\bar \rho}_k  r_k  - (\lambda_\nu +\alpha_e)  m \eta_1   
 \right ) \\ 
 &\approx & r_l  \frac{ (1-\theta) \lambda_\nu     }{  ( {\bar \alpha}-\alpha_e) {\bar \alpha} }     \lambda_\nu \Delta_r \theta  (1 - m \eta_1) \mbox{ for all } l, \ \ \    \    \mbox{ when }  r_N \approx 0
\end{eqnarray*} and these  or equivalently vector $\{r_l \}$ should satisfy   equation (\ref{Eqn_ol}) with appropriate ${\bar \alpha}$, i.e., we will require
\begin{eqnarray*}
r_l &=&  \left (\theta \sum_{k \in {\cal M}_x   } a_{l, k} r_k + (1-\theta) \sum_{k \in {\cal M}_x  }  m_{l,k}    r_k    
\right )  \frac{ \lambda_\nu}{\lambda_\nu  +{\bar \alpha} }  ,
\end{eqnarray*} which for social network example (special case) translates to satisfy the following:
\begin{eqnarray*}
r_l &=&  \left (\theta \Delta_r  r_l  + (1-\theta) r_l  \sum_{k \in {\cal M}_x  } c_k^m     r_k    
\right )  \frac{ \lambda_\nu}{\lambda_\nu  +{\bar \alpha} }  .
\end{eqnarray*} Thus  the following ${\bar \alpha}$ satisfies all equations
\begin{eqnarray*}
{\bar \alpha }  =  \lambda_\nu   \left ( \theta \Delta_r - 1   + (1-\theta) m\eta_1 \delta \eta_2 \sum_k {\bar \rho}_k r_k \right ),
\end{eqnarray*}except for the case with $l = N$, for which it is an approximation  as in Case 1.  Again as seen above $\{o_l\}$ (or equivalently $\{r_l\}$) is a vector of all positive or all negative entries, thus ${\bar \alpha}$ 
would be the unique Perron root of  ${A}_{mx}$ (see equation (\ref{Eqn_ol})), when it is positive definite (i.e., when ${\bar \alpha} > 0$).
This is the second approximate solution for OSNs, when  one starts with  $(l, l+1)$ types. It is easy to verify that the two approximations coincide when $\theta = 0$; easy to verify that $o_l = 0$ for all $l$ and $h_l$ with case 2 coincides with $h_l$ of case 1. In this case the solution in fact is exact, as there is no difference between $l < N$ and $l = N$  in  (\ref{Eqn_gl_hl_zero}). For this case, 
$$
h_l =  \frac{  \lambda_\nu  (1-\theta) m \eta_1 r_l  -   o_l    {\bar \alpha}}  { \alpha_e } \mbox{ for any }  l \in {\cal M}_x.
$$

 It is easy to check that the \underline{sub-case (\ref{Eqn_Special_sub_Case})  (starting with $(l+1,l$ types)} mentioned at the end of Case I can also satisfy the conditions of this case and has the  solution with same ${\bar \alpha}$, but with    
\begin{eqnarray*}
o_l  &=& r_l    \frac{ (1-\theta) \lambda_\nu     }{  ( {\bar \alpha}-\alpha_e) {\bar \alpha} }  \left (     \lambda_\nu \Delta_r \theta + (1-\theta)   \lambda_\nu  m^2 \eta_1^2 \delta (1-\eta_2 + \eta_2  \delta  )     \sum_{k   } {\bar \rho}_k  r_k  - (\lambda_\nu +\alpha_e)  m \eta_1    \delta
 \right ) \mbox{ and }  \\
 h_l &=&  \frac{  \lambda_\nu  (1-\theta) m \eta_1 \delta  r_l  -   o_l    {\bar \alpha}}  { \alpha_e } \mbox{ for any }  l \in {\cal M}_x. \hspace{10mm}  \mbox{\eop}
\end{eqnarray*}

 }

\section*{Appendix D: Expected number of shares in non-viral scenario}
We have  $$y_{l,k}^j = E[\lim_{t \to \infty }Y^j (t) | {\bf X}(0)  = {\bf e}_{l,k}].$$
 $y_{l,k}^j $ can  be obtained by solving appropriate FP equations.  These FP equations are obtained by conditioning on the events of the first transition, as before. Here we have some additional events depending upon the starting TL. When a mixed type TL is subjected to the `share transition',  then we can have shares exclusively of the post of one of the CPs, and or shares of both the posts.   Whereas when an exclusive CP-type  TL  is subjected to a `share transition',    only exclusive types are engendered, as in single CP.  With the `shift' transition, we have similar changes as in the single CP case.

Below we obtain the expected shares for CP-1 without loss of generality and hence suppress the superscript $^j$ for remaining discussions.

Let $Y_{l,k} = \lim_{t \to \infty} Y_{l,k} (t)$ be the total  number of shares of post-${\bf P}$, before extinction, when started with one TL of $(l,k)$ type  (with $k = l+1$ or $l-1$). Let 
$y_{l,k} := E\left[Y_{l,k} \right]$ be its expected value.   The total  number of shares  of any CP post is finite on the extinction paths. Thus,  by conditioning on the events of first transition epoch,  one can write the following recursive equations for any $l < N$: 
\begin{eqnarray}
\label{Eqn_yll1}
y_{l,l+1} & = & \theta \Big ( \mathbb{1}_{\{l < N-1\}} y_{l+1,l+2} + \mathbb{1}_{\{l = N-1\}} y_{N,0}  \Big )  +(1-\theta)r_l (1-\delta) m\eta_1 (1 + \textbf{y}_{ex1}.\bm{\bar{\rho} })  \nonumber \\
&& + \ (1-\theta)r_l\delta  m\eta_1\Bigg [(1-\eta_2)(1+ \textbf{y}_{ex1}.\bm{\bar{\rho}})  \nonumber  
+\eta_2(1+ p\textbf{y}_{mx1}.\bm{\bar{\rho}} +  (1-p)\textbf{y}_{mx2}.\bm{\bar{\rho}})\Bigg ] \nonumber 
\end{eqnarray}
where ${\bf  y}_{ex1} = \{ y_{1,0}^1, y_{2, 0}^1, \cdots, y_{N-1,0}\} $. And again for any $l < N$, 
{\footnotesize{
\begin{eqnarray}
\label{Eqn_yl1l}
y_{l+1,l} & = & \mathbb{1}_{\{l < N-1\}}\theta y_{l+2,l+1} 
 + \ (1-\theta)r_l\delta  m\eta_1\Bigg [(1-\eta_2)(1+ \textbf{y}_{ex1}.\bm{\bar{\rho}})  \nonumber 
 \   \eta_2 \big (1+ p\textbf{y}_{mx1}.\bm{\bar{\rho}} +  (1-p)\textbf{y}_{mx2}.\bm{\bar{\rho}} \big ) \Bigg ].\nonumber 
\end{eqnarray}}}

One can easily solve the above set of linear equations to obtain the  fixed point solution, by first obtaining the solutions for 
\begin{eqnarray*}
{\bf y}_{mx1} \idot {\bm {\bar \rho}}+{\bf y}_{mx2}  \idot {\bm {\bar \rho}}\mbox{ with } \hspace{-10mm}\\
 {\bf y}_{mx1 } & := &\{ y_{1,2}^1, y_{2, 3}^1, \cdots, y^1_{N-1,N}\}  \mbox{ and }  \\
{\bf y}_{mx2 } &:= &\{ y_{2,1}^1, y_{3,2}^1, \cdots, y^1_{N, N-1}\}  . \end{eqnarray*}
 We carry out the analysis for the special case with reading probabilities: $r_i = d_1 d_2^i$. 
Recall $c_j =  (1-\theta) m\eta_j$. Define the following which will be used only in this part:
\begin{eqnarray}
\label{Eqn_consts}
B_{ex1, \delta} & = & c_1\delta(1-\eta_2)\left(1+ \textbf{y}_{ex1}.\bm{\bar{\rho}} \right), \nonumber
\\
B_{ex1, 1-\delta}& = & c_1 (1-\delta) (1 + \textbf{y}_{ex1}.\bm{\bar{\rho} }), \nonumber
\\
\bar{C}_{mx1} & = & c_{mx}\left(1+ \left(p \textbf{y}_{mx1} + (1-p)\textbf{y}_{mx2}\right).\bm{\bar{\rho}} \right) \mbox{ and } \nonumber
\\
\bar{C}_{mx2} & = & c_{mx}\left(1+ \left((1-p) \textbf{y}_{mx1} + p\textbf{y}_{mx2}\right).\bm{\bar{\rho}} \right)\end{eqnarray}
The first two quantities can be computed from the expected number of  shares given as in Part-I \cite{Ranbir2} (in the single CP model), while the remaining are obtained by solving the above FP equations. Now we can rewrite the   equations  (\ref{Eqn_yll1})-(\ref{Eqn_yl1l}) in the following manner for the special case\footnote{One can easily write down the equations for general case, but are avoid to simplify the expressions.} with $r_i = d_1 d_2^i$

\begin{eqnarray*}
y_{l,l+1} & = & \theta y_{l+1,l+2}  +  \left( B_{ex1, 1-\delta} + \bar{C}_{mx2} + B_{ex1, \delta} \right) d_1 d_2^l   \\
&&   \hspace{35mm}\mbox{ for  }  l < N-1  \mbox{ and }
\\ 
 y_{N-1,N} & = & \theta y_{N,0} + \left( B_{ex1, 1-\delta} + \bar{C}_{mx2} + B_{ex1, \delta} \right) d_1 d_2^{N-1} .
\end{eqnarray*}
Solving these equations using backward recursion: 

 \begin{eqnarray*}
y_{N-2,N-1}  =   \theta^2 y_{N,0} +  \left( B_{ex1, 1-\delta} + \bar{C}_{mx2} + B_{ex1, \delta} \right)  \left(\theta  d_1 d_2^{N-1} + d_1 d_2^{N-2}\right).
 \end{eqnarray*}
 and then continuing in a similar way
 \begin{eqnarray*}
y_{N-l,N-l+1}  =  \theta^l y_{N,0} +  \left( B_{ex1, 1-\delta} + \bar{C}_{mx2} + B_{ex1, \delta} \right) d_1 d_2^{N-l} \left[ \sum_{i = 0}^{l-1} \left(\theta d_2 \right)^{i} \right]. 
 \end{eqnarray*} 
 One can rewrite it as the following for any $l < N$:

  \begin{eqnarray}
y_{l,l+1}  & =& \theta^{N-l} y_{N,0} +  \left( B_{ex1, 1-\delta} + \bar{C}_{mx2} + B_{ex1, \delta} \right) d_1 d_2^l \left[ \sum_{i = 0}^{N-l-1} \left(\theta d_2 \right)^{i} \right] 
\nonumber  \\
 & =& \theta^{N-l} y_{N,0} +  \left( B_{ex1, 1-\delta} + \bar{C}_{mx2} + B_{ex1, \delta} \right) d_1 d_2^l \frac{1-\left(\theta d_2 \right)^{N-l}}{1-\theta d_2}. \hspace{6mm}  
  \label{Eqn_simpl1}
\end{eqnarray}
In exactly similar lines, for any $l < N$, we have:
\begin{eqnarray*}
y_{l+1,l}  =  \theta y_{l+2,l+1}  + \left( \bar{C}_{mx1} + B_{ex1, \delta} \right) d_1 d_2^l.
 \end{eqnarray*} This simplifies to the following for any $l < N$:
  \begin{eqnarray}
  \label{Eqn_simpl2}
y_{l+1,l}  =  \theta^{N-l}  +   \left( \bar{C}_{mx1} + B_{ex1, \delta} \right) d_1 d_2^l \frac{1-\left(\theta d_2 \right)^{N-l}}{1-\theta d_2}.
\end{eqnarray}
Multiplying the left hand sides of the equations (\ref{Eqn_simpl1}) and (\ref{Eqn_simpl2}) with ${\bar \rho}_l$ and summing it up we obtain 
$\textbf{y}_{mx1}.\bm{\bar{\rho}}$ and $\textbf{y}_{mx2}.\bm{\bar{\rho}}$ respectively:
\begin{eqnarray}
\textbf{y}_{mx1}.\bm{\bar{\rho}} & = &    \sum_{l < N} {\bar \rho}_l  \theta^{N-l}  y_{N, 0}  +  \left( B_{ex1, 1-\delta} + \bar{C}_{mx2} + B_{ex1, \delta} \right)  O_{mx} \hspace{6mm}  \label{TwoCPShares12}
\\
\textbf{y}_{mx2}.\bm{\bar{\rho}} & = &\sum_{l < N} {\bar \rho}_l  \theta^{N-l}    +  \left(  \bar{C}_{mx1} + B_{ex1, \delta} \right)  O_{mx}  \\
O_{mx} &:=& d_1 \sum_l \frac{\big (d_2^l {\bar \rho}_l \big) -  (\theta d_2)^N  \Big ( \bar{\rho}_l / \theta^l \Big ) }{\left(1-d_2 \bar{\rho}\right)\left(1-\theta d_2\right)} . \label{TwoCPShares21}
\end{eqnarray}
Note that for general $r_l$ which need not be $d_1 d_2^l$, we will have
\begin{eqnarray}
O_{mx} := \sum_{l < N}  {\bar \rho}_l  \sum_{i=0}^{N-l-1}  \theta^i r_{l+i} . \label{Eqn_Omx}
\end{eqnarray}

On adding equations (\ref{TwoCPShares12}) and (\ref{TwoCPShares21})
\begin{eqnarray*}
\textbf{y}_{mx1}.\bm{\bar{\rho}} + \textbf{y}_{mx2}.\bm{\bar{\rho}} & = & 
\sum_{l < N} {\bar \rho}_l  \theta^{N-l}
(1 +  y_{N,0})    +
\left( B_{ex1, 1-\delta} + \bar{C}_{mx1} + \bar{C}_{mx2} + 2 B_{ex1, \delta} \right) O_{mx} .
\end{eqnarray*}
This implies   using (\ref{Eqn_consts})

\begin{eqnarray*}
\textbf{y}_{mx1}.\bm{\bar{\rho}} + \textbf{y}_{mx2}.\bm{\bar{\rho}} & = &
\sum_{l < N} {\bar \rho}_l  \theta^{N-l}
(1 +  y_{N,0})   +
 \left( B_{ex1, 1-\delta} +c_{mx} \big( 2 + \textbf{y}_{mx1}.\bm{\bar{\rho}} + \textbf{y}_{mx2}.\bm{\bar{\rho}} \big)+ 2 B_{ex1, \delta} \right) O_{mx} .
\end{eqnarray*}
We have unique fixed point solution (when  $c_{mx} O_{mx} < 1$) for  $\textbf{y}^1_{mx} \idot\bm{\bar{\rho}}  := \textbf{y}_{mx1}.\bm{\bar{\rho}} + \textbf{y}_{mx2}.\bm{\bar{\rho}}$, which equals 

\begin{eqnarray*}
\textbf{y}^1_{mx} \idot\bm{\bar{\rho}}     =   \frac{\sum_{l < N} {\bar \rho}_l  \theta^{N-l}
(1 +  y_{N,0})   + \left(2\left( B_{exj, \delta} + c_{mx} \right)+ B_{exj, 1-\delta}\right)O_{m_x}}{1-c_{mx}O_{m_x}}.
\end{eqnarray*}
In the above $y_{N,0}$ and ${\bf y}_{ex1}$ of equation (\ref{Eqn_consts}) can be obtained as in Part-I\cite{Ranbir2} 

We   obtain further simpler expressions   for the special case, when  $\bar{\rho}_l = \tilde{\bar{\rho}} \bar{\rho}^l$ (with ${\bar \rho} < 1$) with  
$$\tilde{\bar{\rho}} = \frac{1}{\sum_{i=1}^{N-1} \bar{\rho}^i } =\frac{ (1- {\bar \rho})} { {\bar \rho}(1-{\bar \rho}^{N-1})} $$ and when $N \to \infty$.   Observe that
\begin{eqnarray*}
O_{mx} = d_1 
\tilde{\bar{\rho}} 
\sum_{i=1}^{N-1} \frac{\left(d_2 \bar{\rho} \right)^i- \left(d_2\theta\right)^{N} \left(\frac{\bar{\rho}}{\theta}\right)^i}{1-\theta d_2} \, \to \, \frac{d_1 d_2  (1-{\bar \rho}) }{\left(1-d_2 \bar{\rho}\right)\left(1-\theta d_2\right)} \ ,
\end{eqnarray*}as $N \to \infty$ because 
$$
\tilde{\bar{\rho}} 
\sum_{i=0}^{N-1}  \left(d_2\theta \right)^{N} \left(\frac{\bar{\rho}}{\theta}\right)^i   
=  \left(d_2\theta \right)^{N} \frac{(1-{\bar \rho}^N)}{(1- {\bar \rho})}  \frac{ \theta^N - {\bar \rho}^N}{ \theta - \rho}  \theta^{-N+1 } \to 0.
$$
In a similar way
  $$
\tilde{\bar{\rho}}  \theta^N  \sum_{l < N} ( {\bar \rho}/  \theta)^{l}  = 
 \tilde{\bar{\rho}}  \theta   \frac{ \theta^N - {\bar \rho}^N} {\theta - {\bar \rho} } \to 0.
$$
 And $y_{N,0}$ can be bounded as $N  \to \infty$ (see Part-I \cite{Ranbir2} for details).
Thus as   $N \to \infty$ for any  $j=1,2$:
 
\begin{eqnarray*}
\textbf{y}^j_{mx} \idot\bm{\bar{\rho}} 
& \to &
 \frac{\Big ( 2 c_j\delta \left[ (1 + \textbf{y}^j_{exj}.\bm{\bar{\rho} })(1-\eta_{-j}) + \eta_{-j}\right]  +  c_j (1-\delta) (1 + \textbf{y}^j_{exj}.\bm{\bar{\rho} }) \Big )O_{mx}}{1-c_{mx}  O_{mx}} 
 \end{eqnarray*}
 \begin{eqnarray*}
\mbox{with  }  O_{m_x} &\to & \frac{d_1  d_2  (1-\bar{\rho})}{\left(1-d_2 \bar{\rho}\right)\left(1-\theta d_2\right)} \mbox{ and where } -j := 2 \mathbb{1}_{\{j=1\}} + 1 \mathbb{1}_{\{j=2\}}.
\end{eqnarray*}

Here $ \textbf{y}^j_{exj}$  is similar to that in Part-I \cite{Ranbir2}, and  $\{ y_{l,k}^j \}$ with $k = l+1$ or $l-1$ 
can be computed  uniquely using $\{\textbf{y}^j_{mx} \idot \bm{\bar{\rho}} \}$.   \eop

 \bibliographystyle{apacite}

\end{document}